\documentclass[aop,MSNbibl,seceqn,dvips]{arximspdf}
\usepackage{mathbh}
\usepackage{graphicx}

\doi{10.1214/14-AOP908}
\volume{42}
\issue{5}
\pubyear{2014}
\firstpage{1809}
\lastpage{1884}

\makeatletter

\newcommand{\rrVert}{\Vert}
\newcommand{\rrvert}{\vert}
\newcommand{\llVert}{\Vert}
\newcommand{\llvert}{\vert}

\newcommand{\dx}[1]{\textup{d} #1}
\newcommand{\bigcupdot}{\charfusion[\mathop]{\bigcup}{\cdot}}
\newcommand{\eps}{\varepsilon}

\def\moverlay{\mathpalette\mov@rlay}
\def\mov@rlay#1#2{\leavevmode\vtop{%
\baselineskip\z@skip \lineskiplimit-\maxdimen
\ialign{\hfil$\m@th#1##$\hfil\cr #2\crcr}}}
\newcommand{\charfusion}[3][\mathord]{
#1{\ifx #1\mathop\vphantom{#2}\fi
\mathpalette\mov@rlay{#2\cr #3}
}
\ifx #1\mathop\expandafter\displaylimits\fi}

\newcommand{\mathds}{\mathbb}
\newcommand{\one}{\mathbh{1}}
\newcommand{\dsOne}{\mathbh{1}}
\newcommand{\N}{\mathds{N}}
\newcommand{\R}{\mathds{R}}

\newcommand{\cA}{\mathcal{A}}
\newcommand{\cE}{\mathcal{E}}
\newcommand{\cH}{\mathcal{H}}
\newcommand{\cM}{\mathcal{M}}
\newcommand{\cN}{\mathcal{N}}
\newcommand{\cP}{\mathcal{P}}
\newcommand{\cS}{\mathcal{S}}
\newcommand{\cT}{\mathcal{T}}

\newcommand{\tts}{\mathtt{s}}
\newcommand{\ttPI}{\mathtt{PI}}
\newcommand{\ttLSI}{\mathtt{LSI}}

\newcommand{\Expect}{\mathds{E}}
\newcommand{\laplace}{\Delta}

\newcommand{\argmin}{\operatorname{\arg\min}}
\newcommand{\cov}{\operatorname{cov}}
\newcommand{\dist}{\operatorname{dist}}
\newcommand{\Ent}{\operatorname{Ent}}
\newcommand{\HWI}{\operatorname{HWI}}
\newcommand{\Id}{\operatorname{Id}}
\newcommand{\LSI}{\operatorname{LSI}}
\newcommand{\dLSI}{\operatorname{dLSI}}
\newcommand{\osc}{\operatorname{osc}}
\newcommand{\PI}{\operatorname{PI}}
\newcommand{\supp}{\operatorname{supp}}
\newcommand{\sym}{\operatorname{sym}}
\newcommand{\var}{\operatorname{var}}
\newcommand{\vspan}{\operatorname{span}}
\newcommand{\WI}{\operatorname{WI}}

\newcommand{\eqref}[1]{(\ref{#1})}
\newcommand{\overset}{\stackrel}
\newcommand{\mathclap}[1]{\mbox{\fontsize{8.36}{10.36}\selectfont{#1}}}

\newtheorem{theo}{Theorem}[section]
\newtheorem{cor}[theo]{Corollary}
\newtheorem{lem}[theo]{Lemma}
\newtheorem{prop}[theo]{Proposition}
\newproclaim{assume}[theo]{Assumption}
\newproclaim{defn}[theo]{Definition}
\newproclaim{rem}[theo]{Remark}

\renewcommand\thelonglist{\roman{longlist}}
\renewcommand\labellonglist{(\thelonglist)}

\newcommand{\fraca}[2]{{#1}/{#2}}
\newcommand{\fracb}[2]{{(#1)}/{#2}}
\newcommand{\fracc}[2]{{#1}/{(#2)}}
\newcommand{\fracf}[2]{({#1})/({#2})}
\makeatother

\begin{document}
\begin{frontmatter}

\title{Poincar\'{e} and logarithmic Sobolev inequalities by
decomposition of the energy landscape}
\runtitle{$\PI$ and $\LSI$ by decomposition of the energy landscape}

\begin{aug}
\author[A]{\fnms{Georg} \snm{Menz}\ead[label=e1]{gmenz@stanford.edu}}
\and
\author[B]{\fnms{Andr\'{e}} \snm{Schlichting}\corref{}\ead[label=e2]{schlichting@iam.uni-bonn.de}\thanksref{t1}}
\runauthor{G. Menz and A. Schlichting}
\thankstext{t1}{Supported by DFG project FOR718, the IMPRS of the
Max-Planck-Institute and the University of Leipzig.}
\affiliation{Stanford University and University of Bonn}
\address[A]{Department of Mathematics\\
Stanford University\\
Building 380\\
Stanford, California 94305\\
USA\\ \printead{e1}} 
\address[B]{Institute for Applied Mathematics\\
University of Bonn \\
Endenicher Allee 60\\
D-53115 Bonn\\
Germany\\ \printead{e2}}
\end{aug}

\received{\smonth{7} \syear{2012}}
\revised{\smonth{12} \syear{2013}}

%
\begin{abstract}
We consider a diffusion on a potential landscape which is given by a
smooth Hamiltonian $H\dvtx \R^n\to\R$ in the regime of low
temperature $\eps$. We proof the Eyring--Kramers formula for the
optimal constant in the Poincar\'{e} (PI) and logarithmic Sobolev
inequality (LSI) for the associated generator $L= \eps\laplace-
\nabla H \cdot\nabla$ of the diffusion. The proof is based on a
refinement of the two-scale approach introduced by Grunewald et al.
 [\textit{Ann. Inst. Henri Poincar\'e Probab. Stat.}
\textbf{45} (2009)
302--351] and of the mean-difference estimate introduced by Chafa\"
{i} and Malrieu [\textit{Ann. Inst. Henri Poincar\'e Probab. Stat.}
\textbf{46} (2010)
 72--96].
The Eyring--Kramers formula follows as a simple corollary from two main
ingredients:
The first one shows that the PI and LSI constant of the diffusion
restricted to metastable regions corresponding to the local minima
scales well in $\eps$. This mimics the fast convergence of the
diffusion to metastable states.
The second ingredient is the estimation of a mean-difference by a
weighted transport distance. It contains the main contribution to the
PI and LSI constant, resulting from exponentially long waiting times of
jumps between metastable states of the diffusion.
\end{abstract}

%

\pdfkeywords{60J60, 35P15, 49R05, Diffusion process, Eyring-Kramers formula, Kramers law, metastability, Poincare inequality, spectral gap, logarithmic Sobolev inequality, weighted transport distance}

\end{frontmatter}

\section{Introduction}\label{sec1}

Let us consider a diffusion on a potential landscape which is given by
a sufficiently smooth \emph{Hamiltonian function} $H\dvtx \R^n\to\R$. We
are interested in the regime of low temperature $\eps>0$. The \emph
{generator} of the diffusion has the following form:
\begin{equation}
\label{edefL} L := \eps\laplace- \nabla H \cdot\nabla.
\end{equation}
The associated \emph{Dirichlet form} is given for a test function
$f\in H^1(\mu)$ by
\[
\label{edefDirichletform} \cE(f) := \int(-Lf)f \,\dx{\mu} = \eps\int\llvert \nabla f\rrvert
^2 \,\dx {\mu} .
\]
The corresponding diffusion $\xi_t$ satisfies the stochastic
differential equation
\begin{equation}
\label{eoverLang} \dx\xi_t = - \nabla H (\xi_t) \,\dx t +
\sqrt{2 \varepsilon} \,\dx B_t,
\end{equation}
where $B_t$ is the \emph{Brownian motion} on $\R^n$. Equation \eqref
{eoverLang} is also called \emph{over-damped Langevin equation}
(cf., e.g., \cite{Legoll2010a}). Under some growth assumption on $H$,
there exists an equilibrium measure of the according stochastic
process, which is called \emph{Gibbs measure} and is given by
\begin{equation}
\label{dgibbsmeasure}\quad  \mu(\dx{x}) = \frac{1}{Z_\mu}\exp \biggl(-\frac{H(x)}{\eps
}
\biggr)\,\dx{x} \qquad\mbox{with } Z_\mu= \int\exp \biggl(-
\frac {H(x)}{\eps} \biggr)\,\dx{x} .
\end{equation}
The evolution \eqref{eoverLang} of the stochastic process $\xi_t$
can be translated into an evolution of the density of the process $\xi
_t$. Namely, under the assumption that the law of the initial state
$\xi_0$ is absolutely continuous w.r.t. the Gibbs measure $\mu$, the
density $f_t \mu$ of the process $\xi_t$ satisfies the \emph
{Fokker--Planck equation} (cf., e.g., \cite{Oksendal2003} or \cite
{Schuss2010})
\[
\label{eevodensity} \partial_t f_t = L f_t = \eps
\laplace f_t - \nabla H \cdot\nabla f_t .
\]
We are particularly interested in the case where $H$ has several local
minima. Then for small $\eps$, the process shows metastable behavior
in the sense that there exists a separation of scales: On the fast
scale, the process converges quickly to a neighborhood of a local
minimum. On the slow scale, the process stays nearby a local minimum
for an exponentially long waiting time after which it eventually jumps
to another local minimum.

This behavior was first described in the context of chemical reactions.
The exponential waiting time follows the \emph{Arrhenius' law} \cite
{Arr} meaning that the mean exit time from one local minimum of $H$ to
another one is exponentially large in the energy barrier between them.
By now, the Arrhenius law is well understood even for nonreversible
systems by the \emph{Freidlin--Wentzell theory \cite{FrWe}}, which is
based on \emph{large deviations}.

\begin{figure}

\includegraphics{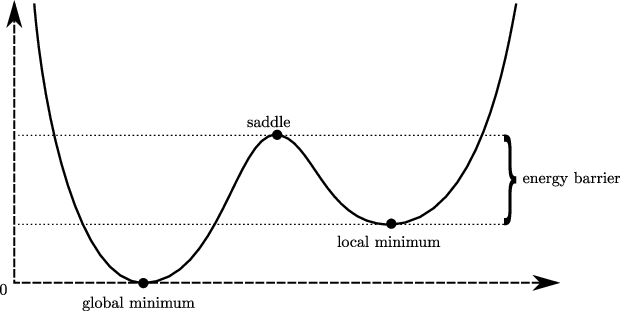}

\caption{General double-well potential $H$ on $\R$.} \label
{figdouble-well-simple}
\end{figure}

A refinement of the Arrhenius law is the \emph{Eyring--Kramers formula}
which additionally considers \emph{pre-exponential} factors. The
Eyring--Kramers formula for the Poincar\'{e} inequality (PI) goes back
to Eyring \cite{Eyring1935} and Kramers \cite{Kramers1940}. Both
argue that also in high-dimensional problems of chemical reactions most
reactions are nearby a single trajectory called \emph{reaction
pathway}. Evaluating the Hamiltonian along this \emph{reaction
coordinate} gives the classical picture of a double well potential
(cf. Figure~\ref{figdouble-well-simple}) in one dimension with an
\emph{energy barrier} separating the two local minima for which
explicit calculations are feasible.
%

However, a rigorous proof of the Eyring--Kramers formula for the
multidimensional case was open for a long time. For a special case,
where all the minima of the potential as well as all the lowest saddle
points in-between have the same energy, Sugiura \cite{Sugiura1995}
defined an exponentially rescaled Markov chain on the set of minima in
such a way that the preexponential factors become the transitions rates
between the metastable regions of the rescaled process. For the generic
case, where the local minima and saddles have different energies, the
group of Bovier et al. \cite{Bovier2004,Bovier2005} obtained
first-order asymptotics that are sharp in the parameter $\eps$. They
also clarified the close connection between \emph{mean exit times},
\emph{capacities} and the exponentially small eigenvalues of the
operator $L$ given by \eqref{edefL}. The main tool of \cite
{Bovier2004,Bovier2005} is potential theory. The small eigenvalues are
related to the mean exit times of appropriate subsets of the
state space. Further, the mean exit times are given by \emph{Newtonian
capacities} which can explicitly be calculated in the regime of low
temperature $\varepsilon$.

Shortly after, Helffer, Klein and Nier \cite
{Helffer2004,Helffer2006,Helffer2005} also deduced the Eyring--Kramers
formula using the connection of the spectral gap estimate of the
Fokker--Planck operator $L$ given by \eqref{edefL} to the one of the
\emph{Witten Laplacian}. This approach makes it possible to get
quantitative results with the help of semiclassical analysis. They
deduced sharp asymptotics of the exponentially small eigenvalues of $L$
and gave an explicit expansion in $\eps$ to theoretically any order.
An overview on the Eyring--Kramers formula can be found in the review
article of Berglund \cite{Berglund2011}.

In this work, we provide a new proof of the Eyring--Kramers formula for
the first eigenvalue of the operator $L$, that is, its \emph{spectral
gap}. The advantage of this new approach is that it extends to the
logarithmic Sobolev inequality (LSI), which was not investigated
before. The LSI was introduced by \cite{Gross} and is stronger than the PI.
Therefore, the LSI is usually harder to deduce than the PI due to its
nonlinear structure.

By deducing the Eyring--Kramers formula for the LSI, we encounter a
surprising effect: In the generic situation of having two local minima
with different energies, the Eyring--Kramers formula for the LSI
differs from the Eyring--Kramers formula for the PI by a term of
inverse order in $\varepsilon$. However, in the symmetric situation of
having local minima with the same energy, the Eyring--Kramers formula
for the LSI coincides with the corresponding formula for the PI
(cf. Corollary~\ref{coreyringkramerscomp}).

We conclude the \hyperref[sec1]{Introduction} with an overview of the article:\vspace*{0.5\baselineskip}

\noindent\hangindent=0pt In Section~\ref{sdefPILSI}, we introduce PI and LSI.

\noindent\hangindent=0pt
In Section~\ref{ssetting}, we discuss the setting and the assumptions
on the Hamiltonian~$H$.

\noindent\hangindent=0pt
In Section~\ref{soutlinemainresults}, we outline the new approach
and state the main results of this work.

\noindent\hangindent=0pt
In Section~\ref{clocal} and Section~\ref{CMEANDIFF}, we proof the
main ingredients of our new approach. Namely, in Section~\ref{clocal}, we deduce a local PI and a local LSI with optimal scaling in
$\varepsilon$, whereas in Section~\ref{CMEANDIFF} we estimate a
mean-difference by using a weighted transport distance.

\noindent\hangindent=0pt
In the \hyperref[app]{Appendices}, we provide for the convenience of the reader some
basic but nonstandard facts that are used in our arguments.

\subsection{Poincar\'{e} and logarithmic Sobolev inequality}\label
{sdefPILSI}

%
\begin{defn}[{[\ref{PIvarrho} and \ref{LSIalpha}]}]\label{defsglsi}
Let $X$ be an Euclidean space. A~Borel probability measure $\mu$ on
$X$ satisfies the \emph{Poincar\'{e} inequality} with constant
$\varrho> 0$, if for all test functions $f\in H^1(\mu)$
{\renewcommand{\theequation}{$\PI(\varrho)$}
\begin{equation}\label{PIvarrho}
\var_\mu(f) := \int \biggl(f-\int f \,\dx{\mu} \biggr)^2
\,\dx{\mu} \leq\frac{1}{\varrho} \int\llvert \nabla f\rrvert ^2 \,
\dx {\mu} .
\end{equation}
}%
In a similar way, the probability measure $\mu$ satisfies the \emph
{logarithmic Sobolev inequality} with constant $\alpha>0$, if for all
test function $f\dvtx X \to\R^+$ with $I(f\mu|\mu)<\infty$ holds
{\renewcommand{\theequation}{$\LSI(\alpha)$}
\begin{equation}\label{LSIalpha}
\Ent_\mu(f) := \int f \log\frac{f}{\int f  \,\dx{\mu} } \,\dx {\mu} \leq
\frac{1}{\alpha} \int\frac{\llvert \nabla f\rrvert ^2}{2 f} \,\dx{\mu} =: I(f\mu|\mu) ,\hspace*{-20pt}
\end{equation}
\setcounter{equation}{3}}%
where $I(f\mu|\mu)$ is called \emph{Fisher information}. The
gradient $\nabla$ is determined by the Euclidean structure of $X$.
Test functions are those functions for which the gradient exists and
the right-hand side in \ref{PIvarrho} and \ref{LSIalpha} is well defined.
\end{defn}
%
%
\begin{rem}[{[Relation between \ref{PIvarrho} and \ref{LSIalpha}]}]\label{remLSIPI}
Rothaus \cite{Rothaus1978} observed that \ref{LSIalpha} implies $\PI
(\alpha)$. This can be seen by setting $f=1+\eta g$ for $\eta$ small
and observing that
\[
\Ent_\mu\bigl(f^2\bigr) = 2\eta^2
\var_{\mu}(g) + O\bigl(\eta^3\bigr) \quad\mbox{as well as}
\quad\int\llvert \nabla f\rrvert ^2 \,\dx{\mu} = \eta^2
\int\llvert \nabla g\rrvert ^2 \,\dx{\mu} .
\]
Hence, if $\mu$ satisfies \ref{LSIalpha} then $\mu$ satisfies $\PI
(\alpha)$, which always implies $\alpha\leq\varrho$.
\end{rem}
%
%

\subsection{Setting and assumptions}\label{ssetting}
This article uses almost the same setting as found in \cite
{Bovier2004,Bovier2005}. Before stating the precise assumptions on the
Hamiltonian $H$, we introduce the notion of a \emph{Morse function}.
%
\begin{defn}[(Morse function)]\label{defmorse}
A smooth function $H\dvtx \R^n\to\R$ is a \emph{Morse function}, if the
Hessian $\nabla^2 H$ of $H$ is nondegenerated on the set of critical
points. More precisely, for some $1 \leq C_H < \infty$ holds
\begin{equation}
\label{emorse} \forall x \in\cS:= \bigl\{x\in\R^n\dvtx  \nabla H = 0 \bigr
\} \dvtx  \frac
{\llvert \xi\rrvert ^2}{C_H} \leq \bigl\langle\xi , \nabla^2 H(x)\xi \bigr
\rangle \leq C_H \llvert \xi\rrvert ^2 .
\end{equation}
\end{defn}
We make the following growth assumption on the Hamiltonian $H$
sufficient to ensure the existence of PI and LSI.
Hereby, we have to assume stronger properties for $H$ if we want to
proof the LSI.
%
\begin{assume}[($\PI$)]\label{assumeenv}
$H\in C^3(\R^n, \R)$ is a nonnegative Morse function, such that for
some constants $C_H>0$ and $K_H \geq0$ holds
{\setcounter{equation}{0}\renewcommand{\theequation}{{A\arabic{equation}}$_\ttPI$}
\begin{eqnarray}
\label{assumegradsuperlinear} \liminf_{\llvert  x\rrvert  \to\infty} \llvert \nabla H\rrvert &\geq&
C_H ,
\\
\label{assumegradlaplace} \liminf_{\llvert  x\rrvert  \to\infty} \bigl(\llvert \nabla H\rrvert
^2 - \laplace H \bigr) &\geq& - K_H .
\end{eqnarray}
}%
\end{assume}
%
%
\begin{assume}[($\LSI$)]\label{assumeenvLSI}
$H\in C^3(\R^n, \R)$ is a nonnegative Morse function, such that for
some constants $C_H>0$ and $K_H \geq0$ holds
{\setcounter{equation}{0}\renewcommand{\theequation}{{A\arabic{equation}}$_\ttLSI$}
\begin{eqnarray}
\label{assumelyapLSI} \liminf_{\llvert  x\rrvert  \to\infty}\frac{\llvert
\nabla H(x)\rrvert ^2 - \laplace
H(x)}{\llvert  x\rrvert ^2} &\geq&
C_H ,
\\
\label{assumehessLSI} \inf_{x} \nabla^2 H(x) &\geq& -
K_H .
\end{eqnarray}
}%
\end{assume}
%
%
\begin{rem}[(Discussion of assumptions)]\label{assumeremdiscussion}
The Assumption~\ref{assumeenv} yields the following consequences for
the Hamiltonian $H$:
\begin{itemize}
\item
The condition \eqref{assumegradsuperlinear} and $H(x) \geq0$
ensures that $e^{-H}$ is integrable and can be normalized to a
probability measure on $\R^n$ (see Lemma~\ref{plineargrowth}).
Hence, the Gibbs measure $\mu$ given by \eqref{dgibbsmeasure} is
well defined.
%
\item The Morse Assumption \eqref{emorse} together with the growth
condition \eqref{assumegradsuperlinear} ensures that the set $\cS$
of critical points is discrete and finite. In particular, it follows
that the set of local minima $\mathcal{M} =  \{m_1,\dots
,m_M \}$ is
also finite, that is, $M := \# \mathcal M < \infty$.
\item The \emph{Lyapunov-type} condition \eqref{assumegradlaplace}
allows to recover the Poincar\'{e} constant of the full Gibbs measure
$\mu$ from the Poincar\'{e} constant of the Gibbs measure $\mu_U$
restricted to some bounded domain $U$ (cf. Section~\ref{clocal}).
Because Gibbs measures with finite support and smooth Hamiltonian
always satisfy a Poincar\'{e} inequality with some unspecified
constant, we get that the Gibbs measure $\mu$ also satisfies a
Poincar\'{e} inequality. Equivalently, this means that there exists a
spectral gap for the operator $L$ given by \eqref{edefL}.%
\end{itemize}
Similarly the Assumption~\ref{assumeenvLSI} has the following
consequences for the Hamiltonian $H$:
\begin{itemize}
\item
One difference between the Assumptions~\ref{assumeenv} and~\ref
{assumeenvLSI} is that \eqref{assumegradsuperlinear} yields linear
growth at infinity for $H$, whereas a combination of condition \eqref
{assumelyapLSI} and \eqref{assumehessLSI} yields quadratic growth;
that is,
{\renewcommand{\theequation}{{A0}$_\ttLSI$}
\begin{equation}
\label{assumegradLSI} \liminf_{\llvert  x\rrvert  \to\infty} \frac{\llvert
\nabla H(x)\rrvert }{\llvert  x\rrvert } \geq
C_H .
\end{equation}
}%
Note that quadratic growth at infinity is a necessary condition to
obtain \ref{LSIalpha} with $\alpha>0$ (cf. \cite{Royer2007}, Theorem~3.1.21).
\item
In addition, \eqref{assumelyapLSI} and \eqref{assumehessLSI}
imply \eqref{assumegradsuperlinear} and \eqref{assumegradlaplace},
which is only an indication that \ref{LSIalpha} is stronger than \ref{PIvarrho} in the sense of Remark~\ref{remLSIPI}.
%
\item
The condition \eqref{assumelyapLSI} is again a \emph{Lyapunov type}
condition.
To enforce it to a LSI, additionally the condition \eqref
{assumehessLSI} has to be assumed (cf. Section~\ref{clocal}).
\end{itemize}
\end{rem}
To keep the presentation feasible and clear, we additionally assume a
nondegeneracy assumption, even if it is not really needed for the proof
of the Eyring--Kramers formula. The saddle height $\widehat H(m_i,m_j)$
between two local minima $m_i,m_j$ is defined by
\[
\widehat H(m_i,m_j) := \inf \Bigl\{\max
_{s\in[0,1]} H\bigl(\gamma(s)\bigr) \dvtx  \gamma\in C\bigl([0,1],
\R^n\bigr), \gamma(0)=m_i,\gamma(1)=m_j
\Bigr\} .
\]

%
\begin{assume}[(Nondegeneracy)] \label{assumenondegenerate}
There exists $\delta>0$ such that:
\begin{longlist}[(iii)]
\item The saddle height between two local minima $m_i, m_j$ is attained
at a unique critical point $s_{i,j} \in\cS$ of index one, that is, it
holds $H(s_{i,j})= \widehat{H}(m_i,m_j)$ and if $ \{\lambda
_1,\dots,\lambda_n \}$ denote the eigenvalues of $\nabla^2
H(s_{i,j})$, then it holds $\lambda_1<0$ and $\lambda_i > 0$ for
$i=2,\dots, n$. The point $s_{i,j}$ is called \emph{communicating
saddle} between the minima $m_i$ and $m_j$.
\item\label{assumorderingofM} The set of local minima $\cM=
\{ m_1, \ldots, m_M  \}$ is ordered such that $m_1$ is a global
minimum and for all $ i \in \{3,\dots,M \}$ yields
\[
\label{assumeordermin} H(s_{1,2}) - H(m_2) \geq
H(s_{1,i}) - H(m_i) + \delta.
\]
\end{longlist}
\end{assume}
%
%
\begin{rem}
The fact, that $s_{i,j}$ is indeed a critical point is explained
in \cite{Jost2008}, Proposition~6.2.1.
Since $H$ is a Morse function after Assumption~\ref{assumeenv} the
critical point $s_{i,j}$ is nondegenerate. Moreover, an indirect
perturbation argument implies that $s_{i,j}$ is a saddle point of index
one, which shows that except for uniqueness, Assumption~\ref
{assumenondegenerate}(i) is already implied by Assumption~\ref
{assumeenv}. This fact is known as Murrell--Laidler theorem in the
chemical literature \cite{Wales2003}.
\end{rem}

\section{Outline of the new approach and main results}\label
{soutlinemainresults}

In this section, we present the new approach to the Eyring--Kramers
formula and formulate the main results of this article. Because the
strategy is the same for the PI and LSI, we consider both cases
simultaneously. The approach uses ideas of the two-scale approach for
LSI \cite{GORV,Otto2007a,Lelievre2009} and the method by \cite
{Chafai2010} to deduce PI and LSI estimates for mixtures of measures.
However, the heuristics outlined in the \hyperref[sec1]{Introduction} provide a good
orientation for our proceeding. Remember that we have a splitting into
two time-scales:
\begin{itemize}
\item the fast scale describes the \emph{fast relaxation} to a local
minima of $H$ and
\item the slow scale describes the \emph{exponentially long
transitions} between local equilibrium states.
\end{itemize}
Motivated by these two time scales, we specify in Section~\ref{spartition} a splitting of the measure $\mu$ into \emph{local
measures} living on a \emph{metastable regions}\vadjust{\goodbreak} around the local
minima of $H$. This splitting is lifted from the level of the measure
to the level of the variance and entropy. In this way, we obtain \emph
{local variances and entropies}, which heuristically should correspond
to the fast relaxation, and \emph{coarse-grained variances and
entropies}, which should correspond to the exponentially long transitions.

Now, we handle each contribution separately. The local variances and
entropies are estimated by local PI (cf. Theorem~\ref{thmlocalPI})
and local LSI, respectively (cf. Theorem~\ref{thmlocalLSI}). The
heuristics suggest that this contribution should be of higher order
because this step only relies on the fast scale.

Before we estimate the coarse-grained variances and entropies, we bring
them in the form of mean-differences. This is automatically the case
for the variances. However, for the coarse-grained entropies one has to
apply a new weighted discrete LSI (cf. Section~\ref{sresdiscLSI}),
which causes the difference between the PI and LSI in the
Eyring--Kramers formula.
The main contribution to the Eyring--Kramers formula
(cf. Corollary~\ref{coreyringkramersPI} and Corollary~\ref
{coreyringkramersLSI}) results from the estimation of the
mean-difference, which is stated in Theorem~\ref{thmcoarsemd}.

At this point, let us shortly summarize the main results of this article:
\begin{itemize}
\item We provide good estimates for the local variances and entropies
(cf. Section~\ref{smainlocal}) and
\item We provide sharp estimates for the mean-differences (cf.
Section~\ref{smainmd}).
\item From these main ingredients, the Eyring--Kramers formulas follow
as simple corollaries (cf. Section~\ref{sEyringKramers}).
\end{itemize}
We close this chapter with a discussion on the optimality of the
Eyring--Kramers formula for the LSI in one dimension (cf. Section~\ref{smainopt}).

\textit{Notational remark}: Almost all of the following definitions
and quantities will depend on $\eps$, for lucidity this dependence is
not expressed in the notation. The arguments and main results hold for
$\eps>0$ fixed and small.

\subsection{Partition of the state space}\label{spartition}

The inspiration to partition the state space comes from the work \cite
{Jerrum2004a} for discrete Markov chains. In order to get sharp
results, the partition of the state space $\mathbb{R}^n$ cannot be
arbitrarily but has to satisfy certain conditions.
%
\begin{defn}[(Admissible partition)]\label{defpartition}
The family $\cP_{\cM} =  \{\Omega_i \}_{i=1}^M$ with
$\Omega_i$
open and connected is called an admissible partition for $\mu$ if the
following conditions hold:
\begin{longlist}[(iii)]
\item For each local minimum $m_i \in\cM$ exists $\Omega_i\in\cP
_{\cM}$ with $m_i \in\Omega_i$ for $i=1,\dots,M$.
\item$ \{\Omega_i \}_{i=1}^M$ is a partition of $\R^n$ up
to sets
of Lebesgue measure zero, which is denoted by $\R^n = \biguplus
_{i=1}^M \Omega_i$.
\item The partition sum of each element $\Omega_i$ of $\cP_\cM$ is
approximately Gaussian, that is, for $i =1,\dots, M$
\begin{equation}
\label{ePartSumAdmPart} \quad \mu(\Omega_i) Z_\mu= \frac{(2\pi\eps)^{\fraca{n}{2}}}{\sqrt{\det
\nabla^2 H(m_i)}}
\exp \biggl(-\frac{H(m_i)}{\eps} \biggr) \bigl(1 + O\bigl(\sqrt{\eps} \llvert \log
\eps\rrvert ^{\fraca {3} {2}}\bigr) \bigr) .
\end{equation}
\end{longlist}
\end{defn}
%
%
\begin{rem}\label{rcanonicalpartition}
A canonical way to obtain an admissible partition for $\mu$ would be
to associate to every local minimum $m_i \in\cM$ for $i=1,\dots,M$
its \emph{basin of attraction $\Omega_i$ w.r.t. $H$} defined by
\[
\label{ebasinofattraction} \Omega_i := \Bigl\{y \in\R^n \dvtx  \lim
_{t\to\infty} y_t = m_i, \dot
{y}_t=-\nabla H(y_t), y_0 =y \Bigr\} .
\]
Unfortunately, this choice would lead to technical difficulties later
on. We get rid of these technical problems by choosing the partition
$\Omega_i$ in a slightly different way. For details, we refer the
reader to Section~\ref{clocal}.
\end{rem}
Using an admissible partition of the state space, one can decompose the
Gibbs measure $\mu$ into a mixture of \emph{local} Gibbs measures
$\mu_i$.
%
\begin{defn}[(Mixture representation of $\mu$)]\label{defmixturemu}
Let $\cP_\cM= \{\Omega_i \}_{i=1}^M$ be an admissible
partition for\vadjust{\goodbreak}
$\mu$. The \emph{local Gibbs measures $\mu_i$} are defined as the
restriction of $\mu$ to $\Omega_i$
\begin{equation}
\label{elocalmu} \mu_i(\dx{x}) := \frac{1}{Z_i Z_\mu}
\one_{\Omega_i}(x) \exp \biggl( - \frac{H(x)}{\varepsilon} \biggr) \,\dx{x}  \qquad
\mbox {where } Z_i := \mu(\Omega_i) .
\end{equation}
The \emph{marginal measure $\bar\mu$} is given by a sum of Dirac measures
\[
\label{emarginalmu} \bar\mu:= Z_1 \delta_1 + \cdots+
Z_M \delta_M .
\]
%
Then the \emph{mixture representation} of $\mu$ w.r.t. $\cP_\cM$
has the form
\begin{equation}
\label{emixturemu} \mu:= Z_1 \mu_1 + \cdots+ Z_M
\mu_M .
\end{equation}
\end{defn}
As was shown in \cite{Chafai2010}, Section~4.1, the decomposition of
$\mu$ yields a decomposition of the variance $\var_\mu(f)$ and
entropy $\Ent_\mu(f)$.
%
\begin{lem}[(Splitting of variance and entropy for partition)]\label
{lemsplitvar-ent}
For a mixture representation \eqref{emixturemu} of $\mu$ holds for
all $f\dvtx  \R^n \to\R$
\begin{eqnarray}
\var_{\mu}(f) &=& \sum_{i=1}^M
Z_i \var_{\mu_i}(f) + \sum_{i=1}^M
\sum_{j>i} Z_i Z_j \bigl(
\Expect_{\mu_i}(f) - \Expect_{\mu
_j}(f) \bigr)^2,
\label{elemsplitvar}
\\
\Ent_{\mu}(f) &=& \sum_{i=1}^M
Z_i\Ent_{\mu_i}(f) + \Ent_{\bar
\mu} (\bar f )
.\label{elemsplitent}
\end{eqnarray}
We call the terms $\var_{\mu_i}(f)$ and $\Ent_{\mu_i}(f)$ \emph
{local variance} and \emph{local entropy}. The term $ (\Expect
_{\mu_i}(f) - \Expect_{\mu_j}(f) )^2$ is called \emph
{mean-difference}.\vadjust{\goodbreak} The term $\Ent_{\bar\mu}(\bar f)$ is called \emph
{coarse-grained entropy} and is given by
\begin{equation}
\label{elemsplitcoarse-entropy} \Ent_{\bar\mu} (\bar f ) := \sum
_{i=1}^M Z_i \bar f_i
\log\frac{ \bar f_i }{\sum_{j=1}^M Z_j \bar f_j },
\end{equation}
where $\bar f_i := \Expect_{\mu_i}(f)$.
\end{lem}

We skip the proof of Lemma~\ref{lemsplitvar-ent} because it only
consists of a straightforward substitution of the mixture
representation \eqref{emixturemu}.
The formula \eqref{elemsplitvar} for estimating the variance $\var
_{\mu}(f)$ is already in its final form. For the relative
entropy $\Ent_{\mu}(f)$, we still have to do some work. The aim is to
get an estimate that only involves the local terms like $\var_{\mu
}(f)$ and $\Ent_{\mu_i}(f)$ and a mean difference $ (\Expect
_{\mu_i}(f) - \Expect_{\mu_j}(f) )^2$. This is achieved in the next
subsection [cf. Corollary~\ref{corsplitentropy} and \eqref
{splitequentropy}].

\subsection{Discrete logarithmic Sobolev type inequalities}\label
{sresdiscLSI}

Starting with the identity \eqref{elemsplitent}, we have to
estimate the \emph{coarse-grained entropy $\Ent_{\bar\mu}(\bar f)$}.
We expect that the main contribution comes\vadjust{\goodbreak} from this term. If $H$ has
only two minima, we can use the following discrete LSI for a Bernoulli
random variable, which was given by Higuchi and Yoshida \cite
{Higuchi1995} and Diaconis and Saloff-Coste \cite{Diaconis1996}, Theorem~A.2, at the same time.
%
\begin{lem}[(Optimal logarithmic Sobolev inequality for Bernoulli
measures)]\label{thmlsi2point}
A Bernoulli measure $\mu_p$ on $X= \{0,1 \}$, that is, a
mixture of
two Dirac measures $\mu_p = p \delta_0 + q \delta_1$ with $p+q=1$
satisfies the discrete logarithmic Sobolev inequality
\begin{equation}
\label{equlsi2point} \Ent_{\mu_p}\bigl(f^2\bigr) \leq
\frac{pq}{\Lambda(p,q)} \bigl(f(0)-f(1)\bigr)^2
\end{equation}
with optimal constant given by the \emph{logarithmic mean}
(cf. Appendix~\ref{clogmean})
\[
\Lambda(p,q) := \frac{p-q}{\log p -\log q}  \qquad \mbox{for $p\ne q$} \quad \mbox{and}
\quad \Lambda(p,p) := \lim_{q\to p} \Lambda (p,q)=p .
\]
\end{lem}
We want to handle the general case with more than two minima.
Therefore, we want to generalize Lemma~\ref{thmlsi2point} to
discrete measures with a state space with more than two elements. An
application of the modified LSI for finite Markov chains of Diaconis
and Saloff-Coste \cite{Diaconis1996}, Theorem~A.1, would not lead to
an optimal results (cf. \cite{AndrePhDthesis}, Section~2.3). Even for
a generic Markov chain on the 3-point space, the optimal logarithmic
Sobolev constant is unknown. In this work, we use the following direct
generalization of Lemma~\ref{thmlsi2point}.
%
\begin{lem}[(Weighted logarithmic Sobolev inequality)]\label{lemlogmeandiff}
For $m\in\N$ let $\mu_m = \sum_{i=1}^m Z_i \delta_i$ be a discrete
probability measure and assume that $\min_i Z_i>0$. Then for a
function $f\dvtx   \{1,\dots,m \}\to\R_0^+$ holds the weighted\vadjust{\goodbreak}
logarithmic Sobolev inequality
\begin{equation}
\label{elogmeandiff} \Ent_{\mu_m}\bigl(f^2\bigr) \leq\sum
_{i=1}^{m-1} \sum_{j=i+1}^m
\frac{Z_i
Z_j}{\Lambda(Z_i,Z_j)} (f_i-f_j )^2 .
\end{equation}
\end{lem}
\begin{pf}
We conclude by induction and find that for $m=2$ the estimate \eqref
{elogmeandiff} just becomes \eqref{equlsi2point}, which shows the
base case. For the inductive step, let us assume that \eqref
{elogmeandiff} holds for $m \geq2$. Then the entropy $\Ent_{\mu
_{m+1}}(f^2)$ can be rewritten as follows:
\[
\Ent_{\mu_{m+1}}\bigl(f^2\bigr) 
= (1-Z_{m+1})
\Ent_{\tilde\mu_m}\bigl(f^2\bigr) + \Ent_{\nu}(\tilde f) ,
\]
where the probability measure $\tilde\mu_m$ lives on $ \{
1,\dots ,m \}$ and is given by
\[
\tilde\mu_m := \sum_{i=1}^m
\frac{Z_i}{1-Z_{m+1}} \delta_i .
\]
Further, $\nu$ is the Bernoulli measure given by $\nu:=
(1-Z_{m+1})\delta_0 + Z_{m+1} \delta_1$ and the function $\tilde f \dvtx
\{0,1 \}\to\R$ is given with values
\[
\tilde f_0 := \sum_{i=1}^m
\frac{Z_i f_i^2}{1-Z_{m+1}} \quad\mbox {and}\quad\tilde f_1 :=
f_{m+1}^2 .
\]
Now, we apply the inductive hypothesis to $\Ent_{\tilde\mu_m}(f^2)$
and arrive at
\begin{eqnarray}
(1-Z_{m+1})\Ent_{\tilde\mu_m}\bigl(f^2\bigr) & \leq&
(1-Z_{m+1}) \sum_{i=1}^m \sum
_{j>i} \frac{Z_i Z_j}{(1-Z_{m+1})^2} \frac
{1-Z_{m+1}}{\Lambda(Z_i, Z_j)}
(f_i - f_j )^2
\nonumber
\\
&=& \sum_{i=1}^m \sum
_{j>i} \frac{ Z_i Z_j}{\Lambda(Z_i, Z_j)} (f_i - f_j
)^2 ,\nonumber
\end{eqnarray}
where we used $\Lambda(\cdot,\cdot)$ being homogeneous of degree one
in both arguments (cf. Appendix~\ref{clogmean}), that is, $\Lambda
(\lambda a, \lambda b) = \lambda\Lambda(a,b)$ for $\lambda,a,b>0$.
We can apply the inductive base to the second entropy $\Ent_{\nu
}(\tilde f)$, which is nothing else but the discrete LSI for the
two-point space \eqref{equlsi2point}
\begin{equation}
\label{elogmeandiffp1} \Ent_{\nu}(\tilde f) \leq\frac{Z_{m+1} (1-Z_{m+1})}{\Lambda
(Z_{m+1}, 1 - Z_{m+1})} \bigl(\sqrt{
\tilde f_0} - \sqrt{\tilde f_1} \bigr)^2 .
\end{equation}
The last step is to apply the Jensen inequality to recover the square
differences $(f_i-f_{m+1})^2$ from
\begin{eqnarray}
\bigl(\sqrt{\tilde f_0} - \sqrt{\tilde f_1} \bigr)^2
&=& \sum_{i=1}^m \frac{Z_i f_i^2}{1-Z_{m+1}} - 2
\underbrace{\sqrt{ \sum_{i=1}^m
\frac{Z_i f_i^2}{1-Z_{m+1}}}}_{\geq\sum_{i=1}^m
\frac{Z_i f_i}{1-Z_{m+1}}} f_{m+1} + f_{m+1}^2
\nonumber
\\
&\leq& \sum_{i=1}^m \frac{Z_i}{1-Z_{m+1}}
(f_i - f_{m+1} )^2 .\nonumber
\end{eqnarray}
We obtain in combination with \eqref{elogmeandiffp1} the following estimate:
\[
\Ent_{\nu}(\tilde f) \leq\frac{Z_{m+1}}{\Lambda(Z_{m+1}, 1 -
Z_{m+1})} \sum
_{i=1}^m Z_i (f_i -
f_{m+1} )^2 .
\]
To conclude the assertion, we first note that $1-Z_{m+1}=\sum_{j=1}^m
Z_j \geq Z_j$ for $j=1,\dots,m$. Further, $\Lambda(a,\cdot)$ is
monotone increasing for $a>0$, that is, $\partial_b \Lambda(a,b)>0$
(cf. Appendix~\ref{clogmean}). Both properties imply that $\Lambda
(Z_{m+1},1-Z_{m+1})\geq\Lambda(Z_{m+1},Z_j)$ for $j=1,\dots, m$,
which finally shows \eqref{elogmeandiff}.
\end{pf}
With the help of Lemma~\ref{lemlogmeandiff} we estimate the
coarse-grained entropy $\Ent_{\bar\mu} (\overline{f^2} )$
occurring in the splitting of the entropy \eqref{elemsplitent}.
This generalizes the approach of \cite{Chafai2010}, Section~4.1, to
the case of finite mixtures with more than two components.
%
\begin{lem}[(Estimate of the coarse-grained entropy)]\label
{lemestcoarse-entropy}
The coarse-grained entropy in \eqref{elemsplitcoarse-entropy} can
be estimated by
\begin{eqnarray}
\label{elemestcoarse-entropy} &&\qquad \Ent_{\bar\mu} \bigl(\overline{f^2} \bigr) \nonumber
\\[-8pt]
\\[-8pt]&&\quad \qquad \leq
\sum_{i=1}^M \biggl(\sum
_{j\ne i} \frac{Z_i Z_j \var_{\mu_i}(f)}{\Lambda
(Z_i,Z_j)} + \sum_{j>i}
\frac{Z_i Z_j}{\Lambda(Z_i,Z_j)} \bigl(\Expect_{\mu _i}(f) - \Expect_{\mu_j}(f)
\bigr)^2 \biggr),
\nonumber
\end{eqnarray}
where $\overline{f^2}\dvtx   \{1,\dots, M \} \to\R$ is given by
$\overline{f_i^2} := \Expect_{\mu_i}(f^2)$.
\end{lem}
\begin{pf}
Since $\bar\mu= Z_1 \delta_1 + \cdots+ Z_M \delta_M$ is finite
discrete probability measure, we can apply Lemma~\ref{lemlogmeandiff}
to $\Ent_{\bar\mu}(\overline{f^2})$
\begin{equation}
\label{esplitentp1} \Ent_{\bar\mu}\bigl(\overline{f^2}\bigr) \leq
\sum_{i=1}^m \sum
_{j>i} \frac{Z_i Z_j}{\Lambda(Z_i,Z_j)} \Bigl(\sqrt{\overline
{f^2_i}}-\sqrt{\overline{f^2_j}}
\Bigr)^2 .
\end{equation}
The square-root-mean-difference on the right-hand side of \eqref
{esplitentp1} can be estimated by using the Jensen inequality
\begin{eqnarray}
\label{eestsqumeandiff} \bigl(\sqrt{\Expect_{\mu_i}\bigl(f^2\bigr)} -
\sqrt{\Expect_{\mu
_j}\bigl(f^2\bigr)} \bigr)^2
&\leq& \Expect_{\mu_i}
\bigl(f^2\bigr) - 2\underbrace{\sqrt{\Expect_{\mu
_i}
\bigl(f^2\bigr)\Expect_{\mu_j}\bigl(f^2
\bigr)}}_{\geq\Expect_{\mu_i}(f) \Expect
_{\mu_j}(f)} + \Expect_{\mu_j}\bigl(f^2\bigr)
\nonumber
\\
&\leq& \Expect_{\mu_i}\bigl(f^2\bigr) - 2
\Expect_{\mu_i}(f)\Expect_{\mu
_j}(f) + \Expect_{\mu_j}
\bigl(f^2\bigr)
\\
&=& \var_{\mu_i}(f) + \var_{\mu_j}(f) + \bigl(
\Expect_{\mu _i}(f) - \Expect_{\mu_j}(f) \bigr)^2 .
\nonumber
\end{eqnarray}
Now, we can combine \eqref{esplitentp1} and \eqref
{eestsqumeandiff} to arrive at the desired result \eqref
{elemestcoarse-entropy}.
\end{pf}
A combination of Lemma~\ref{lemsplitvar-ent} and Lemma~\ref
{lemestcoarse-entropy} yields the desired estimate of the entropy in
terms of local variances, local entropies and mean-differences.

\begin{cor}\label{corsplitentropy}
Let $\mu$ have a mixture representation according to Definition~\ref
{defmixturemu}, then the entropy of $f$ w.r.t. $\mu$ can be
estimated by
\begin{eqnarray}
\label{splitequentropy} \Ent_{\mu}\bigl(f^2\bigr) &\leq&\sum
_{i=1}^M Z_i
\Ent_{\mu_i}\bigl(f^2\bigr) + \sum
_{i=1}^M \sum_{j\ne i}
\frac{Z_i Z_j}{\Lambda(Z_i,Z_j)} \var_{\mu
_i}(f)
\nonumber
\\[-8pt]
\\[-8pt]
&&{}+ \sum_{i=1}^M \sum
_{j>i} \frac{Z_i Z_j}{\Lambda(Z_i,Z_j)} \bigl(\Expect_{\mu_i}(f) -
\Expect_{\mu_j}(f) \bigr)^2 .
\nonumber
\end{eqnarray}
\end{cor}

\subsection{Main results}\label{smainres}

The main results of this work are good estimates of the single terms on
the right-hand side of \eqref{elemsplitvar} and \eqref
{splitequentropy}. In detail, we need the \emph{local PI} and the
\emph{local LSI} provided by Theorem~\ref{thmlocalPI} and
Theorem~\ref{thmlocalLSI}. Furthermore, we need good control of the
mean-differences, which will be the content of Theorem~\ref
{thmcoarsemd}. Finally, the Eyring--Kramers formulas of
Corollary~\ref{coreyringkramersPI} and Corollary~\ref
{coreyringkramersLSI} are simple consequences of these
representations and estimates.

\subsubsection{Local Poincar\'{e} and logarithmic Sobolev
inequalities}\label{smainlocal}

Let us now turn to the estimation of the local variances and entropies.
From the heuristic understanding of the process $\xi_t$ given
by \eqref{eoverLang}, we expect a good behavior of the local
Poincar\'{e} and logarithmic Sobolev constant for the local Gibbs
measures $\mu_i$ as it resembles the fast convergence of $\xi_t$ to a
neighborhood of the next local minimum.
Therefore, the local variances and entropies should not contribute to
the leading order expansion of the total Poincar\'{e} and logarithmic
Sobolev constant of $\mu$. This idea is quantified in the next two theorems.
%
\begin{theo}[(Local Poincar\'{e} inequality)]\label{thmlocalPI}
Under Assumption~\ref{assumeenv}, there exists an admissible
partition $\cP_\cM= \{\Omega_i \}_{i=1}^M$ for $\mu$
(cf. Definition~\ref{defpartition}) such that the associated local
Gibbs measures $ \{\mu_i \}_{i=1}^M$, obtained by
restricting $\mu$
to $\Omega_i$ [cf. \eqref{elocalmu}], satisfy $\PI(\varrho_i)$ with
\[
\varrho_i^{-1} = O(\eps) .
\]
\end{theo}
%
%
\begin{theo}[(Local logarithmic Sobolev inequality)]\label{thmlocalLSI}
Under Assumption~\ref{assumeenvLSI} and for the same admissible
partition $\cP_\cM= \{\Omega_i \}_{i=1}^M$ for $\mu$ as in
Theorem~\ref{thmlocalPI}, the associated local Gibbs measures $
\{\mu_i \}_{i=1}^M$, obtained by restricting $\mu$ to $\Omega_i$
[cf. \eqref{elocalmu}], satisfy $\LSI(\alpha_i)$ with
\[
\alpha_i^{-1} = O(1) .
\]
\end{theo}
Even if Theorem~\ref{thmlocalPI} and Theorem~\ref{thmlocalLSI}
are very plausible, their proof is not easy. The reason is that our
situation goes beyond the scope of the standard tools for PI and LSI:
\begin{itemize}
\item The Bakry--\'{E}mery criterion (cf. Theorem~\ref
{localthmBakryEmery}) cannot be applied because we do not have a
convex Hamiltonian.
\item A naive application of the Holley--Stroock perturbation principle
(cf. Theorem~\ref{localthmHolleyStroock}) would yield an
exponentially bad dependence on the parameter~$\varepsilon$.
\item One cannot apply a simple \emph{Lyapunov argument}, because one
cannot impose a drift condition on the boundary of all elements of the
partition $\cP_\cM$, simultaneously.
\end{itemize}
For the proof we apply a subtle combination of a Lyapunov and a
perturbation argument. The core of the argument is an explicit
construction of a \emph{Lyapunov function}. This Lyapunov function has
to satisfy Neumann boundary conditions on the sets $\Omega_i$. By
using the canonical partition $\Omega_i$ into the basins of attraction
of the gradient flow w.r.t. $H$ (see Remark~\ref
{rcanonicalpartition}), the construction of the Lyapunov function
would be technically very demanding. We avoid these difficulties by
choosing another partition $\Omega_i$ such that the Lyapunov function
will automatically satisfy Neumann boundary conditions on $\Omega_i$.
We outline the argument for Theorem~\ref{thmlocalPI} and
Theorem~\ref{thmlocalLSI} in Section~\ref{clocal}.
%
\begin{rem}[(Optimality of Theorem~\ref{thmlocalPI} and Theorem~\ref
{thmlocalLSI})]
The one-dimensional case indicates that the results of Theorem~\ref
{thmlocalPI} and Theorem~\ref{thmlocalLSI} are the best behavior
in $\varepsilon$, which one can expect in general. The optimality in
the one-dimensional case was investigated in \cite{AndrePhDthesis}, Section~3.3, by using the \emph{Muckenhoupt
functional} \cite{Muckenhoupt1972} and \emph{Bobkov--G\"{o}tze
functional} \cite{Bobkov1999a}.
\end{rem}

\subsubsection{Mean-difference estimate}\label{smainmd}

Let us now turn to the estimation of the mean-difference $
(\Expect_{\mu_i}(f) - \Expect_{\mu_j}(f) )^2$.
From the heuristics and the splitting of the variance \eqref
{elemsplitvar} and entropy \eqref{splitequentropy}, we expect to
see in the estimation of the mean-difference the exponential long
waiting times of the jumps of the diffusion $\xi_t$ given by \eqref
{eoverLang} between the elements of the partition $\cP_\cM$.
We have to find a good upper bound for the constant $C$ in the inequality
\[
\bigl(\Expect_{\mu_i}(f)-\Expect_{\mu_j}(f) \bigr)^2
\leq C \int \llvert \nabla f\rrvert ^2 \,\dx{\mu} .
\]
For this purpose, we introduce in Section~\ref{sectranspmd} a \emph
{weighted transport distance} between probability measures which yields
a variational bound on the constant $C$. By an approximation argument
(cf. Section~\ref{sapproxmd}), we give an explicit construction of a
transport interpolation (cf. Section~\ref{safftransp}), which allows
for asymptotically sharp estimates of the constant $C$.
%
\begin{theo}[(Mean-difference estimate)]\label{thmcoarsemd}
Let $H$ satisfy Assump-\break tion~\ref{assumenondegenerate} and let $\cP
_\cM= \{\Omega_i \}_{i=1}^M$ be an admissible partition
for $\mu$
(cf. Definition~\ref{defpartition}). Moreover, assume that each local
Gibbs measure $\mu_i$ of the mixture representation of $\mu$ (cf.
Definition~\ref{defmixturemu}) satisfy $\PI(\varrho_i)$ with
$\varrho_{i}^{-1}=O(\eps)$. Then the mean-differences between the
local Gibbs measures $\mu_i$ and $\mu_j$ for $i=1,\dots, M-1$ and
$j=i+1,\dots, M$ satisfy
\begin{eqnarray}
\label{thmcoarsemde}\quad  &&\bigl(\Expect_{\mu_i}(f) - \Expect_{\mu_j}(f)
\bigr)^2 \nonumber
\\[-8pt]
\\[-8pt]&&\qquad \lesssim \frac{Z_\mu}{(2\pi\eps)^{\fraca{n}{2}}} \frac{2\pi\eps\sqrt {\llvert \det\nabla^2 H(s_{i,j})\rrvert }}{\llvert
\lambda^-(s_{i,j})\rrvert } \exp \biggl(
\frac{H(s_{i,j})}{\eps} \biggr) \int\llvert \nabla f\rrvert ^2\,\dx{\mu} ,
\nonumber
\end{eqnarray}
where $\lambda^-(s_{i,j})$ denotes the negative eigenvalue of the
Hessian $\nabla^2 H(s_{i,j})$ at the communicating saddle $s_{i,j}$
defined in Assumption~\ref{assumenondegenerate}. The symbol $\lesssim
$ means $\leq$ up to a multiplicative error term of the form
\[
1+ O\bigl(\sqrt{\eps} \llvert \log{\eps}\rrvert ^{\fraca
{3}{2}}\bigr) .
\]
\end{theo}
The proof of Theorem~\ref{thmcoarsemd} is carried out in full detail
in Section~\ref{CMEANDIFF}.
%
\begin{rem}[(Multiple minimal saddles)]
In Assumption~\ref{assumenondegenerate}, we demand that there is
exactly one minimal saddle between the local minima $m_i$ and $m_j$.
The technique developed in Section~\ref{CMEANDIFF} is flexible enough
to handle also cases, in which there exists more than one minimal
saddle between local minima. The according adaptions and the resulting
theorem can be found in \cite{AndrePhDthesis}, Section~4.5.
\end{rem}
%
%
\begin{rem}[(Relation to capacity)]
The quantity on the right-hand side of \eqref{thmcoarsemde} is the
inverse of the capacity of a small neighborhood around $m_i$ w.r.t. to
a small neighborhood around $m_j$. The capacity is the crucial
ingredient of the works \cite{Bovier2004} and \cite{Bovier2005}.
\end{rem}

\subsubsection{Eyring--Kramers formulas}\label{sEyringKramers}

Now, let us turn to the Eyring--Kramers formula. Starting from the
splitting obtained in Lemma~\ref{lemsplitvar-ent} and Corollary~\ref
{corsplitentropy} a combination of Theorem~\ref{thmlocalPI},
Theorem~\ref{thmlocalLSI} and Theorem~\ref{thmcoarsemd}
immediately leads to the multidimensional Eyring--Kramers formula for
the PI (cf. \cite{Bovier2005}, Theorem~1.2) and LSI.
%
\begin{cor}[(Eyring--Kramers formula for Poincar\'{e} inequality)]\label
{coreyringkramersPI}
Under Assumptions~\ref{assumeenv} and~\ref{assumenondegenerate},
the measure $\mu$ satisfies \ref{PIvarrho} with
\begin{equation}
\label{eeyringkramersPI} \frac{1}{\varrho}\lesssim Z_1 Z_2
\frac{Z_\mu}{(2\pi\eps)^{\fraca
{n}{2}}} \frac{2\pi\eps\sqrt{\llvert \det\nabla
^2(H(s_{1,2}))\rrvert }}{\llvert \lambda^-(s_{1,2})\rrvert } \exp \biggl(\frac {H(s_{1,2})}{\eps} \biggr) ,
\end{equation}
where $\lambda^-(s_{1,2})$ denotes the negative eigenvalue of the
Hessian $\nabla^2 H(s_{1,2})$ at the communicating saddle $s_{1,2}$.
Further, the order is given such that $H(m_1)\leq H(m_i)$ and
$H(s_{1,2})-H(m_2)$ is the energy barrier of the system in the sense of
Assumption~\ref{assumenondegenerate}. The prefactors $Z_i$ are given
by the relation
\begin{equation}
\label{eeyringkramersPIpre} Z_i Z_\mu\approx\frac{(2\pi\eps)^{\fraca{n}{2}}}{\sqrt{\det\nabla
^2 H(m_i)}} \exp
\biggl(-\frac{H(m_i)}{\eps} \biggr).
\end{equation}
\end{cor}
\begin{pf}
Using the admissible partition $\cP_\cM$ from Theorem~\ref
{thmlocalPI} we decompose the variance into local variances and
mean-differences given by Lemma~\ref{lemsplitvar-ent}. An
application of Theorem~\ref{thmlocalPI} and Theorem~\ref
{thmcoarsemd} yields the estimate
\begin{eqnarray}
\label{coreyringkramersPIp1} \var_\mu(f) &\leq& \sum_{i}
Z_i \var_{\mu_i}(f) + \sum_{i}
\sum_{j<i} Z_i Z_j \bigl(
\Expect_{\mu_i}(f) - \Expect_{\mu
_j}(f) \bigr)^2
\nonumber
\\
&\lesssim& \biggl(O(\eps) + \sum_{i} \sum
_{j>i} \frac{Z_i Z_j
Z_\mu}{(2\pi\eps)^{\fraca{n}{2}}} \frac{2\pi\eps\sqrt{\llvert \det \nabla^2 H(s_{i,j})\rrvert }}{\llvert \lambda
^-(s_{i,j})\rrvert } \exp
\biggl(\frac{H(s_{i,j})}{\eps} \biggr) \biggr)\\
&&{}\times \int\llvert \nabla f\rrvert
^2 \,\dx {\mu} .
\nonumber
\end{eqnarray}
The final step is to observe that by Assumption~\ref
{assumenondegenerate} the exponential dominating term in \eqref
{coreyringkramersPIp1} is given for $i=1$ and $j=2$. The precise
form of the prefactors $Z_i$ is obtained from \eqref{ePartSumAdmPart}
in Definition~\ref{defpartition}.
\end{pf}
In \cite{Bovier2005}, Theorem~1.2, it is also shown that the upper
bound of \eqref{eeyringkramersPI} is optimal by an approximation of
the harmonic function. Therefore, in the following we can assume
that \eqref{eeyringkramersPI} holds with $\approx$ instead of
$\lesssim$.
%
\begin{rem}[(Higher exponentially small eigenvalues)] The main result of
\cite{Bovier2005}, Theorem~1.2, does not only characterize the second
eigenvalue of $L$ but also the higher exponentially small eigenvalues.
In principle, these characterizations can be also obtained in the
present approach: The dominating exponential modes in \eqref
{coreyringkramersPIp1}, that is, those obtained by setting $i=1$,
correspond to the inverse eigenvalues of $L$ for $j=2,\dots,M$. By
using the variational characterization of the eigenvalues of the
operator $L$, the other exponentially small eigenvalues may be obtained
by restricting the class of test functions $f$ to the orthogonal
complement of the eigenspaces of smaller eigenvalues.
\end{rem}

%
\begin{cor}[(Eyring--Kramers formula for logarithmic Sobolev
inequalities)]\label{coreyringkramersLSI}
Under Assumptions~\ref{assumeenvLSI} and~\ref{assumenondegenerate},
the measure $\mu$ satisfies \ref{LSIalpha} with
\begin{eqnarray}
\label{eeyringkramersLSI} \frac{2}{\alpha}&\lesssim&\frac{Z_1 Z_2}{\Lambda(Z_1,Z_2)} \frac
{Z_\mu}{(2\pi\eps)^{\fraca{n}{2}}}
\frac{2\pi\eps\sqrt{\llvert \det\nabla^2(H(s_{1,2}))\rrvert }}{\llvert \lambda
^-(s_{1,2})\rrvert } \exp \biggl(\frac{H(s_{1,2})}{\eps} \biggr)\nonumber
\\[-8pt]
\\[-8pt] &\approx&
\frac{1}{\Lambda
(Z_1,Z_2)} \frac{1}{\varrho} ,
\nonumber
\end{eqnarray}
where the occurring constants are like in Corollary~\ref
{coreyringkramersPI} and $\Lambda(Z_1,Z_2)$ denotes the logarithmic
mean (cf. Appendix~\ref{clogmean})
\[
\Lambda(Z_1, Z_2) = \frac{Z_1 - Z_2}{\log Z_1 - \log Z_2}.
\]
%
\end{cor}
\begin{pf}
Using the admissible partition $\cP_\cM$ from Theorem~\ref
{thmlocalPI} and Theorem~\ref{thmlocalLSI}, we decompose the
Entropy according to Corollary~\ref{corsplitentropy}. From there, we
estimate the local entropies and variances as well as the
mean-differences by using Theorem~\ref{thmlocalPI}, Theorem~\ref
{thmlocalLSI} and Theorem~\ref{thmcoarsemd}. Overall, this yields
the estimate
\begin{eqnarray}
\label{ecoreyringkramersLSIp1} \Ent_{\mu}\bigl(f^2\bigr) &\leq& O(1) \sum
_{i=1}^M Z_i \int\llvert
\nabla f\rrvert ^2 \,\dx{\mu_i} + O(\eps) \sum
_{i=1}^M \sum_{j\ne i}
\frac{Z_i
Z_j}{\Lambda(Z_i, Z_j)} \int\llvert \nabla f\rrvert ^2 \,\dx {
\mu_i}
\nonumber
\\
&&{} + \sum_{i=1}^M \sum
_{j>i} \frac{Z_i Z_j}{\Lambda(Z_i, Z_j)} \frac{Z_\mu}{(2\pi\eps)^{\fraca{n}{2}}}
\frac{2\pi\eps\sqrt {\llvert \det\nabla^2 H(s_{i,j})\rrvert }}{\llvert
\lambda^-(s_{i,j})\rrvert } \exp \biggl(\frac{H(s_{i,j})}{\eps} \biggr)\\
&&\hphantom{{}+{}}{}\times \int\llvert \nabla f
\rrvert ^2 \,\dx {\mu} .
\nonumber
\end{eqnarray}
The first term on the right-hand side of \eqref
{ecoreyringkramersLSIp1} can be rewritten as\break  $O(1) \int\llvert
\nabla f\rrvert ^2\,\dx{\mu}$. For estimating the second term
in \eqref
{ecoreyringkramersLSIp1}, we argue that its prefactor can be
estimated as
\begin{equation}
\label{eEKLSIsecondprefactor} \sum_{i=1}^M \sum
_{j\ne i} \frac{Z_i Z_j}{\Lambda(Z_i,Z_j)} \stackrel{\mathclap{
\eqref{ecoreyringkramersLSIp4}}} {\lesssim} M \sum
_{i=1}^M Z_i O\bigl(
\eps^{-1}\bigr)= O\bigl(\eps^{-1}\bigr) .
\end{equation}
Indeed, using the one-homogeneity of $\Lambda(\cdot,\cdot)$
(cf. Appendix~\ref{clogmean}) yields
\[
\label{ecoreyringkramersLSIp2} \frac{Z_i Z_j}{\Lambda(Z_i, Z_j)} = Z_i \frac{\log(\fraca
{Z_i}{Z_j})}{\fraca{Z_i}{Z_j}-1} =
Z_i P \biggl(\frac{Z_i}{Z_j} \biggr)   \qquad\mbox{where } P(x) :=
\frac{\log x}{x-1} .
\]
The function $P(x)$ is decreasing and has a logarithmic singularity at
$0$. Therefore, using the characterization of the partitions sums $Z_i$
from \eqref{eeyringkramersPIpre} yields the identity
\begin{equation}
\label{ecoreyringkramersLSIp3b} \frac{Z_i}{Z_j} = \frac{Z_i Z_\mu}{Z_j Z_\mu} \stackrel{\mathclap {
\eqref{eeyringkramersPIpre}}} {\approx} \frac{\sqrt{\nabla^2
H(m_j)}}{\sqrt{\nabla^2 H(m_i)}} \exp \biggl(-
\frac{H(m_i)-
H(m_j)}{\eps} \biggr) ,
\end{equation}
which becomes exponentially small provided that $H(m_i)>H(m_j)$. Hence,
the logarithmic mean can be estimated as
\begin{equation}
\label{ecoreyringkramersLSIp4} \frac{Z_i Z_j}{\Lambda(Z_i, Z_j)} = Z_i P{ \biggl(\frac
{Z_i}{Z_j}
\biggr)} \lesssim Z_i O\bigl(\eps^{-1}\bigr)
\end{equation}
implying the desired estimate \eqref{eEKLSIsecondprefactor}.
Therefore, the second term in \eqref{ecoreyringkramersLSIp1} can
be estimated by $O(1) \int\llvert \nabla f\rrvert ^2 \,\dx{\mu
}$. The third
term dominates the first two terms on an exponential scale. This leads
to the estimate
\begin{eqnarray*}
\Ent_{\mu}\bigl(f^2\bigr) &\lesssim&\sum
_{i=1}^M \sum_{j>i}
\frac{Z_i
Z_j}{\Lambda(Z_i, Z_j)} \frac{Z_\mu}{(2\pi\eps)^{\fraca{n}{2}}} \frac{2\pi\eps\sqrt{\llvert \det\nabla^2 H(s_{i,j})\rrvert }}{\llvert \lambda^-(s_{i,j})\rrvert } e^{\fraca
{H(s_{i,j})}{\eps}}\\
&&{}\times \int
\llvert \nabla f\rrvert ^2 \,\dx {\mu} .
\end{eqnarray*}
From Assumption~\ref{assumenondegenerate} together with \eqref
{ecoreyringkramersLSIp4} follows that the exponentially leading
order term is given for $i=1$ and $j=2$.
\end{pf}
The Eyring--Kramers formula for the PI and LSI stated in Corollary~\ref
{coreyringkramersPI} and Corollary~\ref{coreyringkramersLSI} are
still implicit. To obtain an explicit formula, one still has insert the
asymptotic expansion for the partition functions $Z_1$, $Z_2$, and
$Z_\mu$. The expression for $Z_\mu$ depends on the number of global
minima of the Hamiltonian $H$. Therefore, one has to consider several
cases in order to obtain the explicit Eyring--Kramer formula. In the
following corollary, we look at two special cases: In the first case,
there is only one unique global minimum. In the second case, there are
two global minima. In both cases, the dominating term scales
exponentially in the saddle height, but it is surprising that the
scaling in $\eps$ of the exponential pre factor for the LSI constant changes.
%
\begin{cor}[(Comparison of $\varrho$ and $\alpha$ in special
cases)]\label{coreyringkramerscomp}
Let us state two specific cases of \eqref{eeyringkramersPI}
and \eqref{eeyringkramersLSI}. Therefore, let $ \{\kappa
_i^2 \}_{i=1}^M$ be given by
\begin{equation}
\label{dkappaI} \kappa_i^2 := \det\nabla^2
H(m_i).
\end{equation}
On the one hand, if one has one unique global minimum,
namely $H(m_1)<H(m_i)$ for $i \in \{ 2, \ldots, M  \}$, it holds
\begin{eqnarray}
\frac{1}{\varrho} &\approx& \frac{1}{\kappa_2} \frac{2\pi\eps
\sqrt{\llvert \det\nabla^2(H(s_{1,2}))\rrvert }}{\llvert \lambda^-(s_{1,2})\rrvert } \exp \biggl(
\frac{H(s_{1,2})-H(m_2)}{\eps} \biggr),\label
{eeyringkramersPInondeg}
\\
\frac{2}{\alpha} &\lesssim& \biggl(\frac{H(m_2)-H(m_1)}{\eps} + \log \biggl(
\frac{\kappa_1}{\kappa_2} \biggr) \biggr) \frac
{1}{\varrho} .\label{eeyringkramersLSInondeg}
\end{eqnarray}
On the other hand, if $H(m_1)=H(m_2)< H(m_i)$ for $i \in \{ 3,
\ldots, M  \}$, it holds
\begin{eqnarray}\qquad
\frac{1}{\varrho}&\approx& \frac{1}{\kappa_1 + \kappa_2} \frac
{2\pi\eps\sqrt{\llvert \det\nabla^2(H(s_{1,2}))\rrvert
}}{\llvert \lambda ^-(s_{1,2})\rrvert } \exp \biggl(
\frac
{H(s_{1,2})-H(m_2)}{\eps} \biggr), \label{eeyringkramersPIdeg}
\\
\frac{2}{\alpha} &\lesssim& \frac{1}{\Lambda (\kappa
_1,\kappa_2 )} \frac{2\pi\eps\sqrt{\llvert \det\nabla
^2(H(s_{1,2}))\rrvert }}{\llvert \lambda^-(s_{1,2})\rrvert } \exp \biggl(
\frac {H(s_{1,2})-H(m_2)}{\eps} \biggr). \label
{eeyringkramersLSIdeg}
\end{eqnarray}
\end{cor}
\begin{pf}
By \eqref{eeyringkramersPI}, we still have to estimate nonexplicit
factor $\frac{ Z_1 Z_2Z_\mu}{(2\pi\eps)^{\fraca{n}{2}}}$. If
$H(m_1)<H(m_2)$, then it holds $Z_1 = 1+O (e^{-\fracb
{H(m_2)-H(m_1)}{\eps}} )$. The factor $Z_2 Z_\mu$ is given
by \eqref
{eeyringkramersPIpre} and we obtain
\[
\frac{Z_1 Z_2 Z_\mu}{(2\pi\eps)^{\fraca{n}{2}}} \approx\frac
{1}{\sqrt{\det\nabla^2 H(m_2)}} \exp \biggl(-\frac{H(m_2)}{\eps
}
\biggr) ,
\]
which leads to \eqref{eeyringkramersPInondeg}. For the LSI, we
additionally have to evaluate the factor $\frac{1}{\Lambda(Z_i,Z_j)}$
which can be done with the help of \eqref{ecoreyringkramersLSIp3b}
\begin{eqnarray}
\frac{1}{\Lambda(Z_i, Z_j)} & =& \log \biggl(\frac{Z_i}{Z_j} \biggr) \biggl(1+O
\biggl(\exp \biggl(-\frac{H(m_2)-H(m_1)}{\eps } \biggr) \biggr) \biggr)
\nonumber
\\
&\overset{\mathclap{\eqref{ecoreyringkramersLSIp3b}}} {\approx}& \log
\biggl(\frac{\sqrt{\nabla^2 H(m_j)}}{\sqrt{\nabla^2 H(m_i)}} \exp \biggl(-\frac{H(m_i)- H(m_j)}{\eps} \biggr) \biggr) .\nonumber
\end{eqnarray}
That is already the estimate \eqref{eeyringkramersLSInondeg}.

Let us turn now to the case $H(m_1)=H(m_2)< H(m_3)$. Then it holds
$Z_1+Z_2=1+O (e^{-\fracb{H(m_2)-H(m_1)}{\eps}} )$. In
particular it holds $Z_\mu\approx Z_1Z_\mu+ Z_2 Z_\mu$. Therewith,
we can evaluate
the factor $Z_1 Z_2 \frac{Z_\mu}{(2\pi\eps)^{\fraca{n}{2}}}$ by
using \eqref{eeyringkramersPIpre}
\begin{eqnarray}
Z_1 Z_2 \frac{Z_\mu}{(2\pi\eps)^{\fraca{n}{2}}} &=&
\frac{(2\pi
\eps)^{\fraca{n}{2}}}{Z_\mu} \frac{Z_1 Z_\mu}{(2\pi\eps)^{\fraca
{n}{2}}} \frac{Z_2 Z_\mu}{(2\pi\eps)^{\fraca{n}{2}}} \nonumber\\
&\approx&
\frac{(2\pi\eps)^{\fraca{n}{2}}}{Z_1 Z_\mu+Z_2 Z_\mu} \frac{Z_1
Z_\mu}{(2\pi\eps)^{\fraca{n}{2}}} \frac{Z_2 Z_\mu}{(2\pi\eps
)^{\fraca{n}{2}}}
\nonumber
\\
&\stackrel{\mathclap{\eqref{eeyringkramersPIpre}}} {=}& \frac
{1}{\fraca{1}{\kappa_1}+\fraca{1}{\kappa_2}}
\frac{1}{\kappa_1} \frac{1}{\kappa_2} = \frac{1}{\kappa_1 + \kappa_2},\nonumber
\end{eqnarray}
which precisely leads to the expression \eqref
{eeyringkramersPIdeg}. By using the homogeneity of $\Lambda(\cdot
,\cdot)$ (cf. Appendix~\ref{clogmean}) and again \eqref
{eeyringkramersPIpre}, it follows for the LSI
\[
\frac{Z_1 Z_2}{\Lambda(Z_1,Z_2)} \frac{Z_\mu}{(2\pi\eps)^{\fraca
{n}{2}}} = \frac{1}{\Lambda (\fracc{(2\pi\eps)^{\fraca
{n}{2}}}{Z_2 Z_\mu}, \fracc{(2\pi\eps)^{\fraca{n}{2}}}{Z_1 Z_\mu
} )} =
\frac{1}{\Lambda(\kappa_2,\kappa_1)} .
\]
Finally, the result \eqref{eeyringkramersLSIdeg} is a consequence
of the symmetry of $\Lambda(\cdot,\cdot)$.
\end{pf}

%
\begin{rem}[(Identification of $\alpha$ and $\varrho$)]
Remark~\ref{remLSIPI} shows that always $\alpha\leq\varrho$. We
want to compare this to the case $H(m_1)=H(m_2)$. Comparing \eqref
{eeyringkramersPIdeg} and \eqref{eeyringkramersLSIdeg}, we observe
\begin{equation}
\label{esglsidegest} 1\leq\frac{\varrho}{\alpha} \lesssim\frac{\fracb{\kappa_1 +
\kappa_2}{2}}{\Lambda(\kappa_1,\kappa_2)} ,
\end{equation}
where the constant $\kappa_1$ and $\kappa_2$ are given by \eqref
{dkappaI}. The right-hand side of \eqref{esglsidegest} consists
of an quotient of the arithmetic and the logarithmic mean. The lower
bound of $1$ can also attained by an application of the
logarithmic-arithmetic mean inequality from Lemma~\ref
{logmeanlemgeologar}. Moreover, equality only holds for $\kappa
_1=\kappa_2$. Hence, only in the symmetric case $\varrho\approx
\alpha$.
\end{rem}
%
%
\begin{rem}[(Relation to mixtures)]\label{remrelmix}
If $H(m_1)<H(m_2)$, then \eqref{eeyringkramersLSInondeg} gives
\begin{equation}
\label{eremrelmixlsi} \frac{\varrho}{\alpha}\lesssim\frac{1}{2} \log \biggl(
\frac
{\kappa_2}{\kappa_1}e^{\fracb{H(m_2)-H(m_1)}{\eps}} \biggr) \approx \frac
{1}{2} \llvert
\log Z_2\rrvert  \qquad \mbox{where } Z_2 = \mu(
\Omega_2)  \hspace*{-30pt}
\end{equation}
which shows an inverse scaling in $\eps$. A different scaling behavior
between the Poincar\'{e} and logarithmic Sobolev constant was also
observed by Chafa\"i and Malrieu \cite{Chafai2010} in a different
context. They consider mixtures of probability measures $\nu_0$ and
$\nu_1$ satisfying $\PI(\varrho_i)$ and $\LSI(\alpha_i)$, that
is, for $p\in[0,1]$ the measure $\nu_p$ given by
\[
\nu_p = p \nu_0 + (1-p) \nu_1 .
\]
They deduce conditions under which also $\nu_p$ satisfies $\PI
(\varrho_p)$ and $\LSI(\alpha_p)$ and give bounds on the constants.
They give one-dimensional examples where the Poincar\'{e} constant
stays bounded, whereas the logarithmic Sobolev constant blows up
logarithmically in the mixture parameter $p$ going to $0$ or $1$. The
common feature of the examples they deal with is $\nu_1\ll\nu_2$ or
$\nu_2\ll\nu_1$. This case can be generalized to the
multidimensional case, where also a different scaling of the Poincar\'
{e} and logarithmic Sobolev constants is observed. The details can be
found in \cite{AndrePhDthesis}, Chapter~6.

In the present case, the Gibbs measure $\mu$ has also a mixture
representation \eqref{emixturemu}. In the two-component case, it has
the form
\[
\mu= Z_1 \mu_1 + Z_2 \mu_2 .
\]
Let us emphasize, that $\mu_1\perp\mu_2$. The estimate \eqref
{eremrelmixlsi} also shows a logarithmic blow-up in the mixture
parameter $Z_2$ for the ratio of the Poincar\'{e} and the logarithmic
Sobolev constant.
\end{rem}

\subsection{Optimality of the logarithmic Sobolev constant in one
dimension}\label{smainopt}

\begin{figure}

\includegraphics{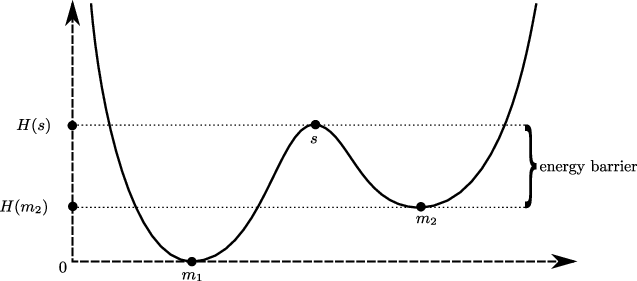}

\caption{Double-well potential $H$ on $\R$ (labeled).}
\label{figdouble-wellR}
\end{figure}

In this section, we give a strong indication that the result of
Corollary~\ref{coreyringkramersLSI} is optimal. We explicitly
construct a function attaining equality in \eqref{eeyringkramersLSI}
for the one-dimensional case.
For this purpose, let $\mu$ be a probability measure on $\R$ having
as Hamiltonian $H$ a generic double-well (cp. Figure~\ref{figdouble-wellR}). Namely, $H$ has two minima $m_1$ and $m_2$ with
$H(m_1)\leq H(m_2)$ and a saddle $s$ in-between. Then Theorem~\ref
{coreyringkramersLSI} shows
\begin{equation}
\label{elsimin1d} \inf_{g\dvtx  \int g^2 \,\dx{\mu}=1 } \frac{\int(g')^2 \,\dx{\mu}}{\int
g^2 \log g^2 \,\dx{\mu}} \gtrsim
\frac{\Lambda(Z_1,Z_2)}{Z_1 Z_2} \frac{\sqrt{2\pi\eps}}{Z_\mu} \frac{\sqrt{\llvert
H''(s)\rrvert }}{2 \pi
\eps} e^{-\fraca{H(s)}{\eps}} .
\end{equation}

We construct a function $g$ attaining the lower bound given by \eqref
{elsimin1d}. We make the following ansatz for the function $g$: We
define $g$ on a small $\delta$-neighborhood around the minima $m_1,
m_2$ and the saddle $s$:
\[
g(x) := \cases{\displaystyle g(m_1) ,&\quad$x\in
B_\delta(m_1) $,
\cr
\displaystyle g(m_1) +
\frac{g(m_2)-g(m_1)}{\sqrt{2\pi\eps\sigma}} \int_{m_1}^x e^{-\fracc{(y-s)^2}{2\sigma\eps}}
\,\dx{y} ,&\quad$x\in B_\delta(s) $,
\cr
\displaystyle g(m_2)
,&\quad$x\in B_\delta(m_2) $. }
\]
The ansatz depends on the parameters $g(m_1)$, $g(m_2)$ and $\sigma$.
In between the $\delta$-neighborhoods, the function $g$ is smoothly
extended in a monotone fashion.
%

The measure $\mu$ is the usual Gibbs measure as in \eqref{dgibbsmeasure}.
We fix $Z_\mu$ by assuming that $H(m_1)=0$. We represent $\mu$ as the mixture
\[
\mu= Z_1 \mu_1 + Z_2 \mu_2
\qquad\mbox{where } \mu_1 := \mu \llcorner\Omega_1
\mbox{ and } \mu_2 := \mu\llcorner \Omega_2,
\]
hereby, $\Omega_1 := (-\infty, s)$ and $\Omega_2 := (s,\infty)$ and
$Z_i := \mu(\Omega_i)$ for $i=1,2$, which implies $Z_1+Z_2=1$.
Using via an asymptotic evaluation of $\int g^2  \,\dx{\mu}$ one gets
\[
\label{eoptlsi1dnorm} \int g^2 \,\dx{\mu} \approx Z_1
g^2(m_1) + Z_2 g^2(m_2)
\stackrel {!} {=} 1 .
\]
This motivates the choice
\[
\label{eoptlsi1dchoice} g^2(m_1) = \frac{\tau}{Z_1} \quad
\mbox{and}\quad g^2(m_2) = \frac
{1-\tau}{Z_2} =
\frac{1-\tau}{1-Z_1} \qquad \mbox{for some $\tau \in[0,1]$.}
\]
%
Let us now calculate the denominator of \eqref{elsimin1d}
\begin{equation}
\label{eoptlsi1denteva} %
\int g^2 \log g^2 \,\dx{\mu} =
\tau\log\frac{\tau}{Z_1} + (1-\tau)\log\frac{1-\tau}{Z_2} .
\end{equation}
The final step is to evaluate the Dirichlet energy $\int(g')^2  \,\dx
{\mu}$. Therefore, we do a Taylor expansion of $H$ around $s$.
Furthermore, since $s$ is a saddle, it holds $H''(s)<0$
\begin{eqnarray}
\label{eoptlsi1ddirichletest} %
\int\bigl(g'\bigr)^2 \,\dx{
\mu} &\approx& \frac{(g(m_2)-g(m_1))^2}{Z_\mu
2\pi\eps\sigma} \int_{B_\delta(s)} e^{-\fracc{(x-s)^2}{\sigma\eps
}- \fraca{H(x)}{\eps}}
\,\dx{x}
\nonumber
\\
&\approx& \frac{(g(m_2)-g(m_1))^2}{ Z_\mu 2 \pi\eps\sigma} \int_{B_\delta(s)} e^{-\fraca{(\fraca{(x-s)^2}{\sigma} +
H(s) + H''(s) \fraca{(x-s)^2}{2} )}{\eps} } \,
\dx{x}
\nonumber
\\[-8pt]
\\[-8pt]
&\approx& \frac{(g(m_2)-g(m_1))^2}{ Z_\mu 2 \pi\eps\sigma} e^{-\fraca{H(s)}{\eps}} \int_{B_\delta(s)}
e^{-(\fracc{(x-s)^2}{2\eps
}) (\fraca{2}{\sigma}+H''(s) )} \,\dx{x}
\nonumber
\\
&\approx& \biggl(\sqrt{\frac{\tau}{Z_1}}-\sqrt{\frac{1-\tau
}{Z_2}}
\biggr)^2 \frac{\sqrt{2\pi\eps}}{Z_\mu} e^{-\fraca
{H(s)}{\eps
}} \frac{1}{ 2 \pi\eps}
\frac{1}{\sigma\sqrt{\fraca{2}{\sigma
}+H''(s)}} , %
\nonumber
\end{eqnarray}
where we assume that $\sigma$ is small enough such that $\frac
{2}{\sigma}+H''(s)>0$. The last step is to minimize the right-hand
side of \eqref{eoptlsi1ddirichletest} in $\sigma$, which means
to maximize the expression $2\sigma+ \sigma^2 H''(s)$ in $\sigma$.
Elementary calculus results in $\sigma= - \frac{1}{H''(s)}= \frac
{1}{\llvert  H''(s)\rrvert }>0$ and, therefore,
\begin{equation}
\label{eoptlsi1ddirichleteval} \int\bigl(g'\bigr)^2 \,\dx{\mu} \approx
\biggl(\sqrt{\frac{\tau
}{Z_1}}-\sqrt{\frac{1-\tau}{Z_2}} \biggr)^2
\frac{\sqrt{2\pi\eps
}}{Z_\mu} \frac{\sqrt{\llvert  H''(s)\rrvert }}{ 2 \pi\eps
} e^{-\fraca
{H(s)}{\eps}} .
\end{equation}
Hence, we have constructed by combining \eqref{eoptlsi1denteva}
and \eqref{eoptlsi1ddirichleteval} an upper bound for the
optimization problem \eqref{elsimin1d} given by
\[
\min_{\tau\in(0,1)} \biggl(\frac{ (\sqrt{\fraca{\tau
}{Z_1}}-\sqrt{\fracb{1-\tau}{Z_2}} )^2}{ \tau\log(\fraca{\tau
}{Z_1}) + (1-\tau)\log(\fracb{1-\tau}{Z_2})}
\biggr) \frac{\sqrt {2\pi\eps
}}{Z_\mu} \frac{\sqrt{\llvert  H''(s)\rrvert }}{ 2 \pi\eps
} e^{-\fraca
{H(s)}{\eps}} .
\]
Note that the parameter $\tau\in(0,1)$ is still free. The minimum in
$\tau$ is attained at $\tau= Z_2$ according to Lemma~\ref
{logmeanlemopt} yielding the desired statement
\[
\min_{\tau\in(0,1)} \frac{ (\sqrt{\fraca{Z_2}{Z_1}}-\sqrt {\fraca{Z_1}{Z_2}} )^2}{ Z_2 \log(\fraca{Z_2}{Z_1}) + Z_1\log
(\fraca
{Z_1}{Z_2})} =
\frac{\Lambda(Z_1,Z_2)}{Z_1 Z_2} .
\]
%

\section{Local Poincar\'{e} and logarithmic Sobolev
inequalities}\label{clocal}

In this section, we proof the local PI of Theorem~\ref{thmlocalPI}
and the local LSI of Theorem~\ref{thmlocalLSI}. Even if the choice
of a specific admissible partition $\Omega_i$ of the space $\R^n$
will be crucial, let us for the moment assume that the partition
$\Omega_i$ is given by the basins of attraction of the deterministic
gradient flow (cf. Remark~\ref{rcanonicalpartition}).

There are standard criteria to deduce the $\PI$ or the $\LSI$.
Unfortunately, these criteria do not apply to our situation. Let us
consider the \emph{Bakry--\'{E}mery criterion} and the \emph
{Holley--Stroock perturbation principle}. The \emph{Bakry--\'{E}mery
criterion} connects convexity of the Hamiltonian to the validity of
the PI and the LSI.
%
\begin{theo}[Bakry--\'{E}mery criterion {\cite{BE}, Proposition~3, Corollaire 2}]\label{localthmBakryEmery}
Let $H\dvtx D \to\R$ be a Hamiltonian with Gibbs measure
\[
\mu(\dx
{x})=Z_\mu^{-1} \exp (-\eps^{-1} H(x) )\,\dx{x}
\]
 on a convex
domain $D$ and assume that $\nabla^2 H(x) \geq\lambda>0$ for all
$x\in\R^n$. Then $\mu$ satisfies \ref{PIvarrho} and \ref{LSIalpha} with
\[
\varrho\geq\frac{\lambda}{\eps} \quad\mbox{and}\quad\alpha \geq
\frac{\lambda}{\eps} .
\]
\end{theo}
One cannot apply the criterion of Bakry--\'Emery \cite{BE} to our
situation, because $H$ is not convex on the elements $\Omega$ of the
admissible partition (cf. Definition~\ref{defpartition}). Moreover,
the elements $\Omega\in\cP_\cM$ are not convex in general.

In nonconvex cases, the standard tool to deduce the $\PI$ and the
$\LSI$ is the \emph{Holley--Stroock perturbation principle}.
%
\begin{theo}[Holley--Stroock perturbation principle {\cite{HS}, p. 1184}]\label{localthmHolleyStroock}
Let $H$ be a Hamiltonian with Gibbs measure $\mu(\dx{x})=Z_\mu^{-1}
\exp (-\eps^{-1} H(x) )\,\dx{x}$. Further, let $\tilde H$
denote a bounded perturbation of $H$ and let $\tilde\mu_\eps$ denote
the Gibbs measure associated to the Hamiltonian $\tilde H$.
If $\mu$ satisfies \ref{PIvarrho} or \ref{LSIalpha} then also
$\tilde\mu$ satisfy $\PI(\tilde\varrho)$ or $\LSI(\tilde\alpha
)$ respectively, where the constants satisfy the bounds
\[
\label{localeHolleyStroockPI-LSI} \tilde\varrho\geq\exp \bigl(-\eps^{-1} \osc\psi \bigr)
\varrho \quad\mbox{and}\quad\tilde\alpha\geq\exp \bigl(-\eps^{-1}
\osc(H- \tilde H) \bigr) \alpha,
\]
where $\osc(H - \tilde H) := \sup(H - \tilde H) - \inf(H - \tilde H)$.
\end{theo}
The perturbation principle of Holley--Stroock \cite{HS} allows to
deduce the PI and the LSI constants of nonconvex Hamiltonians from the
PI and the LSI of an appropriately convexfied Hamiltonian. However due
to its perturbative nature, a naive application Theorem~\eqref
{localthmHolleyStroock} would yield an exponential dependence of the
PI and the LSI constant on $\varepsilon$.
%

An important observation for our argument is that the perturbation
principle of Holley--Stroock can still be useful, if applied in a
careful way: Assume for a moment that the perturbed Hamiltonian $\tilde
H_\eps$ only differs slightly from the original Hamiltonian $H$, that
is, $\osc(H - \tilde H_\eps) = O(\eps)$. Because the perturbation is
small w.r.t. $\eps$, the $\PI$ and $\LSI$ constants of $\mu$ and
$\tilde\mu$ only differ up to an $\eps$-independent factor. This
observation is summarized in the following definition and subsequent
Lemma~\ref{locallemmodification}.
%
\begin{defn}[($\eps$-modification $\tilde H_\eps$ of $H$)]\label
{localdeftildeH}
The family of Hamiltonians $\{\tilde H_\eps\}_{\eps>0}$ is an $\eps
$-modification of $H$, if there exists an $\eps$-independent constant
$C_{\tilde H}>0$ such that for all $\eps$ small enough holds
{\renewcommand{\theequation}{{\~H}$_\eps$}
\begin{equation}
\label{localdeftildeHepsclose} \bigl\llvert \tilde H_\eps(x) - H(x)\bigr\rrvert \leq
C_{\tilde H}\eps \qquad\mbox{for all } x \in\Omega.
\end{equation}
\setcounter{equation}{0}}%
To each $\eps$-modification of $H$ we associate the family of
$\varepsilon$-modified Gibbs measures $\tilde\mu_\varepsilon$ by setting
\[
\label{localdefmodmu} {\tilde\mu}_\varepsilon(\dx{x}) := \frac{1}{Z_{{\tilde\mu
}_\varepsilon}} \exp
\bigl(-\eps^{-1}\tilde H_\eps(x) \bigr) \,\dx {x}\qquad
\mbox{with } Z_{\tilde\mu_\varepsilon} := \int\exp \bigl(-\eps^{-1} \tilde
H_\eps(x) \bigr) \,\dx{x} .
\]
\end{defn}
%
%
\begin{lem}[(Perturbation by an $\eps$-modification)]\label
{locallemmodification}
If the $\eps$-modified Gibbs measures $\tilde\mu_\varepsilon$
satisfy $\PI(\tilde\varrho)$ or $\LSI(\tilde\alpha)$, then the
measure $\mu$ also satisfies \ref{PIvarrho} or \ref{LSIalpha},
respectively, where the constants fulfill the estimate
\[
\label{localemodificationPI-LSI} \varrho\geq\exp (-2 C_{\tilde H} ) \tilde\varrho\quad
\mbox{and}\quad\alpha\geq\exp (-2 C_{\tilde H} ) \tilde \alpha,
\]
where $C_{\tilde H}$ is from \eqref{localdeftildeHepsclose}.
\end{lem}
\begin{pf}
The statement directly follows from an application of Theorem~\ref
{localthmHolleyStroock} by considering the estimate \eqref
{localdeftildeHepsclose}.
\end{pf}

Our approach to Theorem~\ref{thmlocalPI} consists of a nonstandard
application of a Lyapunov argument developed by Bakry, Barthe,
Cattiaux, Guillin, Wang and Wu (cf. \cite{Bakry2008,Bakry2008b,Wang2007} and \cite{Wu2008}), which is reminiscent of the spectral gap
characterization by Donsker and Varadhan \cite{Donsker1976}. Compared to these
works on the Lyapunov approach, we have to explicitly elaborate the
dependence of the PI and LSI constants on $\varepsilon$. Moreover, the
theory is only established for Gibbs measure on the whole space.
Therefore, the Lyapunov approach of the present work has two main ingredients:
\begin{itemize}
\item a Lyapunov function that has to satisfy Neumann boundary
conditions on $\Omega$ and certain estimates (cf. Definition~\ref
{localdefnlyap} and Theorem~\ref{localthmlyapLSI} below), and
\item a $\PI$ for a truncated Gibbs measure (cf. Definition~\ref
{localdefhatmu} and Lemma~\ref{locallemhatmuPI-LSI} below).
\end{itemize}
With the Lyapunov function, we are able to compare the scaling behavior
of the PI constant of $\mu$ with the behavior of the PI constant of
the truncated Gibbs measure $\hat\mu_a$ (cf. Theorem~\ref
{localthmlyapPI} and Theorem~\ref{localthmlyapLSI} below).
%
\begin{defn}[(Truncated Gibbs measure)]\label{localdefhatmu}
For a given number $a >0$, the truncated Gibbs measures $ \{\hat
\mu _{a,i} \}_{i=1}^M$ are obtained from the Gibbs measure $\mu
$ by
restriction to balls of radius $a \sqrt{\eps}$ around $ \{
m_i \}_{i=1}^M$, that is,
\begin{eqnarray}
\label{localedefhatmu}%
\hat\mu_{a,i}(\dx{x}) := \frac{\one_{B_{a\sqrt{\eps
}}(m_i)}(x)}{Z_{\hat\mu_{a,i}}}
\exp\bigl(- \eps^{-1} H(x)\bigr) \,\dx{x}
\nonumber\\
\eqntext{\displaystyle \mbox{with } Z_{\hat\mu_{a,i}} := \int_{B_{a\sqrt
{\eps}}(m_i)}\exp
\bigl(- \eps^{-1} H(x)\bigr) \,\dx{x}.}
\end{eqnarray}
\end{defn}
Because the domain and the Hamiltonian of the truncated Gibbs
measure $\hat\mu_{a,i}$ is convex, one can deduce the scaling
behavior of the truncated Gibbs measure $\hat\mu_{a,i}$ from the
Bakry--\'{E}mery criterion. More precisely, it holds the following.
%
\begin{lem}[($\PI$ and $\LSI$ for truncated Gibbs measure)]\label
{locallemhatmuPI-LSI}
For any $a>0$ and $i=1,\dots,M$ the measures $\hat\mu_{a,i}$ satisfy
$\PI(\hat\varrho)$ and $\LSI(\hat\alpha)$ for $\eps$ small
enough, where
\begin{equation}
\label{localehatmuPI-LSI}
\frac{1}{\hat\varrho} =O(\eps) \quad\mbox{and}\quad
\frac
{1}{\hat\alpha} = O(\eps) .
\end{equation}
\end{lem}
\begin{pf}
In the local minimum $m_i$ the Hessian of $H$ is nondegenerated by
Assumptions~\ref{assumeenv} or~\ref{assumeenvLSI}. Therefore,
for $\eps$ small enough, $H$ is strictly convex in $B_{a\sqrt{\eps
}}(m_i)$ and satisfies by the Bakry--\'{E}mery criterion
(cf. Theorem~\ref{localthmBakryEmery}) $\PI(\hat\varrho)$ and
$\LSI(\hat\alpha)$ with $\hat\varrho$ and $\hat\alpha$ obeying
the relation \eqref{localehatmuPI-LSI}.
\end{pf}

The standard ansatz $\exp ( \frac{H}{2\varepsilon}  )$
for a Lyapunov function has the nice feature that it automatically
satisfies Neumann boundary conditions on the basins of attraction
w.r.t. $H$, which would be also a canonical choice of the partition $\cP
_\cM$ (cf. Remark~\ref{rcanonicalpartition}). Unfortunately, one
cannot guarantee that the necessary estimates for $\exp ( \frac
{H}{2\varepsilon}  )$ hold because there is no control on the
sign of $\Delta H(x)$ close to\vadjust{\goodbreak} saddles [see \eqref{locallyapcondexp} below].
We circumvent this technical problem in the following way: By the
observation from above it suffices to consider an $\eps$-modification
$\tilde H_\eps$ of $H$. We explicitly construct an $\eps
$-modification $\tilde H_\eps$ on the whole space $\R^n$ with the
property that the standard ansatz $\exp ( \frac{\tilde H_\eps
}{2\varepsilon}  )$ satisfies the necessary estimates for being
a Lyapunov function. However in general, the function $\exp (
\frac{\tilde H_\eps}{2\varepsilon}  )$ does not satisfy
Neumann boundary conditions on the basins of attraction w.r.t. $H$.
This problem is solved by the following two observations.
\begin{itemize}
\item The first one is that $\exp ( \frac{\tilde H_\eps
}{2\varepsilon}  )$ satisfies Neumann boundary conditions on the
basin of attraction w.r.t. the deterministic gradient flow defined by
$\tilde H_\eps$, that is,
\begin{equation}
\label{erightchoiceofanpartition} \Omega_i := \Bigl\{y \in\R^n \dvtx  \lim
_{t\to\infty} y_t = m_i, \dot
{y}_t=-\nabla\tilde H_\eps(y_t),
y_0 =y \Bigr\}.
\end{equation}
\item The second observation is that this partition $ \{\Omega
_i \}_{i=1}^M$ of $\R^n$ is admissible in the sense of
Definition~\ref
{defpartition} (see Lemma~\ref{pLyapunovverification} below). This
fact is intuitively clear from the fact that $\tilde H_\eps$ is only a
small perturbation of $H$.
\end{itemize}
Hence, we choose the partition $\cP_\cM:= \{\Omega_i \}
_{i=1}^M$ of
$\R^n$ according to \eqref{erightchoiceofanpartition} and apply
the Lyapunov approach to the local Gibbs measures $\tilde\mu_{\eps
,i}$ given by
\begin{eqnarray}
\label{elocalepsmodlocalGibbsmeasure} {\tilde\mu}_{\varepsilon,i} (\dx{x}) := \frac{\one_{\Omega
_i}(x)}{Z_{{\tilde\mu}_{\varepsilon,i}}} \exp
\bigl(-\eps ^{-1}\tilde H_\eps(x) \bigr) \,\dx{x}\nonumber
\\[-8pt]
\\[-8pt]
\eqntext{\displaystyle \mbox{with } Z_{\tilde\mu_{\varepsilon,i}} := \int_{\Omega_i} \exp \bigl(-
\eps^{-1} \tilde H_\eps(x) \bigr) \,\dx{x} . }
\end{eqnarray}
We get that the local Gibbs measures $\tilde\mu_{\eps,i}$ satisfy a
local PI and LSI with the desired scaling behavior in $\eps$. This
scaling behavior of the PI and LSI constant is then transferred to the
original Gibbs measure $\mu$ restricted to the sets $\Omega_i$ by
using the perturbation Lemma~\ref{locallemmodification}.

The remaining part of this section is organized in the following way.
\begin{itemize}
\item In Section~\ref{slyap}, we present the abstract framework
how the Lyapunov approach is used for deriving the local PI. We
additionally motivate the perturbative nature of the construction of
the Lyapunov function. Under the assumption of the existence of a
Lyapunov function, we also state the proof Theorem~\ref{thmlocalPI}.
\item In Section~\ref{sconstructionlyapunov}, we provide the
central ingredient for the Lyapunov approach, namely the existence of a
Lyapunov function. We also show that the partition obtained by \eqref
{erightchoiceofanpartition} is admissible.
\item In Section~\ref{slyapunovLSI}, we present the abstract
framework how the Lyapunov approach is used for deriving the local LSI.
We show that one can use the same Lyapunov function for the local PI as
for the local LSI. We also state the proof of Theorem~\ref
{thmlocalLSI} deducing the local LSI.
\end{itemize}

\subsection{Lyapunov approach for the Poincar\'{e} inequality}\label{slyap}

We start with explaining the Lyapunov approach for deducing a PI. The
central notion for the Lyapunov approach is the following definition.

%
\begin{defn}[(Lyapunov function for Poincar\'{e} inequality)]\label
{localdefnlyap}
Let $H\dvtx  \Omega\to\R$ be a Hamiltonian with Gibbs measure $\mu(\dx
{x}) = \one_\Omega(x) Z_\mu^{-1} \exp (-\eps^{-1} H(x)
)\,\dx{x}$. Then $W\dvtx  \Omega\to[1,\infty)$ is a \emph{Lyapunov
function} for $H$ provided that:
\begin{longlist}[(ii)]
\item There exist a domain $U\subset\Omega$ and constants $b>0$ and
$\lambda>0$ such that
\begin{equation}
\label{locallyapcond} \varepsilon^{-1} L W \leq- \lambda W + b
\one_{U} \qquad\mbox {a.e. in $\Omega$}.
\end{equation}
\item$W$ satisfies Neumann boundary conditions on $\Omega$ such that
the integration by parts formula holds
\begin{equation}
\label{elyapintegrationbyparts} \forall f\in H^1(\mu|_\Omega)\dvtx  \int
_\Omega f (-L W) \,\dx{\mu} = \eps\int_\Omega
\langle\nabla f , \nabla W \rangle \,\dx{\mu} .
\end{equation}
\end{longlist}
\end{defn}
Compared to the Lyapunov function of \cite{Bakry2008} the condition
$(ii)$ in Definition~\ref{localdefnlyap} is new. The reason is that
we work on the domain $\Omega$ and not on the whole space $\R^n$. The
next statement shows that a Lyapunov function and a $\PI$ for the
truncated measure can be combined to get a $\PI$ for the whole measure.
%
\begin{theo}[(Lyapunov condition for $\PI$ on domains $\Omega$)] \label
{localthmlyapPI}
Suppose that $H$ has a Lyapunov functions in the sense of
Definition~\ref{localdefnlyap} and that the restricted measure $\mu
_U$ given by
\[
\mu_U (\dx{x}) := \mu(\dx{x})\llcorner{U} = \frac{\one_{U}
(x)}{\mu(U)}
\mu(\dx{x}) ,
\]
satisfies $\PI(\varrho_U)$. Then the associated Gibbs measure $\mu$
also satisfies \ref{PIvarrho} with constant
\[
\varrho\geq\frac{\lambda}{b + \varrho_U} \varrho_U .
\]
\end{theo}
The content of the last theorem is standard (cf. \cite{Bakry2008}),
except that we work on the domain $\Omega$ and not on the whole space
$\R^n$. For the convenience of the reader, we state the short proof.
\begin{pf*}{Proof of Theorem~\ref{localthmlyapPI}}
Let us rewrite the Lyapunov condition \eqref{locallyapcond} and observe
\begin{equation}
\label{lyapPIp1} 1\leq- \frac{LW}{\eps\lambda W}+\frac{b}{\lambda}
\frac{\dsOne
_{U}}{W} \leq- \frac{LW}{\eps\lambda W}+\frac{b}{\lambda} \dsOne
_{U} ,
\end{equation}
since $W\geq1$ by Definition~\ref{localdefnlyap}.
By the integration by parts rule \eqref{elyapintegrationbyparts},
we obtain following estimate
which is due to Definition~\ref{localdefnlyap}(ii). Therewith, we
deduce the estimate
\begin{eqnarray}
\label{lyapPIp2} %
\int f^2 \frac{(-LW)}{\eps W}\,\dx{\mu}
&=& \int \biggl\langle\nabla \biggl(\frac{f^2}{W} \biggr) , \nabla W \biggr
\rangle \,\dx{\mu}
\nonumber
\\
&=& 2 \int\frac{f}{W} \langle\nabla f , \nabla W \rangle \,\dx{\mu} -
\int \frac{f^2 \llvert \nabla W\rrvert ^2 }{W^2}\,\dx{\mu}
\nonumber
\\[-8pt]
\\[-8pt]
&=& \int\llvert \nabla f\rrvert ^2 \,\dx{\mu} - \int\biggl\llvert
\nabla f - \frac {f}{W} \nabla W\biggr\rrvert ^2 \,\dx{\mu }
\nonumber
\\
&\leq& \int\llvert \nabla f\rrvert ^2 \,\dx{\mu} . %
\nonumber
\end{eqnarray}
Let us now turn this estimate into one for the variance $\var_{\mu
}(f)$. Due to fundamental properties of the variance, it holds $\var
_{\mu}(f) \leq\int (f - m )^2\,\dx{\mu}$, for any $m\in
\R
$. Hence, applying the estimates \eqref{lyapPIp1} and \eqref
{lyapPIp2} yields
\begin{eqnarray}
\label{lyapPIp3} %
\var_{\mu}(f) & \leq& \int (f-m
)^2 \,\dx{\mu} \stackrel {\mathclap{\eqref{lyapPIp1}}} {\leq} \int
(f-m )^2 \frac
{(-LW)}{\eps\lambda W} + \frac{b}{\lambda} \int
_{U} (f-m )^2 \,\dx{\mu}
\nonumber
\\[-8pt]
\\[-8pt]
&\stackrel{\mathclap{\eqref{lyapPIp2}}} { \leq }& \frac
{1}{\lambda} \int
\llvert \nabla f\rrvert ^2 \,\dx{\mu} + \frac{b \mu
(U)}{\lambda} \int (f-m
)^2 \,\dx{\mu_U} . %
\nonumber
\end{eqnarray}
We set $m=\int f\,\dx{\mu_U}$, then the last integral in the right-hand
side of \eqref{lyapPIp3} becomes $\var_{\mu_U}(f)$, to which we
apply the assumption $\PI(\varrho_U)$.
\end{pf*}
Considering the last theorem, it is only left to construct a Lyapunov
function in the sense of Definition~\ref{localdefnlyap} in order to
deduce the local PI of Theorem~\ref{thmlocalPI}. An ansatz (cf.
\cite{Bakry2008}) for a Lyapunov function is the function $W =\exp
(\frac{1}{2\eps}H )$. Why is this in general a good
candidate for an Lyapunov function?

First note that because by our Assumptions~\ref{assumeenv} or~\ref
{assumeenvLSI} it holds $H \geq0$ hence $W \geq1$ as desired. The
second reason is that this choice satisfies Neumann boundary conditions
on the boundary of the basin of attraction $\Omega$ (see Theorem~\ref
{thmNeumann}).

The third reason is that for this choice of $W$ the Lyapunov
condition \eqref{locallyapcond} is already almost satisfied. One
only has to have a special look at critical points. To be more precise,
let us consider the condition \eqref{locallyapcond} which becomes
\begin{equation}
\label{locallyapcondexp} %
\frac{LW}{\eps W} = \frac{1}{2\eps}\Delta H(x)
- \frac{1}{4\eps
^2}\bigl\llvert \nabla H(x)\bigr\rrvert ^2
\stackrel{!} {\leq} -\lambda + b \dsOne_{U}(x). %
\end{equation}
We investigate under which circumstances this condition is satisfied:
\begin{itemize}
\item At infinity: The assumption \eqref{assumegradlaplace} ensures
that \eqref{locallyapcondexp} is satisfied outside of a fixed large
ball $B_{\tilde R}(0)$ [cf. \eqref{elyapunovoutsidecritical} below].
\item Away from critical points: The Morse assumption ensures $H$ to be
quadratic around critical points, that is, there exists a global
constant $c_H>0$ such that $\llvert \nabla H(x)\rrvert \geq
c_H \dist(x,\cS)$
in a neighborhoods of critical points $\cS$. This estimate
yields \eqref{locallyapcondexp} for $x$ outside of neighborhoods of
order $\sqrt{\eps}$ around critical points (see proof of Lemma~\ref
{plyapunovPIoutside} below).
\end{itemize}
The gradient term cannot help to establish the estimate \eqref
{locallyapcondexp}, if one is close to critical points. More
precisely, it holds:
\begin{itemize}
\item If $x$ is in an $\sqrt{\eps}$-neighborhood around the
minimum $0$, then $\Delta H(x) \approx\sum_i \lambda_i >0$, where
$ \{\lambda_i \}_{i=1}^n$ are the eigenvalues of the
Hessian at $0$.
Additionally, the gradient can be estimated as $\llvert \nabla
H(x)\rrvert ^2
\gtrsim\lambda_{\min}^2 \llvert  x\rrvert ^2$, where
$\lambda_{\min} =\min_i \lambda_i$. Hence, one cannot compensate the positive Laplacian by
the gradient of $H$. Therefore, one has to choose $U = B_{a\sqrt{\eps
}}(0)$ to guarantee the Lyapunov condition \eqref{locallyapcondexp}
around the minimum at $0$.
\item If $x$ is close a local maximum, the Laplacian $\Delta H(x)$ is
negative. Hence, the Lyapunov condition is \eqref{locallyapcondexp}
is satisfied in this region.
\item Assume that $x$ is in an $\sqrt{\eps}$-neighborhood around a
saddle, that is, a critical point $s\in\cS$ of order $1 \leq k < n$.
Again, the gradient term cannot help to establish the estimate \eqref
{locallyapcondexp}. Hence, the condition \eqref{locallyapcondexp} becomes
\[
\Delta H(x) \approx\lambda_1^- +\cdots+ \lambda_k^- +
\lambda _{k+1}^+ + \cdots+ \lambda_{n}^+ \stackrel{!} {\leq}
-\lambda,
\]
where $\lambda^-_i$ are the negative eigenvalue of the Hessian at $s$
and $\lambda^+_j$ are the positive eigenvalues of the Hessian at $s$.
However, for a general Hamiltonian $H$ it may hold that
\[
\lambda_1^- +\cdots+ \lambda_k^- +
\lambda_{k+1}^+ + \cdots+ \lambda _{n}^+ \geq0
\]
implying that $W =\exp (\frac{1}{2\eps}H )$ is not always a
Lyapunov function.
\end{itemize}
Nevertheless, these observations show that $W =\exp (\frac
{1}{2\eps}H )$ is a pretty good guess for a Lyapunov function:
One only
has to change $W$ close to saddles of $H$. This leads to the following
strategy (cf. Lemma~\ref{pLyapunovinside} from below):
\begin{itemize}
\item We construct a perturbation $\tilde H_\eps$ of the Hamiltonian
$H$, which coincides with $H$ except of $\sqrt{\eps}$-neighborhoods
around saddles.
\item In a $\sqrt{\eps}$-neighborhood around a saddle, the
perturbation is constructed in such a way that on the one hand the
Laplacian of $\tilde H_\eps$ is strictly negative. This implies that
the function $W =\exp (\frac{1}{2\eps} \tilde H_\eps )$
satisfies the estimate \eqref{locallyapcondexp}, which is necessary
for being a Lyapunov function.
\item To assure that $W =\exp (\frac{1}{2\eps} \tilde H_\eps
)$ satisfies Neumann boundary condition, we choose $\Omega$ as
a basin of attraction w.r.t. the gradient flow of $\tilde H_\eps$
[cf. \eqref{erightchoiceofanpartition}].
\end{itemize}
%
After these considerations, let us summarize how the Lyapunov approach
is used.
%
\begin{prop}\label{plyapapplied}
Assume that an Hamiltonian $\tilde H_\eps$ satisfies the
Assumption~\ref{assumeenv} uniformly in $\eps$. Let $\mathcal{M}=
\{ m_1,\ldots,m_M  \}$ denote the local minima of $\tilde
H_\eps$.
Assume that there are constants $a>0$ and $\lambda_0>0$ such that for
all $\eps>0$ small enough holds
\begin{equation}
\label{elyapverificationdrif} \qquad \frac{1}{2\varepsilon} \Delta\tilde H_\eps(x) -
\frac{1}{4
\varepsilon^2} \bigl\llvert \nabla\tilde H_\eps(x)\bigr\rrvert
^2 \leq- \frac{\lambda
_0}{\varepsilon} \qquad\mbox{for all } x \notin\bigcup
_{m
\in\mathcal{M} } B_{a \sqrt{\varepsilon}} (m).
\end{equation}
Consider the partition $\cP_\cM= \{\Omega_i \}_{i=1}^M$
into the
basins of attraction of the gradient flow of $\tilde H_\eps$
[cf. \eqref{erightchoiceofanpartition}].
Then the associated local Gibbs measures $ \{\tilde\mu_{\eps
,i} \}_{i=1}^M$ given by \eqref{elocalepsmodlocalGibbsmeasure}
satisfy $\PI(\tilde\varrho_i)$ with constant
\[
\tilde\varrho_i^{-1} = O(\eps) .
\]
\end{prop}
\begin{pf}
%
The function $W=\exp (\frac{1}{2\eps} \tilde H_\eps )$
satisfies Neumann boundary conditions on each domain of attraction
$\Omega_i$ in the sense of \eqref{elyapintegrationbyparts} by
Theorem~\ref{thmNeumann}. Indeed, the gradient of $W$ is
\[
\label{egradientlyapunovfunction} \nabla W = \frac{1}{2 \eps} (\nabla\tilde H_\eps)
\exp \biggl(\frac
{1}{2\eps}\tilde H_\eps \biggr).
\]
Hence, $\nabla W \parallel \nabla\tilde H_\eps$ everywhere. Moreover,
$\tilde H_\eps\in C^3$ is Morse and proper by Assumption~\ref
{assumeenv}, which shows all the assumptions of Theorem~\ref{thmNeumann}.

Let $\Omega_i$ be fixed. Then the estimate \eqref
{elyapverificationdrif} is just a translation of the estimate \eqref
{locallyapcond} with constants $\lambda= \frac{\lambda_0}{\eps}$
and $b=\frac{b_0}{\eps}$ for some $b_0>0$. Moreover, we choose
$U=B_{a\sqrt{\eps}}(m_i)$. Therefore, the function $W$ is a Lyapunov
function in the sense of Definition~\ref{localdefnlyap} on $\Omega
_i$. Theorem~\ref{localthmlyapPI} yields that the measure $\tilde
\mu_{\eps,i}$ satisfies $\PI(\tilde\varrho_i)$ with
\[
\label{localproofPI2} \tilde\varrho_i \geq\frac{\lambda_0 \hat\varrho}{b_0 + \eps\hat
\varrho},
\]
where $\hat\varrho$ denotes the PI constant of the truncated Gibbs
measure $\hat\mu_{a,i}$ from Definition~\ref{localdefhatmu}. By
Lemma~\ref{locallemhatmuPI-LSI} holds $\hat\varrho^{-1} = O(\eps
)$, which yields $\tilde\varrho_i^{-1}=O(\eps)$.
\end{pf}

Following our strategy, the main ingredient of the proof of the local
PI is the existence of an $\eps$-modified Hamiltonian $\tilde H_\eps$
satisfying assumption \eqref{elyapverificationdrif} of
Proposition~\ref{plyapapplied}. 
%
\begin{lem}[(Lyapunov function for $\PI$)]\label{pLyapunovverification}
There exits an $\varepsilon$-modification $\tilde H_\varepsilon$
of $H$ in the sense of Definition~\ref{localdeftildeH} such that the
Lyapunov estimate \eqref{elyapverificationdrif} holds for $\tilde
H_\varepsilon$. The corresponding partition $\cP_\cM= \{\Omega
_i \}_{i=1}^M$ into the basins of attraction of the gradient flow of
$\tilde H_\eps$ [cf. \eqref{erightchoiceofanpartition}] is
admissible in the sense of Definition~\ref{defpartition}.
\end{lem}
The proof of Lemma~\ref{pLyapunovverification} is not complicated
but a bit lengthy. It is stated in full detail in Section~\ref{sconstructionlyapunov}. Now, we only have to put together the parts
in order to proof the first main result Theorem~\ref{thmlocalPI}.
\begin{pf*}{Proof of Theorem~\ref{thmlocalPI}}
By a combination of Lemma~\ref{plyapapplied} and Lem-\break ma~\ref
{pLyapunovverification} we know that the $\eps$-modified Gibbs
measures $\tilde\mu_{\eps,i}$ restricted to $\Omega_i$ satisfy a PI
with the desired scaling behavior $\tilde\varrho_i^{-1}=O(\eps)$.
Lemma~\ref{locallemmodification} implies that then the unmodified
Gibbs measure $\mu_i$ restricted to $\Omega_i$ also satisfies a PI
with the same scaling behavior $\varrho_i^{-1}=O(\eps)$.
\end{pf*}

\subsection{Construction of a Lyapunov function}\label
{sconstructionlyapunov}

This section is devoted to the proof of Lemma~\ref
{pLyapunovverification}. We have to construct an $\eps$-modified
Hamiltonian $\tilde H_\eps$ that satisfies the estimate \eqref
{elyapverificationdrif}. Following the motivation of Section~\ref{slyap}, we set $\tilde H_\eps= H$ away from critical points.
Therefore, we have to show that $H$ satisfies the estimate \eqref
{elyapverificationdrif} away from critical points, which is the
content of the next statement.
%
\begin{lem}\label{plyapunovPIoutside}
Assume that the Hamiltonian $H$ satisfies the Assumption~\ref
{assumeenv}. Recall that $\mathcal{S}$ denotes the set of all
critical points of $H$ in $\Omega$; that is,
\[
\mathcal{S}= \bigl\{y \in\Omega\mid \nabla H (y)=0 \bigr\}.
\]
Then for $a>0$ large enough exists $\lambda_0>0$ and $\eps_0>0$ such
that for all $\eps< \eps_0$
\begin{equation}
\label{elyapunovoutsidecritical}\quad  \frac{\Delta H(x)}{2\varepsilon} - \frac{|\nabla H(x)
|^2}{4\varepsilon^2} \leq-
\frac{\lambda_0}{\varepsilon} \qquad\mbox{for all } x \in\R^n \Bigm\backslash\bigcup
_{y \in\mathcal
{S}} B_{ a \sqrt{\varepsilon}} (y) .
\end{equation}
\end{lem}
\begin{pf}
The proof basically consists only of elementary calculations based on
the nondegeneracy assumption on $H$. We consider two cases: One in
which we verify \eqref{elyapunovoutsidecritical} for $\llvert x\rrvert  \geq
\tilde R$ with $\tilde R < \infty$ large enough. In the second case,
we verify   \eqref{elyapunovoutsidecritical} for $ \llvert x\rrvert  \leq
\tilde R$.

Let us turn to the first case. We use the assumptions \eqref
{assumegradsuperlinear} and \eqref{assumegradlaplace} and we define
$\tilde R$ such that
\begin{equation}
\label{localdeftildeR} \forall\llvert x\rrvert \geq\tilde R\dvtx  \llvert \nabla H\rrvert
\geq\frac
{C_H}{2} \quad\mbox{and}\quad\llvert \nabla H\rrvert ^2
- \laplace H(x) \geq - 2 K_H.
\end{equation}
Therewith, it is easy to show that for $\llvert  x\rrvert \geq
\tilde R$
\begin{eqnarray}
\label{egradientconditionfaroutside}%
\frac{\Delta H(x)}{2\varepsilon} - \frac{|\nabla H(x)
|^2}{4\varepsilon^2} & \stackrel{
\mathclap{\eqref {localdeftildeR}}} {\leq}& \frac{1}{\eps}
\biggl(K_H - \frac
{\llvert \nabla H(x)\rrvert }{2} \biggl(\frac{\nabla
H(x)}{2\eps}-1 \biggr)
\biggr)
\nonumber
\\[-8pt]
\\[-8pt]
& \leq& \frac{1}{\eps} \biggl(K_H-\frac {C_H^2}{8} \biggl(
\frac{C_H^2}{8\eps}-1 \biggr) \biggr) \leq- \frac{\lambda
_0}{\eps},
\nonumber
\end{eqnarray}
if $\eps\leq\frac{C_H^2}{8}  (1+8/C_H^2(K_H+\lambda _0)
)^{-1}=:\eps_0$. The latter shows the desired statement in this
case, with $\lambda_0>0$ arbitrary for $\eps\leq\eps_0(\lambda
_0)$.

Let us consider the second case. Because $\llvert x\rrvert \leq\tilde R$ it holds
$\llvert \Delta H (x)\rrvert  \leq C_{\tilde R}$. Therefore, the desired
estimate \eqref{elyapunovoutsidecritical} follows, if we show that
there is a constant $0< c_H$ such that
\begin{equation}
\label{egradientconditionawaycriticalpoints}\quad    \bigl\llvert \nabla H(x)\bigr\rrvert \geq c_H a
\sqrt{\eps} \qquad\forall x \in B_{\tilde
R}(0) \Bigm\backslash\bigcup
_{y \in\mathcal{S}} B_{ a \sqrt{\varepsilon
}} (y) \mbox{ and } \forall a\in\bigl[0,
\eps^{-1/2}\bigr].\hspace*{-25pt}
\end{equation}
Because, then it follows
\[
\frac{\Delta H(x)}{2\varepsilon} - \frac{|\nabla H(x)
|^2}{4\varepsilon^2} \stackrel{\mathclap{\eqref
{egradientconditionawaycriticalpoints}}}
{\leq} \frac{1}{\eps
} \biggl(\frac{C_{\tilde R}}{2} - \frac{c_H a}{4}
\biggr) =: - \frac
{\lambda_0}{\eps},
\]
with $\lambda_0>0$ by choosing $a> 2C_{\tilde R}/c_H =: a_0$. Hence,
we can choose first $a>a_0$, which gives rise to some $\lambda
_0(a)>0$, by the last estimate under the assumption $a<\eps
_0^{-1/2}\leq\eps^{-1/2}$. Hence, we have to choose $\eps_0 < \min
\{\eps_0(\lambda_0(a)), a^{-2} \}$ with $\eps_0(\lambda_0(a))$
defined after \eqref{egradientconditionfaroutside}.

Finally, the estimate \eqref
{egradientconditionawaycriticalpoints} is a consequence of the
fact that $H$ is a Morse function (cp. Definition~\ref{defmorse} and
Assumption~\ref{assumeenv}) and, therefore, nondegenerate quadratic
around critical points. That means, there exists a global constant
$c_H>0$ such that $\llvert \nabla H(x)\rrvert \geq c_H \min
\{\dist(x,\cS ),1 \}$, which implies \eqref
{egradientconditionawaycriticalpoints}.
\end{pf}
Now, we consider the $\eps$-modification $\tilde H_\eps$ near
critical points. The verification of the following statement represents
the core of the construction of the Lyapunov function.
%
\begin{lem}\label{pLyapunovinside}
Let $\mathcal{M}= \{m_1, \ldots, m_M  \}$ denote the set
containing the minima of $H$. Then there are constants $C>0$, $a>0$ and
$\lambda_0 > 0 $ such that for $\eps<C$ there exists an $\eps
$-modification $\tilde H_\eps$ of $H$ in the sense of Definition~\ref
{localdeftildeH} satisfying
\[
\label{elyapunovlocalpert} \tilde H_\eps(x)= H (x) \qquad\mbox{for all } x \notin
\bigcup_{y
\in\mathcal{S}\setminus\mathcal{M} } B_{a \sqrt{\varepsilon}} (y)\vadjust{\goodbreak}
\]
and
\begin{equation}
\label{elyapunovinsidecritical} \frac{\Delta\tilde H_\varepsilon(x)}{2\varepsilon} - \frac{|\nabla
\tilde H_\varepsilon(x) |^2}{4\varepsilon^2} \leq-
\frac{\lambda
_0}{\varepsilon} \qquad\mbox{for all } x \in\bigcup
_{y \in
\mathcal{S}\setminus\mathcal{M}} B_{a \sqrt{\varepsilon}} (y) .
\end{equation}
As a direct consequence of Lemma~\ref{plyapunovPIoutside}, the
estimate \eqref{elyapunovinsidecritical} is satisfied for all
\[
x \notin\bigcup_{m \in\mathcal{M} } B_{a \sqrt{\varepsilon}} (m).
\]
\end{lem}
\begin{pf}
It is sufficient to construct the $\varepsilon$-modification $\tilde
H_\varepsilon$ only locally on a small neighborhood of any critical
point $y\in\cS\setminus\mathcal{M}$. By translation, we may assume
w.l.o.g. that $y=0$.

Because the Hamiltonian $H$ is a Morse function in the sense of
Definition~\ref{defmorse}, we may assume that $u_i$, $i \in \{
1, \ldots, n  \}$ are orthonormal eigenvectors w.r.t. the
Hessian $\nabla^2 H (0)$. The corresponding eigenvalues are denoted by
$\lambda_i$, $i \in \{1, \ldots, n  \}$ labeled such that
$\lambda_1, \ldots, \lambda_\ell< 0$ and $\lambda_{\ell+1},
\ldots, \lambda_n > 0$ for some $\ell\in \{1, \ldots, n
\}$.
If $\ell=n$, hence $\lambda_i < 0 $ for $i=1,\dots,n$, we are nearby
a local maximum and set $\tilde H_\varepsilon(x) = H(x)$ on $B_{a
\sqrt{\varepsilon}}(0)$ and the desired estimate \eqref
{elyapunovinsidecritical} follows directly for $x \in B_{a \sqrt
{\varepsilon}}(0)$.

Otherwise, that is, $\ell<n$, let us choose a constant $\delta>0$
small enough such that
\begin{equation}
\label{elyapPIdeftildedelta}\quad  -\tilde\delta:= (n-2\ell) \delta+ \sum
_{i=1}^\ell\lambda_i < 0 \quad
\mbox{and}\quad\delta\leq\frac{1}{2}\min \{\lambda _i \dvtx  i=
\ell+1,\dots,n \} .\hspace*{-25pt}
\end{equation}
Because $u_1, \ldots, u_n $ is an orthonormal basis of $\R^n$, we
introduce a norm $\llvert \cdot\rrvert _\delta$ on $\R^n$ by
\begin{equation}
\label{elyapdefnorm} \llvert x\rrvert _\delta^2:= \sum
_{i=1}^\ell\frac{1}{2} \delta\bigl\llvert
\langle u_i , x \rangle\bigr\rrvert ^2 + \sum
_{i=\ell+1}^n \frac{1}{2} (\lambda_i
- \delta) \bigl\llvert \langle u_i , x \rangle \bigr\rrvert
^2.
\end{equation}
The norm $\llvert \cdot\rrvert _\delta$ is equivalent to the standard
Euclidean norm $\llvert \cdot\rrvert $ and satisfies the estimate
\begin{equation}
\label{elyapequivnorms} \frac{\delta}{2} \llvert x\rrvert ^2 \leq \llvert
x\rrvert _\delta^2 \leq \frac{\lambda_{\max}^+ - \delta}{2} \llvert x\rrvert
^2 \leq \frac
{\lambda_{\max}^+}{2} \llvert x\rrvert ^2 ,
\end{equation}
where $\lambda_{\max}^+ = \max \{\lambda_i \dvtx  i=\ell+1,\dots
,n \}$. The last ingredient for the construction of $\tilde
H_\varepsilon
$ is a smooth cut-off function $\xi\dvtx  [0, \infty) \to\R$ satisfying
for $a>0$ to be specified later
\begin{eqnarray}
\label{edefauxxi} \qquad \xi'(r) &=& -1 \qquad \mbox{for } r\leq
\tfrac{1}{4}a^2 \eps, \qquad -1\leq \xi'(r) \leq0\qquad
\mbox{for } r \geq\tfrac{1}{4}a^2 {\eps}, \nonumber
\\[-8pt]
\\[-8pt]  \xi(r) &=& 0\qquad
\mbox{for } r\geq a^2 {\eps}
\nonumber
\end{eqnarray}
and in addition for some $C_\xi>0$,
\begin{equation}
\label{epropauxxipert} 0\leq\xi(r)\leq C_\xi a^2 \eps\quad
\mbox{and}\quad\bigl\llvert \xi ''(r)\bigr\rrvert \leq
\frac{C_\xi}{{a^2 \eps}} .
\end{equation}
With the help of the norm $\llvert \cdot\rrvert _\delta$ and the function $\xi$ we
define the function $\tilde H_\eps$ by
\begin{equation}
\label{elyapPIdefHb} \tilde H_\varepsilon(x) = H(x) + H_b (x) \qquad
\mbox{where } H_b (x) := \xi \bigl( \llvert x\rrvert
_\delta^2 \bigr) .
\end{equation}
%
Note that by definition of $H_b$ holds $\tilde H_\eps(x)= H(x)$ for
all $\llvert x\rrvert  \geq a \sqrt\varepsilon$. Because $\xi(r)=O(\eps)$, it
follows that $\tilde H_\eps$ is an $\eps$-modification of $H$ in the
sense of Definition~\ref{localdeftildeH}.

Let us now turn to the verification of the estimate \eqref
{elyapunovinsidecritical}. It is sufficient to deduce the following
two facts: The first one is the estimate
\begin{equation}
\label{elyapPIconstr1} \Delta\tilde H_\eps(x) \leq- \frac{\tilde\delta}{2}
\qquad\mbox {for all } \llvert x\rrvert _\delta\leq\frac{a}{2}
\sqrt{\varepsilon}.
\end{equation}
The second one is that there is a constant $\lambda_0>0$ such that for
$a$ large enough and $\varepsilon$ small enough it holds
\begin{equation}
\label{elyapPIconstr2} \frac{\Delta\tilde H_\eps(x)}{2} - \frac{|\nabla\tilde H_\eps
(x)|^2}{4\varepsilon} \leq-
\lambda_0 \qquad\mbox{for all } \frac{a}{2} \sqrt{\varepsilon}
\leq\llvert x\rrvert _\delta\leq a \sqrt {\varepsilon} .
\end{equation}
Let us first derive the estimate \eqref{elyapPIconstr1}. Using
that $\xi'(x)= -1$ for $\llvert x\rrvert _\delta\leq\frac{a}{2} \sqrt\varepsilon
$, one obtains that $\Delta H_b(x) = - \Delta\llvert x\rrvert _\delta^2$ for
$\llvert x\rrvert _\delta\leq\frac{a}{2} \sqrt\varepsilon$. Hence, by Taylor
expansion we get for $\llvert x\rrvert _\delta\leq\frac{a}{2} \sqrt\varepsilon$ that
\begin{eqnarray*}
\Delta\tilde H (x) & =& \Delta H(0) - \Delta\llvert x\rrvert
_\delta^2 + O(\sqrt{\eps}) \leq\sum
_{i=1}^n \lambda_i -\sum
_{i=\ell+1}^n \lambda_i + (n-2\ell) \delta+
O(\sqrt{\eps})
\\
&=& \sum_{i=1}^\ell\lambda_i
+ (n-2\ell) \delta+ O(\sqrt{\eps}) \overset{\mathclap{\eqref{elyapPIdeftildedelta}}}
{\leq} - \tilde\delta+ O(\sqrt\varepsilon) \leq - \frac{\tilde\delta}{2},
\end{eqnarray*}
for $\eps$ small enough, which yields the desired statement \eqref
{elyapPIconstr1}.

Let us turn to the verification of \eqref{elyapPIconstr2}. We need
that there exists a constant $0<C_{\Delta}<\infty$ independent of
$\eps$ and $a$ such that
\begin{equation}
\label{eboundCDelta} \Delta\tilde H (x) \leq C_\Delta\qquad\mbox{for all }
\frac{a}{2} \sqrt{\varepsilon} < \llvert x\rrvert _\delta< a
\sqrt{\varepsilon}.
\end{equation}
Indeed, observe that
\begin{eqnarray}
\Delta\tilde H_\eps(x) &=& \Delta H(x) +
\xi''\bigl(\llvert x\rrvert _\delta^2
\bigr) \bigl\llvert \nabla \llvert x\rrvert _\delta^2\bigr
\rrvert ^2 +\underbrace{\xi'\bigl(\llvert x\rrvert
_\delta^2\bigr)}_{\leq0} \underbrace{\Delta \llvert
x\rrvert _{\delta
}^2}_{\geq0}
\nonumber
\\
&\overset{\mathclap{\eqref{epropauxxipert}}}
{\leq}& \Delta H(x) + \frac{C_\xi}{{a^2 \eps}} \Biggl\llvert \sum
_{i=1}^\ell\delta \langle u_i , x
\rangle u_i + \sum_{i=\ell+1}^n (
\lambda_i-\delta) \langle u_i , x \rangle
u_i \Biggr\rrvert ^2
\nonumber
\\
&\leq& \Delta H(x) + \frac{C_\xi}{{a^2 \eps}} \Biggl(\sum
_{i=1}^\ell\delta^2 \bigl\llvert \langle
u_i , x \rangle \bigr\rrvert ^2 + \sum
_{i=\ell+1}^n (\lambda_i-
\delta)^2 \bigl\llvert \langle u_i , x \rangle\bigr
\rrvert ^2 \Biggr)
\nonumber
\\
&\overset{\mathclap{\eqref{elyapdefnorm}}}
{\leq}& \Delta H(x) + \frac{C_\xi}{a^2 {\eps}} 2\lambda_{\max}^+ \llvert x
\rrvert _\delta^2 
\leq C_H +
2C_\xi\lambda_{\max}^+ =: C_\Delta,\nonumber
\end{eqnarray}
where $C_\Delta$ is independent of $\eps$ and $a$, which
yields \eqref{eboundCDelta}.

Additionally, we need that there is a constant $0 < c_\nabla< \infty$
such that
\begin{equation}
\label{eboundCnabla} \bigl\llvert \nabla\tilde H_\eps(x)\bigr\rrvert
^2 \geq c_\nabla a^2 \varepsilon\qquad \mbox{for
all } \frac{a}{2} \sqrt{\varepsilon} < \llvert x\rrvert _\delta<
\tilde a \sqrt{\varepsilon}.
\end{equation}
Before deducing \eqref{eboundCnabla}, we want to show that the
observations \eqref{eboundCDelta} and \eqref{eboundCnabla}
already yield the desired statement \eqref{elyapPIconstr2}: For
$a^2 \geq4 \frac{C_\Delta}{c_\nabla}$, one gets
\begin{eqnarray*}
\frac{\Delta\tilde H_\eps(x)}{2 \varepsilon} - \frac{|\nabla
\tilde H_\eps(x)|^2}{4\varepsilon^2} \leq\frac{C_\Delta}{2
\varepsilon} -
\frac{c_\nabla a^2 }{4 \varepsilon}\leq- \frac{C_\Delta}{2
\varepsilon} \qquad\mbox{for all }
\frac{a}{2} \sqrt{\varepsilon} < \llvert x\rrvert _\delta< a
\sqrt{\varepsilon},
\end{eqnarray*}
which is the desired statement \eqref{elyapPIconstr2}. Therefore,
it is only left to deduce the estimate \eqref{eboundCnabla}. By the
definition of $\tilde H_\eps$ from above, we can write
\begin{equation}
\label{elyapPInablatildeHnorm} \bigl\llvert \nabla\tilde H_\eps(x)\bigr\rrvert
^2 = \bigl\llvert \nabla H (x)\bigr\rrvert ^2 + \bigl
\llvert \nabla H_b (x)\bigr\rrvert ^2 + 2 \bigl\langle
\nabla H (x) , \nabla H_b(x) \bigr\rangle.
\end{equation}
Let us have a closer look at each term on the right-hand side of the
last identity and let us start with the first term. By applying
Taylor's formula to $\nabla H(x)$, we obtain
\begin{eqnarray}
\bigl\llvert \nabla H (x) - \nabla^2 H (0) x\bigr\rrvert \leq\tilde
C_\nabla\llvert x\rrvert \stackrel {\mathclap{\eqref{elyapequivnorms}}}
{\leq} C_\nabla\llvert x\rrvert _\delta \label{elyapPILaplacenablaH}
\end{eqnarray}
for some $\tilde C_\nabla, C_\nabla>0$. Therefore, we can estimate
\begin{eqnarray}
\label{elyapPIestimatenablaH} \bigl\llvert \nabla H (x)\bigr\rrvert ^2 \geq\bigl
\llvert \nabla^2 H (0) x\bigr\rrvert ^2 -
C_\nabla^2 a^4 \varepsilon^2 \qquad
\mbox{for } \llvert x\rrvert _\delta\leq a \sqrt\varepsilon.
\end{eqnarray}
By the definition of $\lambda_1,\ldots,\lambda_n$, we also know
\begin{equation}
\label{elyapPInabla^2Hidentity} \bigl\llvert \nabla^2 H (0)x\bigr\rrvert
^2= \sum_{i=1}^n
\lambda_i^2 \bigl\llvert \langle u_i , x
\rangle\bigr\rrvert ^2 .
\end{equation}
Let us have a closer look at the second term in \eqref
{elyapPInablatildeHnorm}, namely $\llvert \nabla H_b (x)\rrvert ^2$. From the
definition \eqref{elyapPIdefHb} of $\llvert \nabla H_b (x)\rrvert ^2$ follows
\begin{eqnarray}
\label{elyapPInablaHbidentity}%
\bigl\llvert \nabla H_b (x)\bigr\rrvert
^2 &=& \bigl\llvert \xi' \bigl(\llvert x\rrvert
_\delta^2\bigr)\bigr\rrvert ^2 \Biggl( \sum
_{i=1}^\ell\delta^2 \bigl\llvert
\langle u_i , x \rangle \bigr\rrvert ^2 + \sum
_{i=\ell+1}^n (\lambda_i-
\delta)^2 \bigl\llvert \langle u_i , x \rangle \bigr
\rrvert ^2 \Biggr)
\nonumber
\\[-8pt]
\\[-8pt]
&\overset{\mathclap{\eqref{elyapdefnorm}}}
{\leq}& 2\lambda _{\max}^+ \llvert x\rrvert _\delta^2.
\nonumber
\end{eqnarray}
Now, we turn   the analysis of the last term, namely $2  \langle
\nabla H (x) , \nabla H_b(x)  \rangle$. By using the
estimates \eqref
{elyapPILaplacenablaH} and \eqref{elyapPInablaHbidentity},
we get for $\llvert x\rrvert _\delta\leq a \sqrt{\varepsilon}$.
\begin{eqnarray}\label{elyapPIskpestimate}
\bigl\langle\nabla H (x) , \nabla H_b(x) \bigr\rangle &=& \bigl
\langle\nabla^2 H(0) x , \nabla H_b(x) \bigr\rangle +
\bigl\langle\nabla H(x) - \nabla^2 H(0) x , \nabla H_b(x)
\bigr\rangle
\nonumber\\
&\mathop{\stackrel{\mathclap{\eqref {elyapPILaplacenablaH}}}
{\geq}}\limits_{\mathclap{\eqref
{elyapPInablaHbidentity}}}& \bigl\langle\nabla^2 H(0) x , \nabla
H_b(x) \bigr\rangle - 2C_\nabla\lambda_{\max}
\llvert x\rrvert _\delta^3\nonumber
\\[-8pt]
\\[-8pt]
&\geq& - \sum_{i=1}^\ell
\lambda_i \delta\bigl\llvert \xi' \bigl(\llvert x
\rrvert _\delta^2 \bigr)\bigr\rrvert \bigl\llvert \langle
u_i , x \rangle\bigr\rrvert ^2 \nonumber\\
&&{}- \sum
_{i=\ell+1}^n \lambda_i (
\lambda_i- \delta) \bigl\llvert \xi' \bigl(\llvert x
\rrvert _\delta^2 \bigr)\bigr\rrvert \bigl\llvert \langle
u_i , x \rangle \bigr\rrvert ^2 -O\bigl(
\varepsilon^{\fraca{3}{2}}\bigr).
\nonumber
\end{eqnarray}
Combining now the estimates and identities \eqref
{elyapPInablatildeHnorm}, \eqref{elyapPIestimatenablaH},
\eqref{elyapPInabla^2Hidentity}, \eqref
{elyapPInablaHbidentity} and \eqref{elyapPIskpestimate}, we
arrive for $\llvert x\rrvert _\delta\leq a \sqrt{\varepsilon}$ at
\begin{eqnarray*}
\bigl\llvert \nabla\tilde H_\eps(x)\bigr\rrvert ^2 &
\geq& \sum_{i=1}^\ell \bigl(\lambda
_i - \delta\bigl\llvert \xi' \bigl(\llvert x\rrvert
_\delta^2 \bigr)\bigr\rrvert \bigr)^2 \bigl
\llvert \langle u_i , x \rangle \bigr\rrvert ^2
\\
&& {} + \sum_{i=\ell+1}^n \bigl(
\lambda_i - (\lambda_i- \delta) \bigl\llvert
\xi' \bigl(\llvert x\rrvert _\delta^2 \bigr)
\bigr\rrvert \bigr)^2 \bigl\llvert \langle u_i , x
\rangle\bigr\rrvert ^2 - O\bigl(\varepsilon^{\fraca{3}{2}}\bigr).
\end{eqnarray*}
By \eqref{edefauxxi} holds $|\xi'(|x|_\delta^2)| \leq1$, which
applied to the last inequality yields
\[
\bigl\llvert \nabla\tilde H_\eps(x)\bigr\rrvert ^2 \geq
\delta^2 \sum_{i=1}^n \bigl
\llvert \langle u_i , x \rangle\bigr\rrvert ^2 - O\bigl(
\varepsilon ^{\fraca{3}{2}}\bigr).
\]
Because $u_1, \ldots, u_n$ is an orthonormal basis of $\R^n$, the
desired statement \eqref{eboundCnabla} follows for $ \frac{a \sqrt {\varepsilon}}{2} \leq\llvert x\rrvert _\delta\leq a \sqrt{\varepsilon}$ from
\begin{eqnarray*}
\bigl\llvert \nabla\tilde H_\eps(x)\bigr\rrvert ^2 &\geq&
 \delta^2 \llvert x\rrvert ^2 - O\bigl(\varepsilon
^{\fraca{3}{2}}\bigr) \overset{\mathclap{\eqref{elyapequivnorms}}}
{\geq} \frac
{2\delta^2}{\lambda_{\max}^+} \llvert x\rrvert _\delta^2 - O
\bigl(\varepsilon^{\fraca
{3}{2}}\bigr)
\\
& \geq& \frac{\delta^2}{2\lambda_{\max}^+} a^2 \varepsilon -O\bigl(
\varepsilon^{\fraca{3}{2}}\bigr) \geq c_{\nabla} a^2 \eps
\end{eqnarray*}
for some $c_{\nabla}< \frac{\delta^2}{2 \lambda_{\max}^+}$ and
$\varepsilon$ small enough.
\end{pf}
%

Considering the statement of Lemma~\ref{pLyapunovinside}, there is
only one thing to show in order to verify Lemma~\ref{pLyapunovverification}.
%
\begin{lem}\label{padmisibleverification}
Let $\cP_\cM= \{\Omega_i \}_{i=1}^M$ be the partition
obtained from
the $\tilde H_\eps$ from Lemma~\ref{pLyapunovinside} by considering
the basins of attraction $\Omega_i$ from \eqref
{erightchoiceofanpartition}. Then $\cP_\cM$ is an admissible
partition in the sense of Definition~\ref{defpartition}.
\end{lem}
Before we turn to the proof of Lemma~\ref{padmisibleverification},
we show the following auxiliary statement.

%
\begin{lem}\label{plineargrowth}
If an Hamiltonian $H\dvtx  \R^n \to\R$ satisfies the Assumption~\ref
{assumeenv}, then there exist numbers $R>0$ and $c_H>0$ such that
\[
H(x) \geq\min_{\llvert z\rrvert = R} H(z) + c_H \bigl(\llvert x
\rrvert -R\bigr).
\]
Because $H\geq0$ by Assumption~\ref{assumeenv}, a direct consequence
is $\int\exp(-H(x))  \,\dx{x} < \infty$.
\end{lem}
\begin{pf}
By the assumption \eqref{assumegradsuperlinear}, we can choose $R>0$
large enough such that
\begin{equation}
\label{elowerboundonnablaHlineargrowth} \bigl\llvert \nabla H (x)\bigr\rrvert \geq\frac{C_H}{2}
\qquad\mbox {for all } \llvert x\rrvert \geq R.
\end{equation}
In particular, this implies that for all critical points $s\in\cS$
holds $\llvert  s\rrvert \leq R$.
Now, let us we consider the following evolution:
\[
\dot{x_t} = - \frac{\nabla H(x_t)}{\llvert \nabla H (x_t)\rrvert }, \qquad x_0=x, 0\leq
t < \infty
\]
with starting point $x$, $\llvert x\rrvert > R$. Because by Lemma~\ref{lempartLebesgueGF}
\[
\R^n = \biguplus_{s \in\mathcal{S}} \Bigl\{y \in\R^n \dvtx  \lim
_{t\to \infty} y_t = s, \dot{y}_t=-\nabla
H(y_t), y_0 =y \Bigr\}
\]
and for all critical points $s \in\mathcal{S}$ of $H$ it holds $\llvert s\rrvert
\leq R$, the gradient line $ \{x_t  \}$ has to hit the ball
$B_R(0)$ after some time $t>0$ at some point $x_t$ for the first time.
It follows
\begin{eqnarray*}
H(x_t)-H(x_0) & =& \int_0^t
{\frac{\textup{d} }{\textup{d} s}} H (x_s) \,\dx{s}
\\
&= &- \int_0^t \nabla H(x_s)
\cdot\frac{\nabla H(x_s)}{|\nabla
H(x_s)|} \,\dx{s} = - \int_0^t
\bigl\llvert \nabla H(x_s) \bigr\rrvert \,\dx{s}.
\end{eqnarray*}
Using the lower bound \eqref{elowerboundonnablaHlineargrowth}
on $\llvert \nabla H(x_t)\rrvert $, we get that
\begin{eqnarray*}
H(x) = H(x_t) + \int_0^t \bigl
\llvert \nabla H(x_s) \bigr\rrvert \,\dx{s} \geq\inf
_{\llvert z\rrvert =R} H(z) + t \frac{c_H}{2}.
\end{eqnarray*}
Because the evolution $x_t$ moves at speed $1$, we know that $t$ is the
length of the gradient-flow line connecting the points $x$ and $x_t$.
However, this length cannot be shorter than $t \geq\llvert x\rrvert  - R$, which
yields the desired statement.
\end{pf}

\begin{pf*}{Proof of Lemma~\ref{padmisibleverification}}
We start with showing that $\tilde H_\eps$ has the same local minima
$\mathcal{M} =  \{m_1, \ldots, m_M  \}$ as the original
Hamiltonian $H$. Because
\[
\label{eperturbednatureofmodification} \tilde H_\eps(x)= H (x) \qquad\mbox{for all } x
\notin\bigcup_{y \in\mathcal{S}\setminus\mathcal{M} } B_{a \sqrt{\varepsilon
}} (y),\vadjust{\goodbreak}
\]
it suffices to show that $\tilde H_\eps$ has no local minima in the set
\[
\bigcup_{y \in\mathcal{S}\setminus\mathcal{M} } B_{a \sqrt
{\varepsilon}} (y).
\]
However, this statement follows directly from the estimate \eqref
{elyapunovinsidecritical}, that is,
\[
\frac{\Delta\tilde H_\varepsilon(x)}{2\varepsilon} - \frac{|\nabla
\tilde H_\varepsilon(x) |^2}{4\varepsilon^2} \leq- \frac{\lambda
_0}{\varepsilon} \qquad
\mbox{for all } x \in\bigcup_{y \in
\mathcal{S}\setminus\mathcal{M}} B_{a \sqrt{\varepsilon}}
(y) .
\]
Indeed, the last estimate shows that either $|\nabla\tilde
H_\varepsilon(x) |\neq0$ or $\Delta\tilde H_\varepsilon(x)<0$.

The fact that $\tilde H_\eps$ has the same local minima as $H$ allows
us to apply Lemma~\ref{lempartLebesgueGF} showing
\[
\R^n = \biguplus_{i=1}^M \Omega_i
= \biguplus_{i=1}^M \Bigl\{y \in \R^n \dvtx  \lim
_{t\to\infty} y_t = m_i,
\dot{y}_t=-\nabla\tilde H_\eps(y_t),
y_0 =m_i \Bigr\} ,
\]
which is already (ii) of Definition~\ref{defpartition}.

The last step in the proof is to show that $\mu(\Omega_i)Z_\mu$
satisfies the asymptotic expansion given by \eqref{ePartSumAdmPart}.
Let us consider one local minimum $m_i\in\cM$. W.l.o.g. we assume
$\tilde H_\eps(m_i)= H(m_i)=0$. We introduce $\Sigma_i :=  (
\nabla^2 H (m_i)  )^{-1}$ and define for $r_0>0$ specified later
the ellipsoid
\[
E_i := \bigl\{ x\in\R^n \dvtx  \bigl\llvert
\Sigma_i^{-\fraca{1}{2}} (x - m_i)\bigr\rrvert \leq
\sqrt{2 r_0 \eps\llvert \log\eps\rrvert } \bigr\},
\]
where the square root of $\Sigma_i^{-1}$ is uniquely defined in the
set of positive symmetric matrices.
Note that for small enough $\eps$ it holds $E_i \subset\Omega_i$ and
$\tilde H_\eps(x) = H(x)$ for $x \in E_i$.
The covariance matrix $\Sigma_i$ is nondegenerate because of $H$
being a Morse function. Therefore, there is a constant $c_H<1$ such that
\begin{equation}
\label{esmallcirclecontainedinEi} B_{\sqrt{c_H 2 r_0 \eps\llvert \log\eps\rrvert }}(m_i) \subset E_i
\subset B_{\sqrt{c_H^{-1} 2 r_0 \eps\llvert \log\eps\rrvert }}(m_i).
\end{equation}
We split the integral into
\[
\mu(\Omega_i) Z_\mu= \int_{E_i}
\exp \biggl(-\frac{\tilde H_\eps
(x)}{\eps} \biggr) \,\dx{x} + \int_{\Omega_i \setminus E_i}
\exp \biggl(-\frac{\tilde H_\eps(x)}{\eps} \biggr) \,\dx{x} =: I_1 +
I_2 .
\]
The results follows from an asymptotic expansion for $I_1$ and an error
estimate for $I_2$.

We start with the error estimate for $I_2$. Let the constant $R>0$ be
chosen as in Lemma~\ref{plineargrowth}. We split the term $I_2$ up into
\begin{eqnarray*}
I_2 & =& \int_{(\Omega_i \setminus E_i) \cap B_R(0)} \exp \biggl(-
\frac{\tilde H_\eps(x)}{\eps} \biggr) \,\dx{x} + \int_{\Omega_i
\setminus B_R(0)} \exp
\biggl(-\frac{\tilde H_\eps(x)}{\eps} \biggr) \,\dx{x}
\\
&=:& I_3 + I_4.
\end{eqnarray*}
Let us estimate the term $I_3$. On a small neighborhood around $m_i$ it
holds $H= \tilde H_\eps$ and $H$ is uniformly convex. Therefore, there
is a constant $\delta>0$ and $\kappa>0$ such that for all $x$ with
$\llvert x- m_i\rrvert  \leq\delta$
\[
\label{elyapadmPartconvex} \bigl\llvert \bigl(\nabla^2 H(x)\bigr)^{\fraca{1}{2}}
\xi\bigr\rrvert ^2 = \bigl\langle\xi , \nabla ^2 H(x)\xi
\bigr\rangle \geq\kappa\llvert \xi\rrvert ^2   \qquad\mbox{for all }
\xi\in \R^n .
\]
Hence, for $x \in\Omega_i \setminus E_i$ we have the lower bound by
additionally considering \eqref{esmallcirclecontainedinEi}
\begin{eqnarray*}
\tilde H_\eps(x) \geq\inf_{z \in\partial E_i} \tilde
H_\eps(z) \geq\frac{\kappa}{2} \inf_{z \in\partial E_i}
\llvert z-m_i\rrvert ^2 \geq \kappa c_H
r_0 \eps\llvert \log\eps\rrvert .
\end{eqnarray*}
Now, we can estimate $I_3$ as
\[
I_3 \leq\exp \bigl(- \kappa c_H r_0
\llvert \log\eps\rrvert \bigr) \bigl\llvert B_R(0)\bigr\rrvert .
\]
Let us turn to the estimation of $I_4$. An application of Lemma~\ref
{plineargrowth} yields
\begin{eqnarray*}
I_4 &\leq& \exp \Bigl(- \eps^{-1}\min_{\llvert z\rrvert =R}
H(z) \Bigr)\int_{\Omega_i \setminus B_R(0)} \exp \biggl(- c_H
\frac{\llvert x\rrvert  - R }{\eps
} \biggr)\,\mathrm{d}x
\\
&\leq& C_H \exp\bigl(- \kappa c_h r_0
\llvert \log\eps\rrvert \bigr).
\end{eqnarray*}
So overall, we have estimated the term $I_2$ as
\begin{eqnarray}
\label{eestimateadmissibleouter} I_2 \leq C_H \exp\bigl(- \kappa
c_h r_0 \llvert \log\eps\rrvert \bigr) =
C_H \eps^{\kappa
c_H r_0 } = O\bigl(\eps^\alpha\bigr) \nonumber
\\[-8pt]
\\[-8pt]
\eqntext{\displaystyle \mbox{for } r_0 > \frac{\alpha
}{\kappa c_H} \mbox{ and } \alpha>0. }
\end{eqnarray}
Hence, $I_2$ becomes smaller than every power of $\eps$ for $r_0$
large enough.

Now, we turn to the asymptotic approximation of the term $I_1$. The
Taylor expansion of $H$ on $E_i$ yields for $x \in E_i$
\[
H(x) = \tfrac{1}{2} \bigl\langle x , \nabla^2
H(m_i)x \bigr\rangle + O \bigl( \bigl(\eps\llvert \log \eps\rrvert
\bigr)^{\fraca{3}{2}} \bigr).
\]
In particular, this implies
\[
\exp \biggl(-\frac{H(x)}{\eps} \biggr) = \exp \biggl(-\frac
{1}{2\eps} \bigl
\langle x , \nabla^2 H(m_i)x \bigr\rangle \biggr) \exp
\bigl(O\bigl(\sqrt {\eps} \llvert \log\eps\rrvert ^{\fraca{3}{2}}\bigr) \bigr) .
\]
For $\eps$ small enough, it holds $\exp (O(\sqrt{\eps}  \llvert \log\eps\rrvert ^{\fraca{3}{2}}) ) =1+O(\sqrt{\eps
}  \llvert \log\eps \rrvert ^{\fraca{3}{2}})$. Therewith,
we get the following expression for $I_1$:
\begin{eqnarray}
I_1 &=& \int_{E_i} \exp \biggl(-
\frac{1}{2\eps} \bigl\langle x , \nabla ^2 H(m_i)x
\bigr\rangle \biggr) \,\dx{x} \bigl(1+O\bigl(\sqrt{\eps }\llvert \log\eps\rrvert
^{\fraca
{3}{2}}\bigr)\bigr)
\nonumber
\\
&=& \frac{(2\pi\eps)^{\fraca{n}{2}}}{\sqrt{\det\nabla^2 H(m_i)}}\nonumber\\
&&{}\times \biggl(1 - \frac{\sqrt{\det\nabla^2 H(m_i)}}{(2\pi\eps)^{\fraca
{n}{2}}} \int_{\R^n \setminus {E_i}}
\exp \biggl(-\frac {1}{2\eps
} \bigl\langle x , \nabla^2
H(m_i)x \bigr\rangle \biggr) \,\dx{x} \biggr)
\nonumber
\\
& &{} \times\bigl(1+O\bigl(\sqrt{\eps}\llvert \log\eps\rrvert ^{\fraca{3}{2}}
\bigr)\bigr).\nonumber %
\end{eqnarray}
Now, we apply the following tail estimate for a Gaussian, which we will
proofed for the convenience of the reader below:
\begin{equation}
\label{etailestimateGaussianintegral}\quad  \frac{\sqrt{\det\nabla^2 H(m_i)}}{(2\pi\eps)^{\fraca{n}{2}}} \int_{\R^n \setminus {E_i}} \exp \biggl(-
\frac{1}{2\eps} \bigl\langle x , \nabla^2 H(m_i)x \bigr
\rangle \biggr) \,\dx{x} = O(\sqrt {\eps}).
\end{equation}
The latter yields the asymptotic expansion
\begin{equation}
\label{eestimateadmissibleinner} I_1 = \frac{(2\pi\eps)^{\fraca{n}{2}}}{\sqrt{\det\nabla^2
H(m_i)}} \bigl(1+O\bigl(\sqrt{\eps}
\llvert \log\eps\rrvert ^{\fraca{3}{2}}\bigr)\bigr).
\end{equation}
Now, the desired asymptotic expansion \eqref{ePartSumAdmPart} for
$\mu(\Omega_i)Z_\mu$ follows form a combination of the expansion for
the term $I_1$ in \eqref{eestimateadmissibleinner} and $I_2$
in \eqref{eestimateadmissibleouter} with $\alpha$ chosen
sufficiently large, that is, $\alpha>(n+1)/2$.

We close the argument by deducing the desired tail estimate \eqref
{etailestimateGaussianintegral}. By the change of variables
$x\mapsto y =  (2\eps\Sigma_i )^{-\fraca{1}{2}}(x-m_i)$ and
by denoting $\omega(\eps) = \sqrt{r_0 \llvert \log\eps\rrvert }$, we deduce
\begin{eqnarray*}
&&\frac{\sqrt{\det\nabla^2 H(m_i)}}{(2\pi\eps)^{\fraca{n}{2}} } \int_{\R^n\setminus E_i} \exp \biggl(-
\frac{1}{2\eps} \bigl\langle x , \nabla^2 H(m_i)x \bigr
\rangle \biggr)\,\dx{x}\\ &&\qquad = \frac
{1}{\pi^{\fraca{n}{2}}} \int_{\R^n\setminus B_{\omega(\eps)}(0)}
e^{-y^2} \,\dx{y}
\\
&&\qquad = \frac{n}{\Gamma(\fraca{n}{2}+1)} \int_{\omega(\eps)}^\infty
r^{n-1} e^{-r^2} \,\dx{r} 
=
\frac{\Gamma(\fraca{n}{2},\omega^2(\eps))}{\Gamma(\fraca{n}{2})} ,
\end{eqnarray*}
where $\Gamma(\frac{n}{2},\omega^2(\eps))$ is the complementary
incomplete Gamma function. It has the asymptotic expansion (cf. \cite{Olver1997}, pp. 109--112)
\[
\Gamma \biggl(\frac{n}{2},\omega^2(\eps) \biggr) = O
\bigl(e^{-\omega
^2(\eps
)}\omega^{n-2}(\eps)\bigr)=O\bigl(\eps^{r_0}
\llvert r_0\log\eps\rrvert ^{\fraca
{n}{2}-1}\bigr) = O(\sqrt{\eps})
\]
for $r_0$ large enough, which yields the desired result.
\end{pf*}

\subsection{Lyapunov approach for the logarithmic Sobolev inequality}
\label{slyapunovLSI}

The goal of this section is to prove Theorem~\ref{thmlocalLSI}
deducing the local LSI. We follow the same strategy as for the proof of
Theorem~\ref{thmlocalPI}, which we outlined in Section~\ref{slyap}. Therefore, we consider the partition $\cP_\cM= \{
\Omega _i \}_{i=1}^M$ into the basins of attraction of the
gradient flow of
$\tilde H_\eps$ [cf. \eqref{erightchoiceofanpartition}].

The Lyapunov condition for proving LSI is stronger than the one for PI.
Nevertheless, the construction of the $\eps$-modified Hamiltonian
$\tilde H_\eps$ from the previous section carries over and we can use
the same Lyapunov function as for the PI, but have to provide
additional estimates. The Lyapunov condition for LSI goes back to the
work of Cattiaux et al. \cite{Wang2007}. We adapt \cite{Wu2008}, Theorem~1.2, to the case for domains $\Omega$. In addition, we will
work out the explicit dependence between the constants of the Lyapunov
condition, the logarithmic Sobolev constant and especially their $\eps
$-dependence.
%
\begin{theo}[(Lyapunov condition for $\LSI$)]\label{localthmlyapLSI}
Suppose that:
\begin{longlist}[(iii)]
\item There exists a $C^2$-function $W\dvtx \Omega\to[1,\infty)$ and
constants $\lambda, b >0$ such that for $L = \eps\laplace- \nabla H
\cdot\nabla$ holds
\begin{equation}
\label{locallyapecondLSI} \forall x\in\Omega\dvtx  \frac{1}{\eps} \frac{L W}{W} \leq-
\lambda \llvert x\rrvert ^2 +b .
\end{equation}
\item$\nabla^2 H \geq- K_H$ for some $K_H > 0$ and $\mu$ satisfies
\ref{PIvarrho}.
\item$W$ satisfies Neumann boundary conditions on $\Omega$
[cf. \eqref{elyapintegrationbyparts}].
\end{longlist}
Then $\mu$ satisfies \ref{LSIalpha} with
\begin{equation}
\label{localthmlyapconstLSI} \frac{1}{\alpha} \leq2\sqrt{ \frac{1}{\lambda} \biggl(
\frac
{1}{2}+\frac{b+\lambda\mu(\llvert  x\rrvert ^2)}{\varrho
} \biggr) }+ \frac
{K_H}{2\eps\lambda}+
\frac{K_H(b+\lambda\mu(\llvert  x\rrvert ^2))+2\eps
\lambda}{\varrho\eps\lambda} ,\hspace*{-20pt}
\end{equation}
where $\mu(\llvert  x\rrvert ^2)$ denotes the second moment of
$\mu$. 
\end{theo}
Before turning to the proof of Theorem \eqref
{localthmlyapconstLSI}, we need the following auxiliary result.
%
\begin{lem}[(\cite{Wu2008}, Lemma~3.4)]\label{locallemlyapLSI}
Assume that $V\dvtx \Omega\to\R$ is a nonnegative locally Lipschitz
function such that:
\begin{longlist}[(ii)]
\item For some lower bounded function $\phi$
\begin{equation}
\label{elyapLSIlemma} \frac{L e^V}{e^V} = L V + \eps\llvert \nabla V\rrvert
^2 \leq-\eps\phi
\end{equation}
in the distributional sense.
\item$V$ satisfies Neumann boundary condition on $\Omega$ [cf. \eqref
{elyapintegrationbyparts}].
\end{longlist}
Then for any $g\in H^1(\mu)$ holds
\[
\int\phi g^2 \,\dx{\mu} \leq\int\llvert \nabla g\rrvert
^2 \,\dx{\mu} .
\]
\end{lem}
\begin{pf}
We can assume w.l.o.g. that $g$ is smooth with bounded support and
$\phi$ is bounded. For the verification of the desired statement, we
need the symmetry of $L$ in $L^2(\mu)$ w.r.t. to $V$:
\begin{equation}
\label{esymofL} \forall f\in H^1(\mu) \dvtx  \int f (-L V) \,\dx{\mu} =
\eps \int\nabla f \cdot\nabla V \,\dx{\mu},
\end{equation}
and the Young inequality:
\begin{equation}
\label{etrivest} 2g \nabla V \cdot\nabla g \leq\llvert \nabla V\rrvert
^2 g^2 + \llvert \nabla g\rrvert ^2 .
\end{equation}
An application of the assumption \eqref{elyapLSIlemma} yields
\begin{eqnarray*}
\eps\int\phi g^2 \,\dx{\mu} &\overset{\mathclap{\eqref
{elyapLSIlemma}}} {\leq}& \int \bigl(-L
V - \eps\llvert \nabla V\rrvert ^2 \bigr)g^2 \,\dx{\mu}
\\
& \overset{\mathclap{\eqref{esymofL}}}
{=}& \eps\int \bigl(2 g \nabla V\cdot\nabla g - \llvert \nabla V\rrvert
^2 g^2 \bigr) \,\dx{\mu} \overset {\mathclap{
\eqref{etrivest}}} {\leq} \eps\int\llvert \nabla g
\rrvert ^2 \,\dx{\mu},
\end{eqnarray*}
which is the desired estimate.
\end{pf}
The proof of Theorem~\ref{localthmlyapLSI} relies on an interplay
of some other functional inequalities, which will not occur anywhere else.
%
\begin{pf*}{Proof of Theorem~\ref{localthmlyapLSI}}
The argument of \cite{Wu2008} is a combination of the Lyapunov
condition \eqref{locallyapecondLSI} leading to a defective $\WI$
inequality and the use of the $\HWI$ inequality of Otto and
Villani \cite{Otto2000}. In the following, we will use the measure
$\nu$ given by $\nu(\dx{x}) = h(x)\mu(\dx{x})$, where we can
assume w.l.o.g. that $\nu$ is a probability measure, that is, $\int h
\,\dx{\mu}=1$. The first step is to estimate the Wasserstein
distance in terms of the total variation \cite{Villani2009}, Theorem~6.15
\begin{equation}
\label{locallemlyapLSIp1} W_2^2(\nu,\mu) \leq2 \bigl\llVert \llvert
\cdot\rrvert ^2 (\nu-\mu)\bigr\rrVert _{\mathrm{TV}} .
\end{equation}
For every function $g$ with $\llvert  g\rrvert \leq\phi(x):=
\lambda\llvert  x\rrvert ^2$, where $\lambda$ is from the
Lyapunov condition \eqref
{locallyapecondLSI} we get
\begin{eqnarray}
\label{locallemlyapLSIp2}%
\int g \,\dx{(\nu- \mu)} &\leq& \int\phi \,\dx{\nu} +
\int\phi \,\dx{\mu}
\nonumber
\\[-8pt]
\\[-8pt]
&=& \int\bigl(\lambda\llvert x\rrvert ^2-b\bigr) h(x) \mu(\dx{x}) +
\int b \,\dx {\nu} + \mu(\phi) .
\nonumber
\end{eqnarray}
We can apply to $\int(\lambda\llvert  x\rrvert ^2-b) h   \,\dx
{\mu}$
Lemma~\ref{locallemlyapLSI}, where the assumptions are exactly the
Lyapunov condition \eqref{locallyapecondLSI} by choosing $V=\log
W$. Moreover, the Neumann condition also translates to $V$ since $W$ is
bounded from below by $1$. Therewith, we arrive at
\begin{equation}
\label{locallemlyapLSIp3} \int\bigl(\lambda\llvert x\rrvert ^2-b\bigr) h \,\dx{
\mu} \leq\int \llvert \nabla\sqrt {h}\rrvert ^2 \,\dx{\mu} = \int
\frac
{\llvert \nabla h\rrvert ^2}{4h} \,\dx{\mu} = \frac{1}{2} I(\nu|\mu) ,
\end{equation}
by the definition of the Fisher information. Taking the supremum over
$g$ in \eqref{locallemlyapLSIp2} and combining the estimate
with \eqref{locallemlyapLSIp1} and \eqref{locallemlyapLSIp3}
we arrive at the defective Wasserstein-information inequality
\begin{equation}
\label{locallemlyapLSIdefW2I} \frac{\lambda}{2} W_2^2(\nu,\mu) \leq
\lambda\bigl\llVert \llvert \cdot \rrvert ^2 (\nu-\mu)\bigr\rrVert
_{\mathrm{TV}} \leq\frac
{1}{2}I(\nu|\mu) + b + \mu(\phi) .
\end{equation}
The next step is to use the $\HWI$ inequality \cite{Otto2000}, Theorem~3, which holds by the assumption $\nabla^2 H \geq-K_H$
\[
\Ent_{\mu}(h) \leq W_2(\nu,\mu) \sqrt{2 I(\nu| \mu)} +
\frac
{K_H}{2\eps} W_2^2(\nu,\mu) .
\]
Substituting inequality \eqref{locallemlyapLSIdefW2I} into the
$\HWI$ inequality and using the Young inequality $ab\leq\frac{\tau
}{2} a^2 + \frac{1}{2\tau}b^2$ for $\tau>0$ results in
\begin{eqnarray}
\label{locallemlyapLSIdefLSI} %
\Ent_{\mu}(h)& \leq& \tau I(\nu|\mu) +
\biggl(\frac {1}{2\tau
} + \frac{K_H}{2\eps} \biggr) W_2^2(
\nu,\mu)
\nonumber
\\[-8pt]
\\[-8pt]
&\overset{\mathclap{\eqref{locallemlyapLSIdefW2I}}} { \leq }& \biggl(\tau+
\frac{1}{2\lambda} \biggl(\frac{1}{\tau}+\frac
{K_H}{\eps} \biggr) \biggr)
I(\nu|\mu) + \frac{1}{\lambda} \biggl(\frac{1}{\tau}+\frac{K_H}{\eps}
\biggr) \bigl(b+\mu(\phi)\bigr) . %
\nonumber
\end{eqnarray}
The last inequality is of the type $\Ent_{\mu}(h)\leq\frac
{1}{\alpha_d} I(\nu|\mu) + B \int h \,\dx{\mu}$\vspace*{1pt} and is often called
defective logarithmic Sobolev inequality $\dLSI(\alpha_d,B)$. It is
well known that a defective logarithmic Sobolev inequality can be
tightened by \ref{PIvarrho} to \ref{LSIalpha} with constant (cf.
Proposition \cite{Ledoux1999a})
\begin{equation}
\label{locallemlyapLSIconstdefLSI} \frac{1}{\alpha} = \frac{1}{\alpha_d} + \frac{B+2}{\varrho} .
\end{equation}
A combination of \eqref{locallemlyapLSIdefLSI} and \eqref
{locallemlyapLSIconstdefLSI} reveals
\begin{eqnarray}
\frac{1}{\alpha} &=& \tau+ \frac{1}{2\lambda} \biggl(
\frac
{1}{\tau}+\frac{K_H}{\eps} \biggr) + \frac{1}{\varrho} \biggl(
\frac {1}{\lambda} \biggl(\frac{1}{\tau}+\frac{K_H}{\eps} \biggr)
\bigl(b+\mu (\phi)\bigr)+2 \biggr)
\nonumber
\\
&=& \tau+\frac{1}{\tau\lambda} \biggl(\frac{1}{2}+\frac{b+\mu
(\phi)}{\varrho}
\biggr) + \frac{K_H}{2\eps\lambda}+\frac
{K_H(b+\mu
(\phi))+2\eps\lambda}{\varrho\eps\lambda} =: \tau+ \frac
{c_1}{\tau}+c_2.\nonumber
\end{eqnarray}
The last step is to optimize in $\tau$, which leads to $\tau=\sqrt {c_1}$ and, therefore, $\frac{1}{\alpha}=2\sqrt{c_1}+c_2$. The final
result \eqref{localthmlyapconstLSI} follows by recalling the
definition of $\phi(x)=\lambda\llvert  x\rrvert ^2$.
\end{pf*}

The crucial ingredient is a Lyapunov function satisfying the
condition \eqref{locallyapecondLSI}. We follow the ideas of
Section~\ref{slyap} and Section~\ref{sconstructionlyapunov}.
We use the same $\eps$-modification $\tilde H_\eps$ as constructed in
the proof of Lemma~\ref{pLyapunovinside}.

%
\begin{lem}[(Lyapunov function for $\LSI$)]\label{pLyapunovverificationLSI}
We consider the $\eps$-modification $\tilde H_\varepsilon$ of $H$
constructed in Section~\ref{sconstructionlyapunov}. Then the
Lyapunov function $W(x)= \exp ( \frac{1}{2\eps} \tilde
H_\varepsilon(x)  ) $ satisfies on $\Omega$ the Lyapunov
condition \eqref{locallyapecondLSI} with constants
\[
\label{locallyapeconstLSI} b = \frac{b_0}{\eps}  \quad \mbox{and} \quad \lambda\geq
\frac{\lambda_0}{
\eps}
\]
for some $b_0, \lambda_0>0$ and Hessian $\nabla^2 \tilde H_\eps(x)
\geq- K_{\tilde H_\eps}$ for some $K_{\tilde H_\eps}\geq0$.
\end{lem}
%

The proof consists of three steps, which correspond to three regions of
interests. First, we will consider a neighborhood of $\infty$, that
is, we will fix some $\tilde R >0$ and only consider $\llvert
x\rrvert \geq
\tilde R$. Then we will look at an intermediate regime for $a\sqrt {\eps}\leq\llvert  x\rrvert \leq\tilde R$, where we will
have to take special
care for the neighborhoods around critical points and use the
construction of Lemma~\ref{pLyapunovinside}. The last regime is for
$\llvert  x\rrvert \leq a\sqrt{\eps}$, which will be the
simplest case.

Therefore, besides the construction done in the proof of Lemma~\ref
{pLyapunovinside}, we need an analogous formulation of Lemma~\ref
{plyapunovPIoutside} under the stronger assumption \eqref{assumelyapLSI}.
%
\begin{lem}\label{plyapunovLSIoutside}
Assume that the Hamiltonian $H$ satisfies assump-\break tion \eqref
{assumelyapLSI}. Then there is a constant $0\leq C_H < \infty$ and
$0 \leq\tilde R < \infty$ such that $H(x) = \tilde H_\eps(x)$ for
$\llvert x\rrvert \geq\tilde R$ and for $\eps$ small enough
\begin{equation}
\label{elyapunovoutsidetildeRLSI} \frac{\Delta H(x)}{2\varepsilon} - \frac{|\nabla H(x)
|^2}{4\varepsilon^2} \leq-
\frac{C_H}{\varepsilon} \llvert x\rrvert ^2 \qquad\mbox{for all } \llvert
x\rrvert \geq\tilde R.
\end{equation}
\end{lem}
We skip the proof of the Lemma~\ref{plyapunovLSIoutside}, because
it would work in the same way as for Lemma~\ref{plyapunovPIoutside}
and only consists of elementary calculations based on the nondegeneracy
assumption on $H$. The only difference, is that we now demand the
stronger statement \eqref{elyapunovoutsidetildeRLSI}, which is a
consequence of the stronger assumption \eqref{assumelyapLSI} in
comparison to assumption \eqref{assumegradlaplace}.

Now, we have collected the auxiliary statements and can proof
Lemma~\ref{pLyapunovverificationLSI}.

\begin{pf*}{Proof of Lemma~\ref{pLyapunovverificationLSI}}
First, let us check the lower bound on the Hessian of $\tilde H_\eps$.
Because we use the same $\tilde H_\eps$ as constructed in Lemma~\ref
{pLyapunovinside}, the support of $\tilde H_\eps- H$ is compact.
Additionally, $\tilde H_\eps$ is smooth. This already implies the
lower bound on the Hessian $\nabla^2 \tilde H_\eps$ for compact
domains. Outside a sufficient large domain, we know that $H=\tilde
H_\eps$. Hence, the lower bound on $\nabla^2 \tilde H_\eps$ follows
directly from assumption \eqref{assumehessLSI}.

Now, we verify the Lyapunov condition \eqref{locallyapecondLSI}.
Recall that $W=\exp (\frac{1}{2\eps} \tilde H_\eps )$.
Hence, straightforward calculation reveals
\[
\label{pLyapunovverificationLSIp0} \frac{1}{\varepsilon}\frac{L W}{W} = \frac{1}{2\eps}
\laplace \tilde H_\eps+ \frac{1}{4\eps^2}\llvert \nabla\tilde
H_\eps \rrvert ^2 - \frac{1}{2\eps^2} \llvert \nabla
\tilde H_\eps\rrvert ^2= \frac{1}{2\eps} \laplace\tilde
H_\eps- \frac{1}{4\eps^2}\llvert \nabla\tilde H_\eps
\rrvert ^2 .
\]
If $\llvert  x\rrvert \geq\tilde R$ with $\tilde R$ given in
Lemma~\ref
{plyapunovLSIoutside}, we apply \eqref
{elyapunovoutsidetildeRLSI} and have the Lyapunov condition
fulfilled with constant $\lambda=\frac{C_H}{\eps}$. This allows us
to only consider $x\in B_{\tilde R}\cap\Omega$, which is a bounded
domain. In this case, Lemma~\ref{pLyapunovinside} yields for $a\sqrt {\eps} \leq\llvert  x\rrvert \leq\tilde R$ the estimate
\begin{equation}
\label{pLyapunovverificationLSIp1} \frac{1}{\eps} \frac{L W}{W} \leq- \frac{\lambda_0}{\eps}
\leq -\frac{\lambda_0}{\tilde R^2 \eps} \llvert x\rrvert ^2 .
\end{equation}
Let us consider the final case $\llvert  x\rrvert  \leq a\sqrt {\eps}$. In this
case, the Hamiltonian $H=\tilde H_\eps$. Additionally, $H$ is smooth
and strictly convex on $B_{a\sqrt{\eps}}(0)$. Therefore, one easily
obtains the bound
\begin{equation}
\label{pLyapunovverificationLSIp2} \frac{1}{\eps} \frac{L W}{W} \leq\frac{1}{2\eps}
\laplace H(x) \leq\frac{b_0}{\eps} .
\end{equation}
A combination of \eqref{pLyapunovverificationLSIp1} and \eqref
{pLyapunovverificationLSIp2} yields the desired estimate \eqref
{locallyapecondLSI}.
\end{pf*}
Before proceeding with the proof of Theorem~\ref{thmlocalLSI}, we
remark, that the Lyapunov condition for the PI and in particular for
the LSI imply an estimate of the second moment of $\mu$.
%
\begin{lem}[(Second moment estimate)]\label{local2ndmomlem}
If $H$ fulfills the Lyapunov condition \eqref{locallyapcond} with
$U=B_R(0)$ for $R>0$, then $\mu$ has finite second moment and it holds
\begin{equation}
\label{local2ndmomgenest} \int\llvert x\rrvert ^2 \mu(\dx{x}) \leq
\frac{1+b
R^2}{\lambda}.
\end{equation}
\end{lem}
\begin{pf}
As it is outlined in the proof of Theorem~\ref{localthmlyapPI}
(cf. also \cite{Bakry2008}), the Lyapunov condition \eqref
{locallyapcond} yields the following estimate: for any function $f$
and $m\in\R$ it holds
\[
\int (f-m )^2 \,\dx{\mu} \leq\frac{1}{\lambda} \int \llvert \nabla
f\rrvert ^2 \,\dx{\mu} + \frac{b}{\lambda
}\int_{B_R(0)}
(f-m )^2 \,\dx{\mu} .
\]
We set $f(x)=\llvert  x\rrvert $ and $m=0$ to observe the
estimate \eqref
{local2ndmomgenest}.
\end{pf}
%
Now, we have collected all auxiliary results to proof the second main
Theorem~\ref{thmlocalLSI}.
\begin{pf*}{Proof of Theorem~\ref{thmlocalLSI}}
For the same reason as in the proof of Theorem~\ref{thmlocalPI}, we
omit the index $i$. The first step is also the same as in the proof of
Theorem~\ref{thmlocalPI}. By Lemma~\ref{locallemmodification}, we
obtain that, whenever $\tilde H_\eps$ is an $\eps$-modification of
$\mu$ in the sense of Definition~\ref{localdeftildeH}, the
logarithmic Sobolev constants $\alpha$ and $\tilde\alpha$ of $\mu$
and $\tilde\mu_\eps$ satisfy $\alpha\geq\exp (-2C_{\tilde
H} ) \tilde\alpha$.

The next step is to construct an explicit $\eps$-modification $\tilde
H$ satisfying the Lyapunov condition \eqref{locallyapecondLSI} of
Theorem~\ref{localthmlyapLSI}, which is provided by Lemma~\ref
{pLyapunovverificationLSI}.

Additionally, the logarithmic Sobolev constant $\tilde\alpha$ depends
on the second moment of $\tilde\mu_\eps$. Since $\tilde H_\eps$
satisfies by Lemma~\ref{pLyapunovverification} in particular the
Lyapunov condition for PI \eqref{locallyapcond} with constants
$\lambda\geq\frac{\lambda_0}{\eps}$, $b\leq\frac{b_0}{\eps}$
and $R=a\sqrt{\eps}$, we can apply Lemma~\ref{local2ndmomlem} and
arrive at
\[
\int\llvert x\rrvert ^2 \,\dx{\tilde\mu_\eps} \leq
\frac
{1+R^2 b}{\lambda
}\leq\frac{1+b_0 a^2}{\lambda_0} \eps= O(\eps).
\]
%
Now, we have control on all constants occurring in \eqref
{localthmlyapconstLSI} and can determine the logarithmic Sobolev
constant $\tilde\alpha$ of $\tilde\mu_\eps$. Let us estimate term
by term of \eqref{localthmlyapconstLSI} and use the fact from
Theorem \eqref{thmlocalPI}, that $\tilde\mu_\eps$ satisfies $\PI
(\tilde\varrho)$ with $\tilde\varrho^{-1} = O(\eps)$
\[
2\sqrt{ \frac{1}{\lambda} \biggl(\frac{1}{2} +\frac{b+\lambda
\tilde\mu_\eps(\llvert  x\rrvert ^2)}{\varrho}
\biggr)}\leq 2 \sqrt{\frac{\eps
}{\lambda_0} \biggl(\frac{1}{2} +O(1) \biggr)}
= O(\sqrt{\eps}) .
\]
The second term evaluates to $\frac{K_H}{2\eps\lambda}=O(1)$ and
finally the last one
\[
\frac{K_H(b+\lambda\tilde\mu_\eps(\llvert  x\rrvert
^2))+2\eps\lambda
}{\varrho\eps\lambda} = O(\eps) \biggl(K_H \biggl(\frac {b_0}{\eps}
+ O(\eps) \biggr)+ O(1) \biggr) = O(1) .
\]
A combination of all the results leads to the conclusion $\tilde\alpha
^{-1}= O(1)$ and since $\tilde H_\eps$ is only an $\eps$-modification
of $H$ also $\alpha^{-1}=O(1)$.
\end{pf*}

\section{Mean-difference estimates---weighted transport
distance}\label{CMEANDIFF}

This section is devoted to the proof of Theorem~\ref{thmcoarsemd}.
We want to estimate the mean-difference $ (\Expect_{\mu_i} f -
\Expect_{\mu_j} f )^2$ for $i$ and $j$ fixed. The proof
consists of
four steps:

In the first step, we introduce the \emph{weighted transport distance}
in Section~\ref{sectranspmd}. This distance depends on the transport
speed similarly to the Wasserstein distance, but in addition weights
the speed of a transported particle w.r.t. the reference measure $\mu
$. The weighted transport distance allows in general for a variational
characterization of the constant $C$ in the inequality
\[
\bigl(\Expect_{\mu_i}(f)-\Expect_{\mu_j}(f) \bigr)^2
\leq C \int \llvert \nabla f\rrvert ^2 \,\dx{\mu} .
\]
The problem of finding good estimates of the constant $C$ is then
reduced to the problem of finding a good transport interpolation
between the measures $\mu_i$ and $\mu_j$ w.r.t. to the weighted
transport distance.

For measures as general as $\mu_i$ and $\mu_j$, the construction of
an explicit transport interpolation is not feasible. Therefore, the
second step consists of an approximation, which is done in Section~\ref{sapproxmd}. There, the restricted measures $\mu_i$ and $\mu_j$ are
replaced by \emph{simpler} measures $\nu_i$ and $\nu_j$, namely
truncated Gaussians. We show in Lemma~\ref{mdapproxcorchi2} that
this approximation only leads to higher order error terms.

The most import step, the third one, consists of the estimation of the
mean-difference w.r.t. the approximations $\nu_i$ and $\nu_j$.
Because the structure of $\nu_i$ and $\nu_j$ is very simple, we can
explicitly construct a transport interpolation between $\nu_i$
and $\nu_j$ (see Lemma~\ref{ptranspcostapprox} in Section~\ref{safftransp}). The last step consists of collecting and controlling
the error (cf. Section~\ref{smdaffcomperror}).

\subsection{Mean-difference estimates by transport}\label{sectranspmd}

At the moment, let us consider two arbitrary measures $\nu_0\ll\mu$
and $\nu_1\ll\mu$.
The starting point of the estimation is a representation of the
mean-difference as a transport interpolation. This idea goes back
to \cite{Chafai2010}. However, they used a similar but nonoptimal
estimate for our purpose. Hence, let us consider a transport
interpolation $(\Phi_s \dvtx \R^n \to\R^n)_{s\in[0,1]}$ between $\nu
_0$ and $\nu_1$, that is, the family $(\Phi_s)_{s\in[0,1]}$ satisfies
\[
\Phi_0 = \operatorname{Id} , \qquad (\Phi_1)_\sharp
\nu_0 = \nu_1   \quad \mbox{and}\quad(
\Phi_s)_\sharp\nu_0 =: \nu_s .
\]
The representation of the mean-difference as a transport interpolation
is attained by using the fundamental theorem of calculus, that is,
\[
\bigl(\Expect_{\nu_0}(f) - \Expect_{\nu_1}(f)
\bigr)^2 
= \biggl(\int_0^1
\int \bigl\langle\nabla f(\Phi_s) , \dot\Phi_s \bigr
\rangle \,\dx{\nu_0} \,\dx{s} \biggr)^2 .
\]
At this point, it is tempting to apply the Cauchy--Schwarz inequality
in $L^2(\dx{\nu_0} \times\dx{s} )$ leading to the estimate in \cite
{Chafai2010}. However, this strategy would
not yield the preexponential factors in the Eyring--Kramers
formula \eqref{eeyringkramersPI} (cf. Remark~\ref{mdcompchafai}).
On Stephan Luckhaus' advice, the authors realized
the fact that it really matters on which integral you apply the
Cauchy--Schwarz inequality. This insight lead to the following proceeding:
\begin{eqnarray}\label{mddefnestimate}
\bigl(\Expect_{\nu_0}(f) - \Expect_{\nu_1}(f)
\bigr)^2 
&=& \biggl(\int
_0^1 \int \bigl\langle\nabla f, \dot
\Phi_s\circ \Phi_s^{-1} \bigr\rangle \,\dx{
\nu_{s}} \,\dx{s} \biggr)^2
\nonumber
\\
&=& \biggl(\int \biggl\langle\nabla f , \int_0^1
\dot\Phi_s \circ \Phi _s^{-1} {
\frac{\textup{d} \nu_{s}}{\textup{d} \mu}} \,\dx {s} \biggr\rangle \,\dx{\mu} \biggr)^2
\\
&\leq& \int\biggl\llvert \int_0^1 \dot
\Phi_s \circ\Phi_s^{-1} {\frac{\textup{d} \nu_{s}}{\textup{d} \mu}}
\,\dx{s}\biggr\rrvert ^2 \,\dx{\mu} \int\llvert \nabla f\rrvert
^2 \,\dx {\mu}.
\nonumber
\end{eqnarray}
Note that in the last step we have applied the Cauchy--Schwarz
inequality only in $L^2(\dx{\mu} )$ and that the desired Dirichlet
integral $\int\llvert \nabla f\rrvert ^2 \,\dx{\mu}$ is
already recovered.

The prefactor in front of the   Dirichlet energy on the right-hand
side of \eqref{mddefnestimate} only depends on the transport
interpolation $(\Phi_s)_{s\in[0,1]}$. Hence, we can minimize over all
possible admissible transport interpolations and arrive at the
following definition.

%
\begin{defn}[(Weighted transport distance $\cT_\mu$)]\label
{mddefnweightedTransp}
Let $\mu$ be an absolutely continuous probability measure on $\R^n$
with connected support. Additionally, let $\nu_0$ and $\nu_1$ be two
probability measures such that $\nu_0\ll\mu$ and $\nu_1\ll\mu$,
then define the \emph{weighted transport distance} by
\begin{equation}
\label{mddefnequweightedTransp} \cT_{\mu}^2(\nu_0,
\nu_1) := \inf_{{\Phi_s}} \int\biggl\llvert \int
_0^1 {\dot\Phi_s \circ
\Phi_s^{-1}} {\frac{\textup{d} \nu
_{s}}{\textup{d} \mu}} \,\dx{s}\biggr\rrvert
^2 \,\dx{\mu}.
\end{equation}

The family $(\Phi_s)_{s\in[0,1]}$ is chosen absolutely continuous in
the parameter $s$ such that $\Phi_0 = \Id$ on $\supp\nu_0$ and
$(\Phi_1)_\sharp\nu_0 = \nu_1$.
For a fixed family and $(\Phi_s)_{s\in[0,1]}$ and a point $x\in\supp
\mu$, the \emph{cost density} is defined by
\begin{equation}
\label{mdcostdensity} \cA(x) := \biggl\llvert \int_0^1
\dot\Phi_s \circ\Phi_s^{-1}(x) \nu
_{s}(x) \,\dx{s} \biggr\rrvert .
\end{equation}
\end{defn}
%
%
\begin{rem}[(Relation of $\cT_\mu$ to \cite{Chafai2010})]\label
{mdcompchafai}
The transport distance $\cT_\mu(\nu_0,\nu_1)$ is always smaller
than the constant obtained in \cite{Chafai2010}, Section~4.6. Indeed,
applying the Cauchy--Schwarz inequality on $L^2(\dx{s})$ in \eqref
{mddefnequweightedTransp} yields
\begin{eqnarray*}
\cT_\mu^2(\nu_0,\nu_1) &\leq&
\inf_{\Phi_s} \int\int_0^1
\bigl\llvert \dot\Phi_s \circ\Phi_s^{-1}\bigr
\rrvert ^2 {\frac{\textup
{d} \nu_{s}}{\textup{d} \mu}} \,\dx {s} \int_0^1
{\frac{\textup{d} \nu_{s}}{\textup{d} \mu}} \,\dx {s} \,\dx{\mu}
\\
&\leq& \inf_{\Phi_s} \biggl(\sup_x \biggl(
\int_0^1 {\frac{\textup
{d} \nu_s}{\textup{d} \mu}}(x) \,\dx{s}
\biggr) \int\int_0^1 \llvert \dot
\Phi_s \rrvert ^2 \,\dx{s} \,\dx{\nu _{0}}
\biggr),
\end{eqnarray*}
where we used the assumption that $\nu_{s}\ll\mu$ for all $s\in
[0,1]$ in the last $L^1$-$L^\infty$-estimate.
\end{rem}
%
%
\begin{rem}[(Relation of $\cT_\mu$ to the $L^2$-Wasserstein
distance $W_2$)]\label{mdcompwasser}
If the support of $\mu$ is convex, we can set the transport
interpolation $(\Phi_s)_{s\in[0,1]}$ to the linear interpolation map
$\Phi_s(x) = (1-s)x + s U(x)$. Assuming that $U$ is the optimal
$W_2$-transport map between $\nu_0$ and $\nu_1$, the estimate in
Remark~\ref{mdcompchafai} becomes
\[
\label{mdestwasserstein} \cT_\mu^2(\nu_0,
\nu_1) \leq \biggl(\sup_x \int
_0^1 {\frac
{\textup{d} \nu_s}{\textup{d} \mu}}(x) \,\dx{s} \biggr)
W_2^2(\nu _0,\nu_1) .
\]
\end{rem}
%
%
\begin{rem}[(Invariance under time rescaling)]\label{mdremtimerescale}
The cost density $\cA$ given by \eqref{mdcostdensity} is
independent of rescaling the transport interpolation in the parameter
$s$. Indeed, we observe that
\[
\cA(x) = \biggl\llvert \int_0^1 \dot
\Phi_s \circ\Phi_s^{-1}(x) \nu
_s(x) \,\dx{s} \biggr\rrvert = \biggl\llvert \int
_0^T \dot\Phi _t^T
\circ\bigl(\Phi _t^T\bigr)^{-1}(x)
\nu_{t}^T(x) \,\dx{t} \biggr\rrvert ,
\]
where $\Phi_t^{T} = \Phi_{t/T}$ and $\nu_t^T = \nu_{t/T}$.
\end{rem}
%
%
\begin{rem}[(Relation to negative Sobolev-norms)]\label{remTransportNegSobolev}
The weighted transport distance is a dynamic formulation for the
homogeneous negative Sobolev norm $\dot H^{-1}(\dx{\mu})$ like
Benamou and Brenier did for the Wasserstein distance \cite
{Brenier2000}. Precisely, for $\nu_0 = \varrho_0 \mu$ and $\nu_1 =
\varrho_1 \mu$ holds
\[
\cT_{\mu}^2(\nu_0,\nu_1) =
\llVert \varrho_0 - \varrho _1\rrVert
_{\dot
H^{-1}(\dx{\mu})}^2 = \inf_{h\in\dot H^1(\mu)} \biggl\{\int
\llvert \nabla h\rrvert ^2 \,\dx{\mu} \dvtx  L h = \varrho_0 -
\varrho _1 \biggr\} .
\]
In fact, it is possible to define a whole class of weighted Wasserstein
type distances interpolating between the negative Sobolev norm and the
Wasserstein distance. Theses transports were introduced in \cite
{Dolbeault2009}.
\end{rem}

\subsection{Approximation of the local measures \texorpdfstring{$\mu
_i$}{mu i}}\label{sapproxmd}

In this subsection, we show that it is sufficient to consider only the
mean-difference w.r.t. some auxiliary measures $\nu_i$ approximating
$\mu_i$ for $i=1,\dots, M$. More precisely, the next lemma shows that
there are nice measures $\nu_i$ which are close to the measures $\mu
_i$ in the sense of the mean-difference.
%
\begin{lem}[(Mean-difference of approximation)]\label{mdapproxcorchi2}
For $i=1,\dots, M$ let $\nu_i$ be a truncated Gaussian measure
centered around the local minimum $m_i$ with covariance matrix $\Sigma
_i= (\nabla^2 H(m_i))^{-1}$, more precisely
\begin{eqnarray}
\label{approxdefnu} \nu_i(\dx{x}) = \frac{1}{Z_{\nu_i}} e^{-\fracc{\Sigma
_i^{-1}[x-m_i]}{2\eps}}
\one_{E_i}(x) \,\dx{x} \nonumber
\\[-8pt]
\\[-8pt] \eqntext{\displaystyle \mbox {where } Z_{\nu_i} = \int
_{E_i} e^{-\fracc{\Sigma
_i^{-1}[x-m_i]}{2\eps}} \,\dx{x} ,}
\end{eqnarray}
where we write $A[x]:=  \langle x , Ax  \rangle$. The
restriction $E_i$ is given by
an ellipsoid
\begin{equation}
\label{approxdeftruncradius} E_i = \bigl\{x \in\R^n\dvtx  \bigl\llvert
\Sigma_i^{-\fraca
{1}{2}}(x-m_i)\bigr\rrvert \leq \sqrt{2
\eps} \omega(\eps) \bigr\} .
\end{equation}
Additionally, assume that $\mu_i$ satisfies $\PI(\varrho_i)$ with
$\varrho_i^{-1}=O(\eps)$.

Then the following estimate holds:
\begin{equation}
\label{mdapproxequclaim} \bigl(\Expect_{\nu_i}(f) - \Expect_{\mu_i}(f)
\bigr)^2 \leq O\bigl(\eps ^{\fraca{3} {2}}\omega^3(
\eps)\bigr) \int\llvert \nabla f\rrvert ^2 \,\dx{\mu} ,
\end{equation}
where the function $\omega(\eps)\dvtx \R^+ \to\R^+$ in \eqref
{approxdeftruncradius} and \eqref{mdapproxequclaim} is smooth and
monotone satisfying
\[
\label{mddefomega} \omega(\eps) \geq\llvert \log\eps\rrvert ^{\fraca{1}{2}} \qquad
\mbox{for } \varepsilon< 1.
\]
\end{lem}
The first step toward the proof of Lemma~\ref{mdapproxcorchi2} is
the following statement.
%
\begin{lem}\label{mdapproxlemcov}
Let $\nu_i$ be a probability measure satisfying $\nu_i \ll\mu_i$.
Moreover, if $\mu_i$ satisfies $\PI(\varrho_i)$ for some $\varrho
_i>0$, then the following estimate holds:
\begin{equation}
\label{mdapproxequcov} \bigl(\Expect_{\nu_i}(f) - \Expect_{\mu_i}(f)
\bigr)^2 \leq\frac
{1}{\varrho_i} \var_{\mu_i} \biggl({
\frac{\textup{d} \nu
_i}{\textup{d} \mu_i}} \biggr) \int \llvert \nabla f\rrvert ^2\,\dx{
\mu_i} .
\end{equation}
\end{lem}
\begin{pf}
The result is a consequence from the representation of the
mean-difference as a covariance. Therefore, we note that $\dx{\nu_i}
= {\frac{\textup{d} \nu_i}{\textup{d} \mu_i}} \,\dx{\mu_i}$ since
$\nu_i \ll\mu_i$ and
use the Cauchy--Schwarz inequality for the covariance
\begin{eqnarray*}
\bigl(\Expect_{\nu_i}(f) - \Expect_{\mu
_i}(f)
\bigr)^2 &=& 
\int f {\frac{\textup{d} \nu_i}{\textup{d} \mu_i}} \,\dx{
\mu_i} - \int f\,\dx{\mu_i} {\int{\frac{\textup{d} \nu_i}{\textup{d} \mu_i}}
\,\dx{\mu_i}}
\\
&=&\cov^2_{\mu_i} \biggl({\frac{\textup{d} \nu_i}{\textup{d} \mu
_i}},f \biggr)
\leq\var _{\mu_i} \biggl({\frac{\textup{d} \nu_i}{\textup{d} \mu_i}} \biggr)
\var_{\mu_i}(f).
\end{eqnarray*}
Using the fact that $\mu_i$ satisfies a PI results in \eqref
{mdapproxequcov}.\vadjust{\goodbreak}
\end{pf}
The above lemma tells us that we only need to construct $\nu_i$
approximating $\mu_i$ in variance for $i=1,\dots, M$. The following
lemma provides exactly this.
%
\begin{lem}[(Approximation in variance)]\label{mdapproxchi2}
Let the measures $\nu_i$ be given by \eqref{approxdefnu}. Then the
partition sum $Z_{\nu_i}$ satisfies for $\eps$ small enough
\begin{equation}
\label{normalizationnu} Z_{\nu_i} = (2\pi\eps)^{\fraca{n}{2}} \sqrt{\det
\Sigma_i} \bigl(1+ O\bigl(\sqrt{\eps} \omega^3(\eps)
\bigr)\bigr).
\end{equation}
Additionally, $\nu_i$ approximates $\mu_i$ in variance, that is,
\begin{equation}
\label{mdapproxequvarnumu} \var_{\mu_i} \biggl({\frac{\textup{d} \nu_i}{\textup{d} \mu
_i}} \biggr) = O
\bigl(\sqrt{\eps} \omega^3(\eps)\bigr) .
\end{equation}
\end{lem}
\begin{pf}
The proof of \eqref{normalizationnu} reduces to an estimate of a
Gaussian integral on the complementary domain $\R^{n}\setminus E_i$.
We deduced this estimate already in the proof of Lemma~\ref
{padmisibleverification}. By the same argument, we deduce
\[
Z_{\nu_i} = (2\pi\eps)^{\fraca{n}{2}} \sqrt{\det\Sigma_i}
\bigl(1+O\bigl(\sqrt{\eps} \omega^3(\eps)\bigr) \bigr) .
\]
Since $\mu_i$ comes from the restriction to an admissible partition
according to Definition~\ref{defpartition}
\begin{equation}
\label{normalizationcomp} Z_{\mu_i} = Z_i Z_\mu\exp
\biggl(\frac{H(m_i)}{\eps} \biggr) \stackrel {\mathclap{\eqref{ePartSumAdmPart}}}
{=} Z_{\nu_i} \bigl(1+O\bigl(\sqrt {\eps} \omega^3(\eps)
\bigr) \bigr).
\end{equation}
The relative density of $\nu_i$ w.r.t. $\mu_i$ can be estimated by
Taylor expanding $H$ around $m_i$. By the definition of $\nu_i$ given
in \eqref{approxdefnu}, we obtain that $\Sigma_i^{-1}[y-m_i]-H_i(y)
= O(\llvert  y-m_i\rrvert ^3)$. This observation together
with \eqref
{normalizationcomp} leads to
\begin{eqnarray*}
{\frac{\textup{d} \nu_i}{\textup{d} \mu_i}}(y) &=& \frac{Z_{\mu
_i}}{Z_{\nu_i}} e^{-\fracc{\Sigma_i^{-1}[y-m_i]}{2\eps}  + \fracc{H_i(y)}{2\eps} }
\one_{E_i}(y)= \frac{Z_{\mu_i}}{Z_{\nu_i}} e^{\fraca{O(\llvert
y-m_i\rrvert ^3)}{\eps}}
\one_{E_i}(y)
\\
&=& 1+O\bigl(\sqrt{\eps} \omega^3(\eps)\bigr).
\end{eqnarray*}
Now, the conclusion directly follows from the definition of the variance
\begin{eqnarray*}
\var_{\mu_i} \biggl({\frac{\textup{d} \nu
_i}{\textup{d} \mu_i}} \biggr) &=& \int
_{E_i} \biggl({\frac{\textup{d} \nu_i}{\textup{d} \mu
_i}} \biggr)^2\,
\dx{\mu_i}- \biggl(\int{\frac{\textup{d} \nu_i}{\textup{d} \mu_i}} \,\dx{\mu
_i} \biggr)^2
\\
&=& \int_{E_i} 1+O\bigl(\sqrt{\eps} \omega^3(
\eps)\bigr) \,\dx{\mu_i} - \biggl(\int_{E_i} \,
\dx{\nu_i} \biggr)^2
\\
&\leq& 1 +O\bigl(\sqrt{\eps} \omega^3(\eps)\bigr) - 1 = O\bigl(
\sqrt{\eps} \omega^3(\eps)\bigr) .
\end{eqnarray*}
\upqed\end{pf}
\begin{pf*}{Proof of Lemma~\ref{mdapproxcorchi2}}
A combination of Lemma~\ref{mdapproxlemcov} and Lemma~\ref
{mdapproxchi2} together with the assumption $\varrho_i^{-1}=O(\eps
)$ immediately reveals
\[
\bigl(\Expect_{\nu_i}(f) - \Expect_{\mu_i}(f)
\bigr)^2 \overset {\mathclap{\eqref{mdapproxequcov}, \eqref{mdapproxequvarnumu}}}
{\leq} O\bigl(\eps^{\fraca{3}{2}} \omega^3(\eps)\bigr) \int\llvert
\nabla f\rrvert ^2 \,\dx {\mu_i} .
\]
\upqed\end{pf*}

\subsection{Affine transport interpolation}\label{safftransp}

The aim of this section is to estimate $ (\Expect_{\nu_i}(f) -
\Expect_{\nu_j}(f)  )^2$ with the help of the weighted transport
distance $\cT_\mu(\nu_i, \nu_j)$ introduced in Section~\ref{sectranspmd} and is formulated in Lemma~\ref{ptranspcostapprox}.
For the proof of Lemma~\ref{ptranspcostapprox}, we construct an
explicit transport interpolation between $\nu_i$ and $\nu_j$
w.r.t. the measure $\mu$. We start with a class of possible transport
interpolations and optimize the weighted transport cost in this class.

Let us state the main idea of this optimization procedure. Therefore,
we recall that the measures $\nu_i$ and $\nu_j$ are truncated
Gaussians by the approximation we have done in the previous
Section~\ref{sapproxmd}. Hence, the measures $\nu_i$ and $\nu_j$
are characterized by their mean and covariance matrix. We will choose
the transport interpolation (cf. Section~\ref{secmdaffdef}) such
that the push forward measures $\nu_s:= (\Phi_s)_\sharp\nu_0$ are
again truncated Gaussians. Hence, it is sufficient to optimize among
all paths $\gamma$ connecting the minima $m_i$ and $m_j$ and all
covariance matrices interpolating between $\Sigma_i$ and $\Sigma_j$.

\subsubsection{Definition of regular affine transport
interpolations}\label{secmdaffdef}
Let us state in this section the class of transport interpolation among
we want to optimize the weighted transport cost.
%
\begin{defn}[(Affine transport interpolations)]
Assume that the measures $\nu_i$ and $\nu_j$ are given by Lemma~\ref
{mdapproxcorchi2}. In detail, $\nu_i = \cN(m_i , \eps^{-1}\Sigma
_i)\llcorner E_i $ and $\nu_j = \cN(m_j , \eps^{-1}\Sigma
_j)\llcorner E_j$ are truncated Gaussians centered in $m_i$ and $m_j$
with covariance matrices $\eps^{-1}\Sigma_i$ and $\eps^{-1}\Sigma
_j$. The restriction $E_i$ and $E_j$ are given for $l=1,\dots,M$ by
the ellipsoids
\[
E_l := \bigl\{x \in\R^n\dvtx  \bigl\llvert
\Sigma_l^{-\fraca
{1}{2}}(x-m_l)\bigr\rrvert \leq \sqrt{2
\eps} \omega(\eps) \bigr\}   \qquad\mbox{where } \omega(\eps) \geq\llvert \log\eps
\rrvert ^{\fraca{1}{2}} .
\]
A transport interpolation ${\Phi_s}$ between $\nu_i$ and $\nu_j$ is
called \emph{affine transport interpolation} if there exists:
\begin{itemize}
\item an interpolation path $(\gamma_s)_{s\in[0, T]}$ between
$m_i=\gamma_0$ and $m_j=\gamma_T$ satisfying
\begin{equation}
\label{edefpropaffinetranspinterpol} \gamma= (\gamma_s)_{s\in[0,T]} \in
C^2\bigl([0,T],\R^n\bigr) \quad\mbox {and}\quad\forall s
\in[0,T]\dvtx  \dot\gamma_s\in S^{n-1},
\end{equation}
\item an interpolation path $(\Sigma_s)_{s\in[0,T]}$ of covariance
matrices between $\Sigma_i$ and $\Sigma_j$ satisfying
\[
\Sigma= (\Sigma_s)_{s\in[0,T]} \in C^2\bigl([0,T],
\R^{n\times n}_{\sym
,+}\bigr) ,\qquad \Sigma_0 =
\Sigma_i  \mbox{ and } \Sigma _T=
\Sigma_j,
\]
\end{itemize}
such that the transport interpolation $(\Phi_s)_{s\in[0,T]}$ is given by
\begin{equation}
\label{mdafftranspdef} \Phi_s(x) = \Sigma_s^{\fraca{1}{2}}
\Sigma_0^{-\fraca
{1}{2}}(x-m_0)+\gamma_s .
\end{equation}
\end{defn}

Since the cost density $\cA$ given by \eqref{mdcostdensity} is
invariant under rescaling of time (cf. Remark~\ref
{mdremtimerescale}), one can always assume that the interpolation
path $\gamma_s$ is parameterized by arc-length. Hence, the condition
$\dot\gamma_s \in S^{n-1}$ [cf. \eqref
{edefpropaffinetranspinterpol}] is not restricting.

We want to emphasize that for an affine transport interpolation $(\Phi
_s)_{s\in[0,T]}$ the push forward measure $(\Phi_s)_\sharp\nu_0=\nu
_s$ is again a truncated Gaussian $\cN(\gamma_s ,\allowbreak  \eps^{-1}\Sigma
_s)\llcorner E_s$, where $E_s$ is the support of $\nu_s$ being again
an ellipsoid in $\R^n$ given by
\begin{equation}
\label{mdafftranspEs} E_s = \bigl\{x \in\R^n\dvtx  \bigl\llvert
\Sigma_s^{-\fraca
{1}{2}}(x-\gamma _s)\bigr\rrvert \leq
\sqrt{2\eps} \omega(\eps ) \bigr\}.
\end{equation}
Therewith, the partition sum of $\nu_s$ is given by [cf. \eqref
{normalizationnu}]
\begin{equation}
\label{normalizationnus} Z_{\nu_s} = (2\pi\eps)^{\fraca{n}{2}} \sqrt{\det
\Sigma_s}\bigl(1 + O(\sqrt{\eps})\bigr) .
\end{equation}
By denoting $\sigma_s = \Sigma_s^{\fraca{1}{2}}$ and using the
definition \eqref{mdafftranspdef} of the affine transport
interpolation $(\Phi_s)_{s\in[0,T]}$, we arrive at the relations
\begin{eqnarray}
\dot\Phi_s(x) &=& \dot\sigma_s\sigma_0^{-1}(x-m_0)
+ \dot\gamma _s,
\nonumber
\\
\Phi_s^{-1}(y) &=& \sigma_0
\sigma_s^{-1}(y-\gamma_s)+m_0,
\nonumber
\\
\dot\Phi_s \circ\Phi_s^{-1}(y) &=& \dot
\sigma_s \sigma _s^{-1}(y-
\gamma_s)+\dot\gamma_s .
\nonumber
\end{eqnarray}
Among all possible affine transport interpolations, we are considering
only those satisfying the following regularity assumption.
%
\begin{assume}[(Regular affine transport interpolations)]\label
{assumeafftranspreg}
An affine transport interpolation $(\gamma_s, \Sigma_s)_{s \in
[0,T  ]}$ belongs to the class of \emph{regular affine transport
interpolations} if the length $T< T^\ast$ is bounded by some uniform
$T^\ast>0$ large enough. Further, for a uniform constant $c_\gamma>0$ holds
\begin{equation}
\label{assumeequafftranspregularpath} \inf \bigl\{r(x,y,z) \dvtx  x,y,z\in\gamma, x\ne y\ne z\ne x \bigr\}
\geq c_\gamma,
\end{equation}
where $r(x,y,z)$ denotes the radius of the unique circle through the
three distinct points $x,y$ and $z$. Furthermore, there exists a
uniform constant $C_\Sigma\geq1$ for which
\begin{equation}
\label{assumeequafftranspregular} C_\Sigma^{-1}\Id\leq\Sigma_s \leq
C_\Sigma\Id\quad\mbox{and}\quad\llVert \dot\Sigma_s
\rrVert \leq C_\Sigma.
\end{equation}
\end{assume}

\begin{figure}

\includegraphics{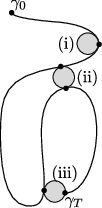}

\caption{Global radius of curvature.}\label{fig3}
\end{figure}

The infimum in condition \eqref{assumeequafftranspregularpath} is
called \emph{global radius of curvature} (cf. \cite
{vonderMosel2002}). It ensures that a small neighborhood of size $\frac
{c_\gamma}{2}$ around $\gamma$ is not self-intersecting, since the
infimum can only be attained for the following three cases (cp. Figure~\ref{fig3}):
\begin{longlist}[(iii)]
\item All three points in a minimizing sequence of \eqref
{assumeequafftranspregularpath} coalesce to a point at which the
radius of curvature is minimal.
\item Two points coalesce to a single point and the third converges to
another point, such that the both points are a pair of closest approach.
\item Two points coalesce to a single point and the third converges to
the starting or ending point of $\gamma$.
\end{longlist}

In the following calculations, there often occurs a multiplicative
error of the form $1+O(\sqrt{\eps} \omega^3(\eps))$.
Therefore, let us introduce for convenience the notation ``$\approx$''
meaning ``$=$'' up to the multiplicative error $1+O(\sqrt{\eps}
\omega^3(\eps))$. The symbols ``$\lesssim$'' and ``$\gtrsim$'' have
the analogous meaning.

%

Now, we can formulate the key ingredient for the proof of Theorem~\ref
{thmcoarsemd}, namely the estimation of the weighted transport
distance $\cT_\mu(\nu_i, \nu_j)$.
%
\begin{lem}\label{ptranspcostapprox}
Assume that $\nu_i$ and $\nu_j$ are given by Lemma~\ref
{mdapproxcorchi2}. Then the weighted transport distance $\cT_\mu
(\nu_i, \nu_j)$ can be estimated as
\begin{eqnarray}
\label{etransportcostsestimateapprox}%
\cT_{\mu}^2(\nu_i,
\nu_j) & =& \inf_{{\Phi_s}} \int \biggl(\int
_0^1 \bigl\llvert \dot\Phi_s \circ
\Phi_s^{-1}\bigr\rrvert {\frac
{\textup{d} \nu _{s}}{\textup{d} \mu}} \,\dx{s}
\biggr)^2 \,\dx {\mu}
\nonumber
\\
& \leq& \inf_{{ \Psi_s}} \int \biggl(\int_0^1
\bigl\llvert \dot \Psi _s \circ\Psi_s^{-1}
\bigr\rrvert {\frac{\textup{d} \nu
_{s}}{\textup{d} \mu}} \,\dx{s} \biggr)^2 \,\dx {\mu}
\nonumber
\\[-8pt]
\\[-8pt]
& \lesssim& \frac{Z_\mu}{(2\pi\eps)^{\fraca{n}{2}}} 2\pi\eps \biggl(\frac{ \sqrt{\llvert \det(\nabla^2 H(s_{i,j})\rrvert }}{\llvert \lambda^-(s_{i,j})\rrvert }+
\frac
{T(C_{\Sigma}) ^{\fracb{n-1}{2}} }{\sqrt{2\pi\eps}}e^{-\omega
^2(\eps)} \biggr)\nonumber\\
&&{}\times e^{\fraca{H(s_{i,j})}{\eps}} ,
\nonumber
\end{eqnarray}
where the infimum over $\Psi_s$ only considers regular affine
transport interpolations $\Psi_s$ in the sense of Assumption~\ref
{assumeafftranspreg}.

In particular, if we choose $\omega(\varepsilon) \geq|\log
\varepsilon|^{\fraca{1}{2}}$, which is enforced by Lemma~\ref
{mdapproxcorchi2}, we get the estimate
\begin{eqnarray}\label
{etransportcostsestimateapproxfinal}
\cT_{\mu}^2(\nu_i, \nu_j) &\leq&
\frac{Z_\mu}{(2\pi\eps)^{\fraca
{n}{2}}} \frac{ 2\pi\eps\sqrt{\llvert \det(\nabla^2
H(s_{i,j})\rrvert }}{\llvert \lambda^-(s_{i,j})\rrvert } e^{\fraca{H(s_{i,j})}{\eps
}}\nonumber
\\[-8pt]
\\[-8pt] &&{}\times\bigl(1+O\bigl(\sqrt{\eps}
\omega^3(\eps)\bigr) \bigr) .
\nonumber
\end{eqnarray}
\end{lem}
%
Before turning to the proof of Lemma~\ref{ptranspcostapprox}, we
want to anticipate the structure of the affine transport interpolation
$(\gamma, \Sigma)$ which realizes the desired estimate \eqref
{etransportcostsestimateapproxfinal}: Having a closer look at the
structure of the weighted transport distance $\cT_{\mu}^2(\nu_i, \nu
_j)$, it becomes heuristically clear that the mass should be
transported from $E_i$ to $E_j$ over the saddle point $s_{i,j}$ into
the direction of the eigenvector to the negative eigenvalue $\lambda
^-(s_{i,j})$ of $\nabla^2 H(s_{i,j})$. There, only the region around
the saddle gives the main contribution to the estimate \eqref
{etransportcostsestimateapproxfinal}. Then we only have one more
free parameter to choose for our affine transport interpolation
$(\gamma, \Sigma)$: It is the covariance structure $\Sigma_{\tau
^*}$ of the interpolating truncated Gaussian measure $\nu_{\tau^*}$
at the passage time $\tau^*$ at the saddle point $s_{
i,j}$. In the proof of Lemma~\ref{ptranspcostapprox} below, we will
see by an optimization procedure that the best $\Sigma_{\tau^*}$ is
given by $\Sigma_{\tau^\ast}^{-1}=\nabla^2 H(s_{i,j})$, restricted
to the stable subspace $\nabla^2 H(s_{i,j})$.

The proof of Lemma~\ref{ptranspcostapprox} presents the core of the
proof of the Eyring--Kramers formulas and consists of three steps
carried out in the following sections:
\begin{itemize}
\item
In Section~\ref{secmdafftubecoord}, we carry out some preparatory work:
We introduce tube coordinates on the support of the transport cost $\cA
$ given by \eqref{mdcostdensity} (cf. Lemma~\ref
{afftransplemcoordchange}), we deduce a pointwise estimate on the
transport cost $\cA$ and we give a rough a priori estimate on the
transport cost $\cA$.
\item In Section~\ref{secmdaffredsad}, we split the transport cost
into a transport cost around the saddle and the complement. We also
estimate the transport cost of the complement yielding the second
summand in the desired estimate \eqref{etransportcostsestimateapprox}.
\item In Section~\ref{secmdaffcostest}, we finally deduce a sharp
estimate of the transport cost around the saddle yielding the first
summand in the desired estimate \eqref{etransportcostsestimateapprox}.
\end{itemize}

%
\begin{figure}[b]

\includegraphics{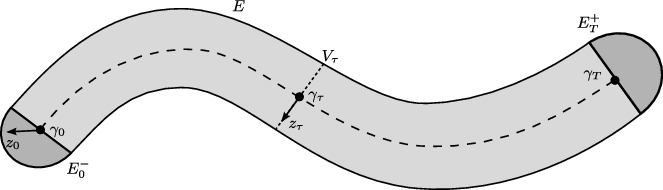}

\caption{The support of $\cA$ in tube coordinates.}\label{fig4}
\end{figure}

\subsubsection{Preparations and auxiliary estimates}\label
{secmdafftubecoord}
The main reason for making the regularity Assumption~\ref
{assumeafftranspreg} on affine transport interpolations is that we
can introduce tube coordinates around the path $\gamma$ as illustrated in Figure~\ref{fig4}. In these
coordinates, the calculation of the cost density $\cA$ given by \eqref
{mdcostdensity} becomes a lot handier.

We start with defining the caps $E_0^-$ and $E_T^+$ as
\[
E_0^- := \bigl\{x\in E_0\dvtx  \langle x-
\gamma_0 , \dot \gamma_0\rangle< 0 \bigr\} \quad
\mbox{and}\quad E_T^+ := \bigl\{ x\in E_T\dvtx  \langle x-
\gamma_T , \dot\gamma_T\rangle> 0 \bigr\}.
\]
The caps $E_0^-$ and $E_T^+$ have no contribution to the total cost but
unfortunately need some special treatment. Further, we define the
slices $V_s$ with $s\in[0,T]$
\[
V_s = \bigl\{x\in\vspan \{\dot\gamma_s
\}^\perp\dvtx  \bigl\llvert \Sigma_s^{-\fraca{1}{2}} x\bigr\rrvert
\leq \sqrt{2\eps} \omega(\eps) \bigr\} .
\]
In $\vspan V_s$, we can choose a basis $e_s^2,\dots,e_s^n$ smoothly
depending on the parameter~$s$. In particular, there exists a family
$(Q_s)_{s\in[0,T]} \in C^2([0,T],SO(n))$ satisfying the same
regularity assumption as the family $(\Sigma_\tau)_{\tau\in[0,T]}$
such that
\begin{equation}
\label{afftranspdefQrot} Q_s e^1 = \dot\gamma_s,
\qquad Q_s e^i = e_s^i \qquad
\mbox{for } i=2,\dots, n ,
\end{equation}
where $(e^1,\dots,e^n)$ is the canonical basis of $\R^n$.

Let use now define the tube $E$ as
\[
E= \bigcup_{s\in[0,T]} (\gamma_s +
V_s ).
\]
The support of the cost density $\cA$ given by \eqref
{mdcostdensity} is now given by
\begin{equation}
\label{echarsuppcA} \supp\cA= E_0^-\cup E \cup E_T^+ .
\end{equation}
By the definition \eqref{mdafftranspEs} of $E_s$ and the uniform
bound \eqref{assumeequafftranspregular} on $\Sigma_s$ holds
\begin{equation}
\label{afftranspdiamVtau} \operatorname{diam} V_s \leq2 \sqrt{2\eps C_\Sigma}
\omega(\eps) .
\end{equation}
Therewith, we find
\[
\label{afftransptubeneighborhood} \supp\cA\subset B_{2\sqrt{2\eps C_\Sigma} \omega(\eps)}\bigl( (\gamma
_\tau)_{\tau\in[0,T]}\bigr) := \bigl\{x \in\R^n\dvtx \llvert
x-\gamma _\tau\rrvert \leq2\sqrt{2\eps C_\Sigma} \omega(\eps)
\bigr\}.
\]
The assumption \eqref{mdafftranspEs} ensures that $B_{2\sqrt{2\eps
C_\Sigma} \omega(\eps)}( (\gamma_\tau)_{\tau\in[0,T]})$ is not
self-intersec\-ting for any $\eps$ small enough. The next lemma just
states that by changing to tube coordinates in $E$ one can
asymptotically neglect the Jacobian determinant $\det J$.
%
\begin{lem}[(Change of coordinates)]\label{afftransplemcoordchange}
The change of coordinates $(\tau,z)\mapsto x= \gamma_\tau+z_\tau$
with $z_\tau\in V_\tau$ satisfies for any function $\xi$ on $E$
\[
\int_E \xi(x) \,\dx{x} \approx\int_0^T
\int_{V_\tau} \xi(\gamma _\tau+ z_\tau) \,
\dx{z_\tau} \,\dx{\tau}.
\]
\end{lem}
\begin{pf}
We use the representation of the tube coordinates via \eqref
{afftranspdefQrot}. Therewith, it holds that $x=\gamma_\tau+Q_\tau
z$, where $z\in \{0 \}\times\R^{n-1}$. Then the
Jacobian $J$ of the
coordinate change $x\mapsto(\tau, Q_\tau z)$ is given by
\[
J=\bigl(\dot\gamma_\tau+ \dot Q_\tau z,
(Q_\tau)_2, \dots, (Q_\tau )_n
\bigr)\in\R^{n\times n},
\]
where $(Q_\tau)_i$ denotes the $i$th column of $Q_\tau$. By the
definition \eqref{afftranspdefQrot} of $Q_\tau$ follows $\dot
\gamma_\tau= (Q_\tau)_1$. Hence, we have the representation $J =
Q_\tau+ \dot Q_\tau z \otimes e_1$. The determinant of $J$ is then
given by
\[
\det (Q_\tau+ \dot Q_\tau z \otimes e_1 ) =
\underbrace {\det (Q_\tau)}_{=1} \det \bigl(\Id+
\bigl(Q_\tau^{\top} \dot Q_\tau z\bigr)\otimes
e_1 \bigr) = 1+ \bigl(Q_\tau^{\top} \dot
Q_\tau z \bigr)_{1} .
\]
By Assumption~\ref{assumeafftranspreg} holds $\llVert \dot
Q_\tau \rrVert \leq C_\Sigma$ implying $ (Q_\tau^{\top}
\dot Q_\tau z )_{1,1}=O(z)$. Since $Q_\tau z\in V_\tau$, we get
$O(z)=O(\sqrt{\eps
} \omega(\eps))$ by \eqref{afftranspdiamVtau}. Hence, we get
\[
\det J = 1 + O\bigl(\sqrt{\eps} \omega(\eps)\bigr),
\]
which concludes the proof.
\end{pf}

An important tool is the following auxiliary estimate.
%
\begin{lem}[(Pointwise estimate of the cost-density $\cA$)]\label
{affinelemdensitytube}
For $x\in\supp\cA$, we define
\begin{equation}
\label{affinelemdeftau} \tau=\argmin\limits_{s\in[0,T]} \llvert x-\gamma_s
\rrvert \quad \mbox{and} \quad z_\tau:= x - \gamma_\tau.
\end{equation}
Then the following estimate holds:
\begin{eqnarray}
\label{afftranspequdensitytube} \cA(x) &\lesssim&(2\pi\eps)^{- \fracb{n-1}{2}} \sqrt{\det
\nolimits_{1,1} \bigl(Q_\tau^\top\tilde
\Sigma_\tau^{-1}Q_\tau \bigr)} e^{-\fracc{\tilde\Sigma_\tau^{-1}[z_\tau]}{2\eps}}\nonumber
\\[-8pt]
\\[-8pt] &=:&
P_\tau e^{-\fracc{\tilde\Sigma_\tau^{-1}[z_\tau]}{2\eps}} ,
\nonumber
\end{eqnarray}
where $Q_\tau$ is defined in $\eqref{afftranspdefQrot}$ and $\tilde
\Sigma_\tau^{-1}$ is given by
\begin{equation}
\label{edeftildesigmatau} \tilde\Sigma_\tau^{-1} =
\Sigma_\tau^{-1} - \frac{1}{\Sigma_\tau
^{-1}[\dot\gamma_\tau]} \Sigma_\tau^{-1}
\dot\gamma_\tau\otimes \Sigma_\tau^{-1}\dot
\gamma_\tau.
\end{equation}
Further, $\det_{1,1} A$ denotes the determinant of the matrix obtained
from $A$ removing the first row and column.
\end{lem}
%
%
\begin{rem}\label{prempointwiseestimatecA}
With a little bit of additionally work, one could show that \eqref
{afftranspequdensitytube} holds with ``$\approx$'' instead of
``$\lesssim$.'' It follows from \eqref{edeftildesigmatau} that
the matrix $\tilde\Sigma_\tau^{-1}$ is positive definite. Hence,
$\cA$ is an $\R^{n-1}$-dimensional Gaussian on the slice $\gamma
_\tau+ V_\tau$ up to approximation errors.
\end{rem}
\begin{pf*}{Proof of Lemma~\ref{affinelemdensitytube}}
By the regularity Assumption~\ref{assumeafftranspreg} on the
transport interpolation, we find that for all $x\in\supp\cA$ holds uniformly
\[
I_T(x):= \{s\dvtx  E_s \ni x \} \qquad \mbox{satisfies } \cH
^{1}\bigl(I_T(x)\bigr) = O\bigl(\sqrt{\eps} \omega(\eps)
\bigr) .
\]
This allows us to linearize the transport interpolation around $\tau$
given in \eqref{affinelemdeftau}. It holds for $s$ such that $x\in E_s$
\begin{eqnarray}
\label{mdcostestlinearizeerror} \Sigma_s^{-1}[x-\gamma_s]
&=&\Sigma_{\tau}^{-1}[\gamma_\tau+
z_\tau- \gamma_{s}] + O\bigl(\eps^{\fraca{3} {2}}
\omega^3(\eps)\bigr)
\nonumber
\\[-8pt]
\\[-8pt]
&= &\Sigma_{\tau}^{-1}\bigl[(\tau-s)\dot\gamma_\tau+
z_\tau\bigr] + O\bigl(\eps ^{\fraca{3} {2}}\omega^3(
\eps)\bigr).
\nonumber
\end{eqnarray}
For similar reasons, we can linearize the determinant $\det\Sigma_s$
and have $\det\Sigma_s = \det\Sigma_{\tau} + O(\sqrt{\eps}
\omega(\eps))$. Finally, we have the following bound on the transport speed:
\begin{eqnarray}
\label{mdcostestspeed} \bigl\llvert \dot\Phi_s \circ\Phi_s^{-1}(x)
\bigr\rrvert \one _{E_s}(x) &=& \bigl\llvert \dot \sigma_s
\sigma_s^{-1}(x-\gamma _s)+\dot
\gamma_s \bigr\rrvert \one_{E_s}(x)
\nonumber
\\
&\leq& \bigl(\bigl\llvert \dot\sigma_s \sigma_s^{-1}(x-
\gamma _s)\bigr\rrvert +\llvert \dot\gamma_s\rrvert
\bigr) \one _{E_s}(x)
\\
&\leq& \bigl(C_\Sigma\llvert x-\gamma_s\rrvert + 1 \bigr)
\one_{E_s}(x) = \bigl(1+ O\bigl(\sqrt{\eps} \omega(\eps)\bigr) \bigr)
\one_{E_s}(x).
\nonumber %
\end{eqnarray}
Let us first consider the case $x\in E$.
We use \eqref{normalizationnus}, \eqref{mdcostestlinearizeerror}
and \eqref{mdcostestspeed} to arrive with $x=\gamma_\tau+z_\tau$
where $z_\tau\in V_\tau$ at
\begin{eqnarray*}
\cA(x) &=& \int_{I_T(x)} \bigl\llvert \dot\Phi_s
\circ\Phi _s^{-1}(x)\bigr\rrvert \frac{1}{Z_{\nu_s}} \exp
\biggl(-\frac{1}{2\eps}\Sigma _s^{-1}[x-
\gamma_s] \biggr) \one_{E_s}(x) \,\dx{s}
\\
&\leq& \frac{1}{(2\pi\eps)^{\fraca{n}{2}}} \int_{I_T(x)} \frac
{1+O(\sqrt{\eps} \omega(\eps))}{\sqrt{\det\Sigma_s}}
\exp \biggl(-\frac{1}{2\eps}\Sigma_s^{-1}[x-
\gamma_s] \biggr) \,\dx{s}
\\
&\lesssim&\frac{1}{(2\pi\eps)^{\fraca{n}{2}}\sqrt{\det\Sigma
_\tau}}\int_\R\exp \biggl(-
\frac{1}{2\eps}\Sigma_{\tau
}^{-1}\bigl[(\tau-s)\dot
\gamma_\tau+ z_\tau\bigr] \biggr)\,\dx{s}
\\
&=&\frac{\sqrt{\det\Sigma_\tau^{-1}}}{(2\pi\eps)^{\fraca{n}{2}}} \frac{\sqrt{2\pi\eps}}{\sqrt{\Sigma_{\tau}^{-1}[\dot\gamma
_{\tau}]}} \exp \biggl(-
\frac{1}{2\eps}\tilde\Sigma_\tau ^{-1}[z_\tau]
\biggr) \bigl(1+O\bigl(\sqrt{\eps} \omega^3(\eps ) \bigr)\bigr) ,
\end{eqnarray*}
where the last step follows by an application of a partial Gaussian
integration (cf. Lemma~\ref{appendlempartGauss}). Finally, by using
the relation \eqref{appendtildeAdetident}, we get that
\[
\frac{\det\Sigma_{\tau}^{-1}}{\Sigma_{\tau}^{-1}[\dot\gamma
_{\tau}]} = \det_{1,1} \bigl(Q_\tau^\top
\tilde\Sigma_\tau ^{-1}Q_\tau \bigr),
\]
and conclude the hypothesis for this case.

Let us now consider the case $x\in E_0^- \cup E_T^+$. For convenience,
we only consider the case $x\in E_0^-$. By the definition of $E_0^-$
holds $\tau=0$. The integration domain $I_T(x)$ is now given by
\begin{equation}
\label{eintegrationdomainE-} I_T(x)=[0,s^\ast) \qquad \mbox{with }
s^\ast= O\bigl(\sqrt{\eps} \omega(\eps)\bigr).
\end{equation}
Therewith, we can estimate $\cA(x)$ in the same way as for $x\in E$
and conclude the proof.
\end{pf*}
We only need one more ingredient for the proof of Lemma~\ref
{ptranspcostapprox}. It is an a priori estimate on the cost density
$\cA$.
%
\begin{lem}[(A priori estimates for the cost density $\cA$)]\label
{afftranspremcostdens}
For $\cA$, it holds:
\begin{eqnarray}
\int\cA(x) \,\dx{x} & \lesssim& T   \quad\mbox{and} \label
{mdcostprobdenest}
\\
\cA(x) & \lesssim& \biggl(\frac{C_\Sigma}{2\pi\eps } \biggr)^{\fracb{n-1}{2}} \qquad
\mbox{for } x \in\supp\cA. \label
{mdcostdensest}
\end{eqnarray}
\end{lem}
\begin{pf}
Let us first consider the estimate \eqref{mdcostprobdenest}. It
follows from the characterization \eqref{echarsuppcA} of the
support of $\cA$ that
\begin{equation}
\label{eaufteilungvonintcA} \int\cA(x) \,\dx{x} = \int_E \cA(x) \,
\dx{x} + \int_{E_0^- \cup
E_T^+} \cA(x) \,\dx{x}.
\end{equation}
Now, we estimate the first term on the right-hand side of the last
identity. Using the change to tube coordinates of Lemma~\ref
{afftransplemcoordchange} and noting that the upper bound \eqref
{afftranspequdensitytube} is a $(n-1)$-dimensional Gaussian density
on $V_\tau$ for $\tau\in[0,T]$, we can easily infer that
\[
\int_E \cA(x) \,\dx{x} \lesssim\llvert \gamma\rrvert =T.
\]
Let us turn to the second term on the right-hand side of \eqref
{eaufteilungvonintcA}. For convenience, we only consider the
integral w.r.t. the cap $E_0^-$. It follows from \eqref
{mdcostestspeed} and \eqref{eintegrationdomainE-} that
\begin{eqnarray*}
\int_{E_0^-} \cA(x) \,\dx{x} & \lesssim& \int
_{E_0^-} \int_0^1 \nu
_s(x) \,\dx{s} \,\dx{x}= \int_0^{s^\ast}
\int_{E_0^-} \nu_s(x) \,\dx{x} \,\dx{s}
\\
& \lesssim&\int_0^{s^\ast} \int
\nu_s(x) \,\dx{x} \,\dx{s} = s^{\ast} = O\bigl(\sqrt{\eps}
\omega(\eps)\bigr),
\end{eqnarray*}
which yields the desired statement \eqref{mdcostprobdenest}.

Let us now consider the estimate (\ref{mdcostdensest}). Note
by Remark~\ref{prempointwiseestimatecA} the matrix $\tilde\Sigma
_\tau^{-1}$ given by \eqref{edeftildesigmatau} is positive
definite and the matrix we subtract is also positive definite.
Therefore, it holds in the sense of quadratic forms
\[
0 < \tilde\Sigma_\tau^{-1} = \Sigma_\tau^{-1}
- \frac{1}{\Sigma
_\tau^{-1}[\dot\gamma_\tau]} \Sigma_\tau^{-1}\dot
\gamma_\tau \otimes\Sigma_\tau^{-1}\dot
\gamma_\tau\leq\Sigma_\tau^{-1}.
\]
Now, the uniform bound \eqref{assumeequafftranspregular} yields
\[
\sqrt{\det_{1,1} \bigl(Q_\tau^\top\tilde
\Sigma_\tau ^{-1}Q_\tau \bigr)} \leq
C_\Sigma^{\fracb{n-1}{2}}.
\]
Then the desired statement \eqref{mdcostdensest} follows directly
from the estimate \eqref{afftranspequdensitytube}.
\end{pf}

\subsubsection{Proof of Lemma~\texorpdfstring{\protect\ref{ptranspcostapprox}}{4.11}:
Reduction to neighborhood around the saddle}\label{secmdaffredsad}
Firstly, observe that from \eqref{mdcostdensest} follows the a priori estimate
\begin{equation}
\label{afftranspaprioriest} \frac{\cA^2(x)}{\mu(x)} \lesssim \biggl(\frac{C_\Sigma}{2\pi \eps
}
\biggr)^{n-1} Z_\mu e^{\fraca{1}{\eps} H(x)} .
\end{equation}
Hence, on an exponential scale, the leading order contribution to the
cost comes from neighborhoods of points where $H(x)$ is large.
Therefore, we want to make the set, where $H$ is comparable to its
value at the optimal connecting saddle $s_{i,j}$, as small as possible.
For this purpose, let us define the following set:
\begin{equation}
\label{mdafftranspdefXi} \Xi_{\gamma,\Sigma} := \bigl\{x\in\supp\cA\dvtx  H(x) \geq
H(s_{i,j}) - \eps\omega^2(\eps) \bigr\} .
\end{equation}
Therewith, we obtain by denoting the complement $\Xi_{\gamma,\Sigma
}^c := \supp\cA\setminus\Xi_{\gamma,\Sigma}$ the splitting
\[
\label{mdafftranspsplitcost} \cT^2_{\mu}(\nu_i,
\nu_j) \leq\int_{\Xi_{\gamma,\Sigma}} \frac
{\cA^2(x)}{\mu(x)} \,
\dx{x} + \int_{\Xi^c_{\gamma,\Sigma}} \frac
{\cA^2(x)}{\mu(x)} \,\dx{x}.
\]
The integral on $\Xi^c_{\gamma,\Sigma}$ can be estimated with the
a priori estimate \eqref{afftranspaprioriest} and Lemma~\ref
{afftranspremcostdens} as follows:
\begin{eqnarray}
\label{mdreducesaddleerror}%
\int_{\Xi_{\gamma,\Sigma}^c}\frac{\cA^2(x)}{\mu(x)} \,
\dx{x} &\overset{\mathclap{\eqref{mdafftranspdefXi}}} {\leq}& Z_\mu
e^{\fraca
{H(s_{i,j})}{\eps}-\omega^2(\eps)} \int_{\Xi_{\gamma,\Sigma
}^c}\cA^2(x) \,\dx{x}
\nonumber
\\
&\overset{\mathclap{\eqref{mdcostdensest}}} {\lesssim} &Z_\mu
e^{\fraca
{H(s_{i,j})}{\eps}-\omega^2(\eps)} \biggl(\frac{C_\Sigma}{2\pi
\eps} \biggr)^{\fracb{n-1}{2}} \int\cA(x) \,
\dx{x}
\\
&\overset{\mathclap{\eqref{mdcostprobdenest}}} {\lesssim} &Z_\mu
e^{\fraca
{H(s_{i,j})}{\eps}-\omega^2(\eps)} \biggl(\frac{C_\Sigma}{2\pi
\eps} \biggr)^{\fracb{n-1}{2}} T .
\nonumber
\end{eqnarray}
We observe that estimate \eqref{mdreducesaddleerror} is the second
summand in the desired bound \eqref{etransportcostsestimateapprox}.

\subsubsection{Proof of Lemma~\texorpdfstring{\protect\ref{ptranspcostapprox}}{4.11}: Cost
estimate around the saddle}\label{secmdaffcostest}
The aim of this subsection is to deduce the estimate
\begin{equation}
\label{ecostestimateatsaddle} \int_{\Xi_{\gamma,\Sigma}} \frac{\cA^2(x)}{\mu(x)} \,\dx{x}
\lesssim\frac{Z_\mu}{(2\pi\eps)^{\fraca{n}{2}}} e^{\fraca
{H(s_{i,j})}{\eps}} \frac{2\pi\eps\sqrt{\llvert \det(\nabla
^2 H(s_{i,j}))\rrvert }}{\llvert \lambda^-(s_{i,j})\rrvert } .
\end{equation}
Note that this estimate would yield the missing ingredient for the
verification of the desired estimate \eqref{etransportcostsestimateapprox}.

By the nondegeneracy Assumption~\ref{assumenondegenerate}, we can
assume that $\eps$ is small enough such that $E_0^-\cup E_T^+ \subset
\Xi_{\gamma,\Sigma}^c$. Hence, it follows that $\Xi_{\gamma,\Sigma
}\subset E$.
We claim that the transport interpolation $\Phi_s$ can be chosen such
that there exists a connected subinterval $I_T \subset[0,T]$ satisfying
\begin{equation}
\label{afftranspconveringxi} \Xi_{\gamma,\Sigma}\subset\bigcup_{s\in I_T}
(V_s +\gamma_s) \quad\mbox{and}\quad
\cH^1(I_T)=O\bigl(\sqrt{\eps} \omega(\eps)\bigr) .
\end{equation}
Indeed, the level set $ \{x\in\R^n\dvtx  H(x)\leq H(s_{i,j})-\eps
\omega^2(\eps) \}$ consists of at least two connected
components $M_i$
and $M_j$ such that $m_i\in M_i$ and $m_j\in M_j$. Further, it holds
\[
\dist(M_i,M_j) = \inf_{{x\in M_i, y\in M_j}}\llvert
x-y\rrvert = O\bigl(\sqrt {\eps} \omega(\eps)\bigr) ,
\]
which follows from expanding $H$ around $s_{i,j}$ in direction of the
eigenvector corresponding to the negative eigenvalue of $\nabla^2
H(s_{i,j})$. We can choose the path $\gamma$ in direction of this
eigenvector in a neighborhood of size $O(\sqrt{\eps} \omega(\eps
))$ around $s_{i,j}$, which shows \eqref{afftranspconveringxi}.

Combining the covering \eqref{afftranspconveringxi} and Lemma~\ref
{afftransplemcoordchange} yields the estimate
\begin{equation}
\label{afftranspreducedtransportsaddle} \int_{\Xi_{\gamma,\Sigma}} \frac{\cA^2(x)}{\mu(x)} \,\dx{x}
\lesssim\int_{I_T} \int_{V_\tau}
\frac{\cA^2(\gamma_\tau+
z_\tau)}{\mu(\gamma_\tau+ z_\tau)} \,\dx{z_\tau} \,\dx{\tau} .
\end{equation}
Recalling the definition \eqref{afftranspdefQrot} of the family of
rotations $(Q_\tau)_{\tau\in[0,T]}$, it holds that $z_\tau=Q_\tau
z$ with $z\in \{0 \}\times\R^{n-1}$. Hence, the following
relation holds:
\begin{eqnarray}
\label{eestimationsaddletubecoordinates} &&\int_{I_T} \int_{V_\tau}
\frac{\cA^2(\gamma_\tau+z_\tau)}{\mu
(\gamma_\tau+z_\tau)} \,\dx{z_\tau} \,\dx{\tau} \nonumber
\\[-8pt]
\\[-8pt]&&\qquad = \int
_{ \{
0 \}\times\R^{n-1}} \int_{I_T} \one_{V_\tau}(Q_\tau
z) \frac
{\cA
^2(\gamma_\tau+Q_\tau z)}{\mu(\gamma_\tau+Q_\tau z)} \,\dx{\tau} \,\dx{z} .
\nonumber
\end{eqnarray}
The next step is to rewrite $H(\gamma_\tau+Q_\tau z)$. We assume,
that $\gamma$ actually passes the saddle $s_{i,j}$ at time $\tau^\ast
\in(0,T)$. Then, by the reason that $\llvert  z_\tau\rrvert
=O(\sqrt{\eps}
\omega(\eps))$ for $z_\tau\in V_\tau$ and the global nondegeneracy
assumption \eqref{emorse}, we can Taylor expand $H(\gamma_{\tau
}+z_{\tau})$ around $s_{i,j}=\gamma_{\tau^\ast}$ for $\tau\in I_T$
and $z_\tau= Q_\tau z \in V_\tau$. More precisely, we get
\begin{eqnarray*}
&&H(\gamma_\tau+ Q_\tau z)-H(s_{i,j})
\\
&&\qquad=\tfrac{1}{2}\nabla^2 H(s_{i,j})[
\gamma_\tau+ Q_\tau z - s_{i,j}] + O\bigl(\llvert
\gamma_\tau+Q_\tau z - s_{i,j}\rrvert
^3\bigr)
\\
&&\qquad= \tfrac{1}{2} \nabla^2 H(s_{i,j})[
\gamma_\tau- \gamma _{\tau^*}] + \tfrac{1}{2}
\nabla^2 H(s_{i,j})[Q_\tau z]
\\
&&\quad\qquad{}+ \bigl\langle Q_\tau z , \nabla^2
H(s_{i,j}) (\gamma_{\tau }-\gamma_{\tau^*}) \bigr\rangle
+ O\bigl(\llvert \gamma_\tau+ Q_\tau z -
\gamma_{\tau^*}\rrvert ^3\bigr).
\end{eqnarray*}
Now, further expanding $\gamma_\tau$ and $Q_\tau$ in $\tau$ leads to
\begin{eqnarray*}
\gamma_\tau&=& \gamma_{\tau^\ast} + \dot\gamma_{\tau^\ast}
\bigl(\tau- \tau^{\ast}\bigr) + O\bigl(\bigl\llvert \tau-
\tau^\ast\bigr\rrvert \bigr)   \quad\mbox{and}
\\
Q_\tau z &=& Q_{\tau^\ast} z + O\bigl(\bigl\llvert \tau-
\tau^\ast\bigr\rrvert \llvert z\rrvert \bigr) .
\end{eqnarray*}
For the expansion of $H$, we arrive at the identity
\begin{eqnarray*}
& &H(\gamma_\tau+ Q_\tau z)-H(s_{i,j})
\\
&&\qquad=\tfrac{1}{2} \nabla^2 H(s_{i,j})\bigl[\dot
\gamma_{\tau^\ast
}\bigl(\tau-\tau^\ast\bigr) +O\bigl(\bigl\llvert
\tau-\tau^\ast\bigr\rrvert ^2\bigr)\bigr]\\
&&\qquad \quad {} +
\tfrac{1}{2} \nabla^2 H(s_{i,j})\bigl[Q_{\tau^\ast}
z + O\bigl(\bigl\llvert \tau-\tau ^*\bigr\rrvert \llvert z\rrvert \bigr)\bigr]
\\
&&\quad\qquad{}+ \bigl\langle Q_{\tau^*} z + O\bigl(\bigl\llvert \tau -
\tau^*\bigr\rrvert \llvert z\rrvert \bigr) , \nabla^2
H(s_{i,j}) \bigl(\dot\gamma_{\tau^*}\bigl(\tau-\tau
^\ast\bigr) +O\bigl(\bigl\llvert \tau-\tau^*\bigr\rrvert ^2
\bigr) \bigr) \bigr\rangle
\\
&&\quad\qquad{}+ O\bigl(\llvert \gamma_\tau+ Q_\tau z -
\gamma _{\tau ^*}\rrvert ^3\bigr)
\\
&&\qquad= \tfrac{1}{2} \nabla^2 H(s_{i,j})[\dot
\gamma_{\tau^\ast
}] \bigl(\tau-\tau^\ast\bigr)^2 +
\tfrac{1}{2} \nabla^2 H(s_{i,j})[ Q_{\tau
^\ast}
z]
\\
&&\quad\qquad{}+ \bigl\langle Q_{\tau^\ast}z , \nabla^2
H(s_{i,j}) \dot \gamma_{\tau^\ast} \bigr\rangle\bigl(\tau- \tau
^\ast\bigr)
\\
&&\quad\qquad{}+ O\bigl(\bigl\llvert \tau-\tau^\ast\bigr\rrvert
^3,\llvert z\rrvert \bigl\llvert \tau -\tau^\ast\bigr
\rrvert ^2,\llvert z\rrvert ^2\bigl\llvert \tau-
\tau^\ast\bigr\rrvert , \llvert z\rrvert ^3\bigr) .
\end{eqnarray*}
Using $\llvert \tau-\tau^\ast\rrvert =O(\sqrt{\eps}
\omega(\eps))$ and
$\llvert  z\rrvert  = O(\sqrt{\eps} \omega(\eps))$, we
obtain for the error
the estimate
\begin{eqnarray*}
&&O\bigl(\bigl\llvert \tau-\tau^\ast\bigr\rrvert ^3,\llvert
z\rrvert \bigl\llvert \tau-\tau^\ast\bigr\rrvert ^2,\llvert
z\rrvert ^2\bigl\llvert \tau-\tau^\ast\bigr\rrvert ,
\llvert z\rrvert ^3\bigr) \\
&&\qquad = O\bigl(\eps^{\fraca{3}{2}}
\omega^3(\eps)\bigr) .
\end{eqnarray*}
The term $ \langle Q_{\tau^\ast}z , \nabla^2 H(s_{i,j}) \dot
\gamma _{\tau^\ast}  \rangle(\tau- \tau^\ast)$ in the
expansion of $H$ has no
sign and has to vanish. This is only the case, if we choose $\dot
\gamma_{\tau^\ast}$ as an eigenvector of $\nabla^2 H(s_{i,j})$ to
the negative eigenvalue $\lambda^-(s_{i,j})$, because then
\[
\label{mdafftranpconstdotgamma} \bigl\langle Q_{\tau^\ast} z , \nabla^2
H(s_{i,j}) \dot\gamma _{\tau^\ast } \bigr\rangle \bigl(\tau-
\tau^\ast\bigr) = \lambda ^-(s_{i,j}) \langle
Q_{\tau^\ast} z , \dot\gamma_{\tau^\ast
} \rangle = 0 .
\]
Additionally, by this choice of $\dot\gamma_{\tau^\ast}$ the
quadratic form $\nabla^2 H(s_{i,j})[\dot\gamma_{\tau^\ast}]$
evaluates to
\[
\nabla^2 H(s_{i,j})[\dot\gamma_{\tau^\ast}] =
\lambda^-(s_{i,j}) \llvert \dot\gamma_{\tau^\ast}\rrvert
^2 = \lambda^-(s_{i,j}) .
\]
Therefore, we deduced the desired rewriting of $H(\gamma_\tau+Q_\tau
z)$ as
\begin{eqnarray}
\label{erewritingHsaddle} H(\gamma_\tau+ Q_\tau z) &=&
H(s_{i,j}) - \bigl\llvert \lambda^-(s_{i,j})\bigr\rrvert
\bigl(\tau- \tau^\ast\bigr)^2 \nonumber
\\[-8pt]
\\[-8pt]&&{}+ \tfrac{1}{2}
\nabla^2 H(s_{i,j})[ Q_{\tau^\ast} z] + O\bigl(
\eps^{\fraca{3}{2}} \omega^3(\eps)\bigr).
\nonumber
\end{eqnarray}
From the regularity assumptions on the transport interpolation, we can
deduce that
\begin{eqnarray}
\tilde\Sigma_\tau^{-1}[Q_\tau z] &=
&\tilde\Sigma_{\tau
^*}^{-1}[Q_{\tau} z] + O\bigl(\bigl
\llvert \tau-\tau^*\bigr\rrvert \llvert z\rrvert ^2\bigr)
\nonumber
\\
&=& \tilde\Sigma_{\tau^*}^{-1}\bigl[Q_{\tau^*} z + O
\bigl(\bigl\llvert \tau -\tau ^*\bigr\rrvert \llvert z\rrvert \bigr)\bigr] + O
\bigl(\bigl\llvert \tau -\tau^*\bigr\rrvert \llvert z\rrvert ^2\bigr)
\nonumber
\\
&= &\tilde\Sigma_{\tau^\ast}^{-1}[Q_{\tau^\ast} z] + O\bigl(
\eps ^{\fraca{3}{2}}\omega^3(\eps)\bigr) .\nonumber
\end{eqnarray}
Then it follows easily from the definition \eqref
{afftranspequdensitytube} of $P_\tau$ that
\begin{equation}
\label{etauistaustar} P_\tau\approx P_{\tau^\ast}  \qquad\mbox{for }
\tau\in I_T .
\end{equation}
Applying the cost estimate \eqref{afftranspequdensitytube} of
Lemma~\ref{affinelemdensitytube}, the representation \eqref
{erewritingHsaddle} and the identity \eqref{etauistaustar}
yields the estimate for $(\gamma_\tau+ Q_\tau z) \in\Xi_{\gamma,
\Sigma}$
\begin{eqnarray}
\label{mdaffrepcostdenstiy} &&\qquad \frac{\cA^2(\gamma_\tau+Q_\tau z)}{\mu(\gamma_\tau+Q_\tau z)} \nonumber
\\[-8pt]
\\[-8pt]&&\qquad \qquad \lesssim Z_\mu e^{\fraca{H(s_{i,j})}{\eps}}
P_{\tau^\ast}^2 e^{-\fracc{(2\tilde\Sigma_{\tau^\ast}^{-1} -\nabla^2
H(s_{i,j}))[Q_{\tau^\ast} z]}{2\eps} - \fracc{\llvert \lambda
^-(s_{i,j})\rrvert  (\tau-\tau^\ast)^2}{2\eps}}.
\nonumber
\end{eqnarray}
The exponentials are densities of two Gaussian, if we put an additional
constraint on the transport interpolation. Namely, we postulate
\[
2\tilde\Sigma_{\tau^\ast}^{-1} - \nabla^2
H(s_{i,j}) > 0 \qquad  \mbox{on } \vspan{V_{\tau^\ast}}
\]
in the sense of quadratic forms. It holds that $\vspan{V_{\tau^\ast
}}=Q_{\tau^\ast}( \{0 \}\times\R^{n-1})=\vspan \{
\dot\gamma _{\tau^\ast} \}^\perp$ is stable subspace of
$\nabla^2 H(s_{i,j})$.
With these preliminary considerations, we finally are able to estimate
the right-hand side of \eqref{eestimationsaddletubecoordinates} as follows:
\begin{eqnarray}\label
{ecalcsaddlebeforeoptimizing}
&&\int_{ \{0 \}\times\R^{n-1}} \int_{I_T}
\one_{V_\tau
}(Q_\tau z) \frac{\cA^2(\gamma_\tau+Q_\tau z)}{\mu(\gamma_\tau+Q_\tau
z)} \,\dx{\tau} \,\dx{z}
\nonumber
\\
&&\qquad\overset{\mathclap{\eqref {mdaffrepcostdenstiy}}} {\lesssim}
Z_\mu e^{\fraca
{H(s_{i,j})}{\eps}}\int_{ \{0 \}\times\R^{n-1}}\int
_{I_T} P_{\tau
^\ast}^2 \nonumber\\
&&\hspace*{60pt}\qquad \quad {}\times e^{-\fracc{(2\tilde\Sigma_{\tau^\ast}^{-1} -\nabla^2
H(s_{i,j}))[Q_{\tau^\ast} z]}{2\eps} - \fracc{ \llvert \lambda
^-(s_{i,j})\rrvert  (\tau-\tau^\ast)^2}{2\eps}} \,\dx{
\tau} \,\dx{z}
\nonumber
\\[-8pt]
\\[-8pt]
&&\qquad\leq Z_\mu e^{\fraca{H(s_{i,j})}{\eps}} \frac{\sqrt {2\pi\eps}}{\sqrt{\llvert \lambda^-(s_{i,j})\rrvert }} \int
_{ \{0 \}\times
\R^{n-1}} P_{\tau^\ast}^2 e^{-\fracc{(2\tilde\Sigma_{\tau^\ast
}^{-1} -\nabla^2 H(s_{i,j}))[Q_{\tau^\ast} z]}{2\eps}} \,\dx
{z}
\nonumber\\
&&\qquad= Z_\mu e^{\fraca{H(s_{i,j})}{\eps}} \frac{\sqrt{2\pi
\eps}}{\sqrt{\llvert \lambda^-(s_{i,j})\rrvert }}
P_{\tau
^\ast}^2 \frac
{(2\pi\eps)^{\fracb{n-1}{2}}}{\sqrt{\det_{1,1} (Q_{\tau
^*}^\top(2\tilde\Sigma_{\tau^\ast}^{-1} -\nabla^2
H(s_{i,j}))Q_{\tau^*} )}}
\nonumber
\\
&&\qquad= \frac{Z_\mu}{(2\pi\eps)^{\fraca{n}{2}}} e^{\fraca
{H(s_{i,j})}{\eps}} \frac{2\pi\eps}{\sqrt{\llvert \lambda
^-(s_{i,j})\rrvert }} \underbrace{
\frac{\det_{1,1}
(Q_{\tau ^*}^\top\tilde\Sigma_{\tau^\ast}^{-1}Q_{\tau^*}
)}{\sqrt{\det
_{1,1} (Q_{\tau^*}^\top(2\tilde\Sigma_{\tau^\ast}^{-1}
-\nabla^2 H(s_{i,j}))Q_{\tau^*} )}}}_{\mathrm{to\ optimize!}} .
\nonumber %
\end{eqnarray}
The final step consists of optimizing the choice of\vspace*{1pt} $\tilde\Sigma
_{\tau^\ast}$. Let us use the notation $A=Q_{\tau^*}^\top\tilde
\Sigma_{\tau^\ast}^{-1}Q_{\tau^*}$ and $B=Q_{\tau^*}^\top
H(s_{i,j})Q_{\tau^*}$. Then the minimization problem has the structure
\begin{equation}
\label{ecalcsaddleminproblem} \inf_{A \in\R^{n\times n}_{\sym,+}} \biggl\{\frac{\det_{1,1}A}{
\sqrt{\det_{1,1} (2A-B )}} \dvtx  2 A -
B > 0 \mbox { on } \{0 \}\times\R^{n-1} \biggr\} .
\end{equation}
In the \hyperref[app]{Appendix}, we show in Lemma~\ref{coarselemmatrixopt} that the
optimal value of \eqref{ecalcsaddleminproblem} is attained at
$\tilde\Sigma_{\tau^\ast}^{-1}=\nabla^2 H(s_{i,j})$ restricted to
$V_{\tau^*}$. The optimal value is given by
\[
\frac{\det_{1,1}A}{ \sqrt{\det_{1,1} (2A-B )}} = \sqrt {\det_{1,1} \bigl(Q_{\tau^*}^\top
\nabla^2 H(s_{i,j})Q_{\tau^*}\bigr)}.
\]
Because $V_{\tau^\ast}$ is the stable subspace of $\nabla^2
H(s_{i,j})$, it holds
\begin{equation}
\label{eoptimizeddeterminant}\quad  \det_{1,1}\bigl(Q_{\tau^*}^\top
\nabla^2 H(s_{i,j})Q_{\tau^*}^\top\bigr) =
\frac{\det(\nabla^2 H(s_{i,j}))}{\lambda^-(s_{i,j})} = \frac{|\det
(\nabla^2 H(s_{i,j}))|}{|\lambda^-(s_{i,j})|} .
\end{equation}
The final step is a combination of \eqref
{afftranspreducedtransportsaddle}, \eqref
{eestimationsaddletubecoordinates}, \eqref
{ecalcsaddlebeforeoptimizing} and \eqref{eoptimizeddeterminant}
to obtain the desired estimate \eqref{ecostestimateatsaddle}. This
together with \eqref{mdreducesaddleerror} concludes \eqref
{etransportcostsestimateapprox} of Lemma~\ref{ptranspcostapprox}.

\subsubsection{Proof of Lemma~\texorpdfstring{\protect\ref{ptranspcostapprox}}{4.11}:
Total error
estimate}
For the verification of Lemma~\ref{ptranspcostapprox}, it is only
left to deduce the estimate \eqref
{etransportcostsestimateapproxfinal}. For that purpose, we analyze
the error terms in the estimate \eqref
{etransportcostsestimateapprox} that is,
\begin{eqnarray*}
&&\cT^2_{\mu}(\nu_i,\nu_j)\\&&\qquad \lesssim
\frac{Z_\mu}{(2\pi\eps)^{\fraca
{n}{2}}} e^{\fraca{H(s_{i,j})}{\eps}} 2\pi\eps \biggl(\underbrace{
\frac{ \sqrt{\llvert \det(\nabla^2
H(s_{i,j}))\rrvert }}{\llvert \lambda^-(s_{i,j})\rrvert
}}_{=O(1)}+\underbrace{\frac{T(C_{\Sigma
})^{\fracb{n-1}{2}} }{\sqrt{2\pi\eps}}e^{-\omega^2(\eps
)}}_{=O(\eps^{-\fraca{1}{2}}e^{-\omega^2(\eps)})}
\biggr) .
\end{eqnarray*}
By the choice of $\omega(\eps)\geq\llvert \log\eps\rrvert
^{\fraca{1}{2}}$,
enforced by Lemma~\ref{mdapproxcorchi2}, we see that
\[
O\bigl(\eps^{-\fraca{1}{2}}e^{-\omega^2(\eps)}\bigr) = O(\sqrt{\eps}).
\]
Recalling, that ``$\lesssim$'' means ``$\leq$'' up to a
multiplicative error of order $1+O(\sqrt{\eps} \omega^3(\eps))$ we
get the desired estimate \eqref{etransportcostsestimateapproxfinal}
\[
\cT^2_{\mu}(\nu_i,\nu_j)\lesssim
\frac{Z_\mu}{(2\pi\eps)^{\fraca
{n}{2}}} e^{\fraca{H(s_{i,j})}{\eps}} 2\pi\eps\frac{ \sqrt {\llvert \det(\nabla^2 H(s_{i,j}))\rrvert }}{\llvert
\lambda^-(s_{i,j})\rrvert } \bigl(1+O\bigl(
\sqrt{\eps} \omega^3(\eps)\bigr) \bigr). 
\]
%

\subsection{Proof of Theorem~\texorpdfstring{\protect\ref{thmcoarsemd}}{2.12}: Conclusion
of the
mean-difference estimate}\label{smdaffcomperror}

With the help of Lemma~\ref{mdapproxcorchi2} and Lemma~\ref
{ptranspcostapprox} the proof of Theorem~\ref{thmcoarsemd} is
straightforward. We can estimate the mean-differences w.r.t. to the
measure $\mu_i$ by introducing the means w.r.t. the
approximations $\nu_i$ and $\nu_j$
\begin{eqnarray*}
&&\bigl(\Expect_{\mu_i}(f) - \Expect_{\mu_j}(f)
\bigr)^2\\&&\qquad  = \bigl(\Expect_{\mu_i}(f) - \Expect_{\nu_i}(f)
+\Expect _{\nu
_i}(f) - \Expect_{\nu_j}(f) + \Expect_{\nu_j}(f)
- \Expect _{\mu_j}(f) \bigr)^2.
\end{eqnarray*}
We apply the Young inequality with a weight that is motivated by the
final total multiplicative error term $R(\eps)$ in Theorem~\ref
{thmcoarsemd}. More precisely,
\begin{eqnarray}
&&\bigl(\Expect_{\mu_i}(f) - \Expect_{\mu_j}(f)
\bigr)^2\nonumber\\
 &&\qquad \leq \bigl(1+\eps^{\fraca{1}{2}}\omega^3(\eps)
\bigr) \bigl(\Expect_{\nu
_i}(f)-\Expect_{\nu_j}(f)
\bigr)^2
\nonumber
\\
& & \qquad \quad {} + 2\bigl(1+{\eps}^{-\fraca{1}{2}}\omega^{-3}(\eps)\bigr) \bigl(
\bigl(\Expect_{\mu_i}(f) -\Expect_{\nu_i}(f) \bigr)^2+
\bigl(\Expect _{\mu_j}(f) -\Expect_{\nu_j}(f)
\bigr)^2 \bigr) .\nonumber %
\end{eqnarray}
Then the estimate \eqref{mdapproxequclaim} of Lemma~\ref
{mdapproxcorchi2} yields
\begin{eqnarray}
\label{mdapproxequ} \bigl(\Expect_{\mu_i}(f) - \Expect_{\mu_j}(f)
\bigr)^2 &\leq& \bigl(1+\sqrt {\eps} \omega^3(\eps)\bigr)
\bigl(\Expect_{\nu_i}(f)-\Expect_{\nu
_j}(f) \bigr)^2\nonumber
\\[-8pt]
\\[-8pt] &&{}+
O(\eps) \int\llvert \nabla f\rrvert ^2 \,\dx{\mu} ,
\nonumber
\end{eqnarray}
which justifies the statement, that the approximation only leads to
higher-order error terms in $\eps$.
An application of \eqref{mddefnestimate}\vadjust{\goodbreak} to the estimate \eqref
{mdapproxequ} transfers the mean-difference to the Dirichlet form
with the help of the weighted transport distance
\[
\bigl(\Expect_{\mu_i}(f) - \Expect_{\mu_j}(f)
\bigr)^2 \leq \bigl( \bigl(1+\sqrt{\eps} \omega^3(\eps)
\bigr) \cT^2_{\mu
}(\nu_i,\nu_j) +
O(\eps) \bigr) \int\llvert \nabla f\rrvert ^2 \,\dx{\mu} .
\]
The weighted transport distance $\cT_{\mu}(\nu_i,\nu_j)$ is
dominating the above estimate. Finally, we arrive at the estimate
\[
\bigl(\Expect_{\mu_i}(f) - \Expect_{\mu_j}(f)
\bigr)^2 \lesssim \cT ^2_{\mu}(
\nu_i,\nu_j) \int\llvert \nabla f\rrvert ^2
\,\dx {\mu} .
\]
Now, the Theorem~\ref{thmcoarsemd} follows directly from an
application of the estimate \eqref
{etransportcostsestimateapproxfinal} of Lemma~\ref
{ptranspcostapprox} and setting $\omega(\varepsilon) = \llvert
\log \varepsilon\rrvert ^{\fraca{1}{2}}$.

%
\begin{appendix}\label{app}
\section{Properties of the logarithmic mean \texorpdfstring{$\Lambda$}{Lambda}}\label{clogmean}

In this part of the \hyperref[app]{Appendix}, we collect some properties of the
logarithmic mean $\Lambda(\cdot,\cdot)$.
A more complete study can be found in \cite{Carlson1972}.

Let us first recall the definition of $\Lambda(\cdot,\cdot)\dvtx \R^+
\times\R^+ \to\R^+$
\begin{equation}
\label{logmeanintrep1} \Lambda(a,b) = \int_0^1
a^s b^{1-s} \,\dx{s}= \cases {\displaystyle
\frac{a-b}{\log a -\log b} ,&\quad$a\ne b $,
\cr
\displaystyle a ,&\quad$a=b$.}
\end{equation}
%
The equation \eqref{logmeanintrep1} justifies the statement, that
$\Lambda(\cdot,\cdot)$ is a mean, since one immediately recovers the
simple bounds $\min \{a,b \}\leq\Lambda(a,b)\leq\max
\{a,b \}$.
Moreover, two other immediate properties are:
\begin{itemize}
\item$\Lambda(\cdot, \cdot)$ is symmetric
\item$\Lambda(\cdot, \cdot)$ is homogeneous of degree one, that is,
for $\Lambda(\lambda a, \lambda b) = \lambda\Lambda(a,b)$ for
$\lambda>0$.
\end{itemize}
The derivatives of $\Lambda(\cdot, \cdot)$ are given by
straight-forward calculus
\[
\partial_a \Lambda(a,b) = \frac{1-\fraca{\Lambda(a,b)}{a}}{\log a-
\log b} >0 \quad\mbox{and}
\quad \partial_b \Lambda(a,b) = \frac{1- \fraca{\Lambda(a,b)}{b}}{\log b
- \log a} >0 .
\]
Hence, $\Lambda(\cdot,\cdot)$ is strictly monotone increasing in
both arguments.

The following result is almost classical and proven for instance in
\cite{Carlson1972}, Theorem~1, \cite{Mielke2011}, Appendix~A, and
\cite{Bhatia2008}.
%
\begin{lem}\label{logmeanlemgeologar}
The logarithmic mean can be bounded below by the geometric mean and
above by the arithmetic mean
\begin{equation}
\label{logmeangeologar} \sqrt{ab}\leq\Lambda(a,b) \leq\frac{a+b}{2} ,
\end{equation}
with equality if and only if $a=b$.
\end{lem}
%
The bounds in \eqref{logmeangeologar} are good, if $a$ is of the same
order as $b$, whereas the following bound is particularly good if
$\frac{a}{b}$ becomes very small or very large.
%
\begin{lem}\label{logmeanlemupperbound}
It holds for $p\in(0,1)$, the following bound:
\begin{equation}
\label{logmeanlemeupperbound} \frac{\Lambda(p,1-p)}{p(1-p)} < \min \biggl\{\frac{1}{p\log(\fraca
{1}{p})},
\frac{1}{(1-p)\log(\fracc{1}{1-p})} \biggr\} .
\end{equation}
\end{lem}
\begin{pf}
Let us first consider the case $0<p< \frac{1}{2}$. Then it is enough
to show that
\begin{equation}
\label{logmeanlemupperboundp1} \frac{\Lambda(p,1-p)}{p(1-p)} p\log\frac{1}{p} =
\frac{(1-2p)
\log(\fraca{1}{p})}{(1-p)\log(\fracb{1-p}{p})} \overset{!} {<} 1.
\end{equation}
This follows easily from the following lower bound on the denominator
\begin{eqnarray*}
(1-p)\log\frac{1-p}{p} &=& (1-2p)\log\frac{1}{p} + p\log
\frac{1}{p} - (1-p)\log\frac{1}{1-p}\\ &>& (1-2p)\log\frac{1}{p}
,
\end{eqnarray*}
since $p\log\frac{1}{p} - (1-p)\log\frac{1}{1-p} > 0$ for $0 < p <
\frac{1}{2}$. The case $\frac{1}{2}<p<1$ follows by symmetry under
the variable change $p\mapsto1-p$. It remains to check the case
$p=\frac{1}{2}$. The left-hand side of \eqref
{logmeanlemupperboundp1} evaluates for $p=\frac{1}{2}$ to
\[
\lim_{p\to\fraca{1}{2}} \frac{\Lambda(p,1-p)}{p(1-p)} p\log\frac
{1}{p} =
\log2 < 1 .
\]
\upqed\end{pf}
The logarithmic mean also occurs in the following optimization problem,
which appears in the proof of the optimality of the Eyring--Kramers
formula for the logarithmic Sobolev constant in one dimension
(cf. Section~\ref{smainopt}).
%
\begin{lem}\label{logmeanlemopt}
For $p\in(0,1)$ and $t\in(0,1)$, we define the function $h_p(t)$
according to
\begin{equation}
\label{logmeanlemdefhp} h_p(t) = \frac{ (\sqrt{\fraca{t}{p}}-\sqrt{\fracf
{1-t}{1-p}} )^2}{t \log(\fraca{t}{p}) + (1-t)\log(\fracf{1-t}{1-p})}.
\end{equation}
Then it holds
\begin{equation}
\label{logmeanlemeopt} \min_{t\in(0,1)} h_p(t) =
\frac{\Lambda(p,1-p)}{p(1-p)} .
\end{equation}
The minimum in \eqref{logmeanlemeopt} is attained for $t = 1-p$.
\end{lem}
\begin{pf}
Let us introduce the function $f_p\dvtx (0,1)\to\R$ and $g_p\dvtx (0,1)\to\R$
given by the nominator and denominator of $h_p$ in \eqref
{logmeanlemdefhp}, namely
\[
f_p(t) := \biggl(\sqrt{\frac{t}{p}}-\sqrt{\frac{1-t}{1-p}}
\biggr)^2 \quad\mbox{and}\quad g_p(t):=t \log
\frac{t}{p} + (1-t)\log \frac{1-t}{1-p}.
\]
It is easy to verify, that the following relations for the derivatives
hold true:
\begin{eqnarray} \label{logmeanlemoptp1}
f_p'(t) &=& \biggl(\sqrt{\frac{t}{p}}-\sqrt{
\frac
{1-t}{1-p}} \biggr) \biggl(\frac{1}{\sqrt{tp}} + \frac{1}{\sqrt {(1-p)(1-t)}}
\biggr), \nonumber\\
  g_p'(t) &=& \log\frac{t}{p} - \log
\frac
{1-t}{1-p},
\\
f_p''(t) &=& \sqrt{\frac{(1-t)t}{(1-p)p}}
\frac{1}{2(1-t)^2 t^2} > 0 , \qquad g_p''(t)
= \frac{1}{(1-t)t} > 0 .
\nonumber
\end{eqnarray}
Hence, both functions $f_p$ an $g_p$ are strictly convex and have a
unique minimum for $t=p$, where they are both zero. The derivative of
the quotient of $f_p$ and $g_p$ has the form
\begin{equation}
\label{logmeanlemoptp3} h_p'(t) := \biggl(\frac{f_p(t)}{g_p(t)}
\biggr)' = \frac{1}{g_p^2(t)} \bigl(f_p'(t)
g_p(t) - f_p(t) g_p'(t)
\bigr).
\end{equation}
The representation \eqref{logmeanlemoptp1} for $g_p'$ leads to
\begin{equation}
\label{logmeanlemoptp4} 
h_p'(t)
g_p^2(t) = \bigl(t f_p'(t) -
f_p(t) \bigr) \log\frac{t}{p} + \bigl((1-t)f_p'(t)+f_p(t)
\bigr) \log\frac{1-t}{1-p} .\hspace*{-30pt}
\end{equation}
Now, we can use \eqref{logmeanlemoptp1} for $f_p'$ to find
\begin{eqnarray}
\label{logmeanlemoptp5} %
&&tf_p'(t)-f_p(t)\nonumber\\
&&\qquad = \biggl(\sqrt{\frac{t}{p}}-\sqrt{\frac
{1-t}{1-p}} \biggr) \biggl(
\sqrt{\frac{t}{p}}+\frac{t}{\sqrt {(1-p)(1-t)}} - \sqrt{\frac{t}{p}} + \sqrt{
\frac
{1-t}{1-p}} \biggr)
\\
&&\qquad = \frac{1}{\sqrt{(1-p)(1-t)}} \biggl(\sqrt{\frac {t}{p}}-\sqrt {
\frac{1-t}{1-p}} \biggr) %
\nonumber
\end{eqnarray}
and likewise
\begin{equation}
\label{logmeanlemoptp6} (1-t)f_p'(t) + f_p(t) =
\frac{1}{\sqrt{tp}} \biggl(\sqrt{\frac
{t}{p}}-\sqrt{\frac{1-t}{1-p}}
\biggr) .
\end{equation}
Using \eqref{logmeanlemoptp5} and \eqref{logmeanlemoptp6}
in \eqref{logmeanlemoptp4} leads by \eqref{logmeanlemoptp3} to
\[
\label{logmeanlemoptp7} h_p'(t) = \frac{1}{g_p^2(t)}
\underbrace{ \biggl(\sqrt{\frac
{t}{p}}-\sqrt{\frac{1-t}{1-p}}
\biggr)}_{=:v_p(t)} \underbrace{ \biggl(\frac{\log(\fraca{t}{p})}{\sqrt{(1-p)(1-t)}} +
\frac{\log(\fracf
{1-t}{1-p})}{\sqrt{tp}} \biggr)}_{=:w_p(t)} .
\]
Since $g_p(p)=g_p'(p)=0$ and $g_p''(p)>0$, the function $\frac
{1}{g_p^2(t)}$ has a pole of order $4$ in $t=p$. Moreover, the function
$v_p(t)$ has a simple zero in $t=p$. We have to do some more
investigations for the function $w_p(t)$. First, we observe that
$w_p(t)$ can be rewritten as
\begin{eqnarray*}
w_p(t) &=& \underbrace{\frac{t-p}{\sqrt{(1-t)t (1-p)p}}}_{=: \hat
w_p(t)}
\\
&&{}\times\underbrace{ \biggl(\frac{\sqrt{tp} \log(\fraca
{t}{p})}{(t-p)} - \frac{\sqrt{(1-t)(1-p)} \log(\fracf
{1-t}{1-p})}{(p-t)}
\biggr)}_{:=\tilde w_p(t)} .
\end{eqnarray*}
The function $\tilde w_p(t)$ can be expressed in terms of the
logarithmic mean
\begin{equation}
\label{logmeanlemoptp8} \tilde w_p(t) = \frac{\sqrt{tp}}{\Lambda(t,p)} -
\frac{\sqrt {(1-t)(1-p)}}{\Lambda(1-t,1-p)}
\end{equation}
and is measuring the defect in the geometric-logarithmic mean
inequality \eqref{logmeangeologar}. Let us switch to exponential
variables and set
\[
x(t) := \log\sqrt{\frac{t}{p}} \quad\mbox{and}\quad y(t) := \log \sqrt{
\frac{1-t}{1-p}} .
\]
Note that either $x(t)\leq0 \leq y(t)$ for $t\leq p$ or $y(t)\leq0
\leq x(t)$ for $t\geq0$ with equality only for $t=p$. Therewith,
\eqref{logmeanlemoptp8} can be rewritten as
\[
\tilde w_p(t) = \frac{x(t)}{\sinh (x(t) )}-\frac
{y(t)}{\sinh
(y(t) )} .
\]
By making use of the fact, that the function $x\mapsto\frac{x}{\sinh
x}$ is symmetric, strictly monotone decreasing in $\llvert  x\rrvert $ and has a
unique maximum in $1$, we can conclude that
\[
\tilde w_p(t)=0 \quad\mbox{if and only if} \quad x(t) = - y(t) .
\]
The solutions to the equation $x(t)=-y(t)$ are given for $t\in \{
p,1-p \}$. Let us first consider the case $t=p$, then
$x(t)=y(t)=0$ and
$w_p(p)$ is a zero of order $2$, since the function $x\mapsto\frac
{x}{\sinh(x)}$ is strictly concave for $t=0$. Now, we can go back to
$h_p'(t)$ and argue with the representation
\[
\lim_{t\to p} h_p'(t) = \lim
_{t\to p} \frac{ v_p(t) \hat w_p(t)
\tilde w_p(t)}{g_p^2(t)} \overset{!} {\ne} 0 .
\]
This is a consequence of counting the zeros for $t=p$ in the nominator
and denominator according to their order; for the denominator
$g_p^2(p)$ is a zero of order $4$. For the nominator, we have $v_p(p)$
is a zero of order $1$, $\hat w_p(p)$ is a zero of order $1$ and
$\tilde w_p(p)$ is a zero of order $2$, which leads in total again to a
zero of order $4$ exactly compensating the zero of the denominator.

The other case is $t=1-p$. Let us evaluate $h_p(1-p)$, which is given by
\begin{eqnarray}
h_p(1-p) &=& \frac{\fracc{(p-(1-p)  )^2}{p(1-p)}  }{(1-p)
\log(\fracb{1-p}{p}) + p \log(\fracc{p}{1-p})}
\nonumber
\\
&=& \frac{1}{p(1-p)} \frac{(p-(1-p))^2}{(p-(1-p)) \log(\fracc
{p}{1-p})} = \frac{\Lambda(p,1-p)}{p(1-p)} .\nonumber
\end{eqnarray}
Since $t=1-p$ is the only critical point of $h_p(t)$ inside $(0,1)$, it
remains to check whether the boundary values are larger than
$h_p(1-p)$. They are given by
\[
\lim_{t\to0} h_p(t) = \frac{1}{(1-p)\log(\fracc{1}{1-p})} \quad
\mbox{and}\quad\lim_{t\to1} h_p(t) =
\frac{1}{p\log(\fraca{1}{p})}.
\]
We observe that the demanded inequality to be in a global minimum
\[
h_p(1-p) = \frac{\Lambda(p,1-p)}{p(1-p)} \overset{!} {<} \min \biggl\{
\frac{1}{p\log(\fraca{1}{p})},\frac{1}{(1-p)\log(\fracc{1}{1-p})} \biggr\}
\]
is just \eqref{logmeanlemeupperbound} of Lemma~\ref
{logmeanlemupperbound}.
\end{pf}

\section{Integration by parts on basins of attraction}\label{cgradsys}
The goal of this Appendix is to proof the integration by parts formula,
which is an ingredient of the Lyapunov approach in Section~\ref{clocal}.
%
\begin{theo}[(Integration by parts)]\label{thmNeumann}
Let $H\in C^3(\R^n,\R)$ be a Morse function (cf. Definition~\ref
{defmorse}) with compact sublevel sets and let $\Omega$ be the basin
of attraction associated to a local minimum of $H$ (cf. Definition~\ref
{defStabMfd}), then it holds
\[
\forall f,g\in H^1(\mu|_\Omega)\qquad \mbox{with } \nabla g
\parallel \nabla H \mbox{ on } \partial\Omega\dvtx  \int_\Omega f
(-L g) \,\dx{\mu } = \eps\int_\Omega \langle\nabla f , \nabla
g \rangle \,\dx{\mu} ,
\]
where $\nabla g \parallel \nabla H$ means $\llvert \nabla g(x) \cdot
\nabla H\rrvert  =
\llvert  \nabla g(x) \rrvert    \llvert \nabla H(x)\rrvert $ for $\cH^{n-1}$-a.e. $x\in\partial
\Omega$.
\end{theo}
%
%
\begin{rem}
The property of $H$ possessing compact sublevel sets is called proper.
This gives enough compactness, that is, the Palais--Smale
condition \cite{Jost2008}, Definition~6.2.1, to apply several results
from Morse theory and dynamical systems. Moreover, if $H$ satisfies
Assumption \eqref{assumeenv}, then $H$ is proper.
\end{rem}

\subsection{Properties of gradient flows}
%
\begin{defn}[(Gradient flow)]
Let $\phi_t(x)$ be the trajectory associated to the negative gradient
flow of $H$ started in $x$, that is,
\[
\partial_t \phi_t = - \nabla H(\phi_t)
\quad\mbox{and}\quad\phi _0(x)=x \in\R^n.
\]
\end{defn}
%
%
\begin{lem}[(Properties of gradient flow trajectories)]\label
{lempropGF}
\begin{longlist}[(iii)]
\item For each $x$, the trajectory $t\mapsto\phi_t(x)$ has a maximal
interval of definition of the form $(-\alpha_x, \infty)$ for $\alpha
_x \in(-\infty,0)\cup \{-\infty \}$.
\item For each $x$:  $\lim_{t\to\infty} \phi_t(x) =: \phi_\infty
(x)\in\cS$.
\item Stability on finite time intervals, that is, for any $T>0$ holds
if $x_n\to x$ also $\phi_T(x_n)\to\phi_T(x)$.
\end{longlist}
\end{lem}
\begin{pf}
Since $H$ is locally Lipschitz, the trajectory $\phi_t(x)$ has a
maximal interval of definition $0\in(\alpha_x,\beta_x)\cup \{
\pm \infty \}$ by the Picard--Lindel\"{o}f theorem. Moreover, since
\begin{equation}
\label{egradsyspartialH} \partial_t H(x_t)= - \bigl\llvert \nabla
H(x_t)\bigr\rrvert ^2 = - \llvert \dot x_t
\rrvert ^2 \leq0
\end{equation}
the trajectory $ \{\phi_t(x) \}_{t\geq0}$ is confined to the
sublevel set $ \{y \dvtx  H(y)\leq H(x) \}$, which is compact,
since $H$ is
proper. On this sublevel set, $H$ is globally Lipschitz and the limit
$\lim_{t\to\infty} \phi_t(x) =: \phi_\infty(x)$ exists
proving (i). In addition, this implies
\[
\int_0^\infty\bigl\llvert \nabla H(
\phi_t)\bigr\rrvert ^2 \,\dx{t} \stackrel {\mathclap{
\eqref{egradsyspartialH}}} {=} - \int_0^\infty
\partial _t H(\phi_t) \,\dx{t} = H(x) - H\bigl(
\phi_\infty(x)\bigr) < \infty.
\]
Therefore, it follows $\phi_\infty(x) \in\cS:= \{x\in\R^n \dvtx
\nabla H(x)=0 \}$ is a critical point proving (ii). The stability
follows from the estimate
\begin{eqnarray}
\label{eqstabest} \bigl\llvert \phi_T(x_n)-
\phi_T(x)\bigr\rrvert &=& \biggl\llvert x_n + \int
_0^T \partial_t
\phi_t(x_n) \,\dx{t} - x - \int_0^T
\partial _t \phi_t(x)\biggr\rrvert
\nonumber
\\[-8pt]
\\[-8pt]
&\leq&\llvert x_n-x\rrvert + \int_0^T
\bigl\llvert \nabla H\bigl(\phi_t(x_n)\bigr)-\nabla H
\bigl(\phi_t(x)\bigr)\bigr\rrvert \,\dx{t}.
\nonumber
\end{eqnarray}
All $\phi_t(x_n)$ are confined to a common compact set by properness
of $H$ and in particular $\nabla H$ is Lipschitz continuous in this
compact set. This leads for some $K>0$ and all $t\in(0,T)$ to the estimate
\[
\bigl\llvert \nabla H\bigl(\phi_t(x_n)\bigr)-\nabla H
\bigl(\phi_t(x)\bigr)\bigr\rrvert \leq K \bigl\llvert \phi
_t(x_n)-\phi_t(x)\bigr\rrvert .
\]
Using this estimate in \eqref{eqstabest}, we can apply the Gronwall
inequality to obtain $\llvert \phi_T(x_n)-\phi_T(x)\rrvert
\leq\llvert  x_n-x\rrvert (1+ e^{KT})$, which proves (iii).
\end{pf}
We want to define a global flow w.r.t. $\nabla H$. Since, $\nabla H$
can have superlinear growth and is in particular not globally Lipschitz
continuous, we use the following reparameterized version for a global flow.
%
\begin{theo}[(Global flow by reparameterization \cite{Weiss2007}, Theorem~4.4)]\label{thmGFrepara}
A glob\-al flow of diffeomorphism $\tilde\phi_t\dvtx \R^n \to\R^n$ w.r.t.
$H$ is defined by
\begin{equation}
\label{edefmodflow} \partial_t \tilde\phi_t(x) = F\bigl(
\tilde\phi_t(x)\bigr) := -\frac{\nabla
H(\tilde\phi_t(x))}{1+\llvert \nabla H(\tilde\phi_t(x))\rrvert } \quad \mbox{and}\quad
\tilde\phi_0(x) =x .
\end{equation}
This flow is equivalent to the negative gradient flow of $H$ upon a
reparameterization of time. The vector field $F$ is globally Lipschitz
and bounded. It defines a global flow on $\R^n$, that is, $\tilde\phi
_{t+s} = \tilde\phi_t \circ\tilde\phi_s$ for all $t,s\in\R$.
\end{theo}
%
%
\begin{cor}\label{coruniqueGF}
Each point $x\in\R^n$ belongs to exactly one trajectory $t\to\phi_t(x)$.
\end{cor}
\begin{pf}
We apply \cite{Jost2008}, Corollary~1.9.1, to the global flow $\tilde
\phi_t$ and by Theorem~\ref{thmGFrepara} translate the result back
to $\phi_t$.
\end{pf}
%
\subsection{The stable manifold}
%
\begin{defn}[(Stable manifold)]\label{defStabMfd}
To each critical point $s\in\cS$, the \emph{stable manifold} is
defined by
\[
W^\tts(s) := \Bigl\{x \in\R^n\dvtx  \lim_{t\to\infty}
\phi_t(x) = s \Bigr\} .
\]
Moreover, we call the dimension $k\in \{0,\dots,n \}$ of the
unstable subspace of $\nabla^2 H(s)$ the index of the saddle point
$s$. If $m$ is a local minimum of $H$, that is, a critical point of
index $0$, we call $W^\tts(m)$ the \emph{basin of attraction for $m$}.
\end{defn}
Lemma~\ref{lempropGF}(ii) and Corollary~\ref{coruniqueGF} ensure
the stable manifold to be well defined and immediately provide the following.
%
\begin{cor}[(Partition of state space)]\label{corpartStateSpace}
Let $\cS$ be all critical points of $H$, then $\R^n$ is the disjoint
union of all stable manifolds denoted by
\[
\label{coreqpartStateSpace} \R^n := \bigcupdot_{s\in\cS}
W^\tts(s) .
\]
\end{cor}
%
%
\begin{theo}[(Local stable manifold theorem \cite{Jost2008}, Theorem~6.3.1)]\label{thmLSMTf}
Let $s\in\cS$ and $\cE^\tts(s)$ be the stable subspace of $\nabla
^2 H(s)$, that is, $\nabla^2 H(s)|_{\cE^\tts}$ has a positive
spectrum. Then there exists a neighborhoods $U, \tilde U$ of $s$, such
that $W^\tts(s)\cap U$ is a $C^1$-graph over $(s+\cE^\tts(s)) \cap
\tilde U$. Especially, the dimension of $W^\tts(s)\cap U$ and $\cE
^\tts(s)$ are equal to $n-k$, where $k$ is the index of $s$.
\end{theo}
The local result can be extended by the reparameterized flow to the
global manifold theorem.
%
\begin{theo}[(Global stable manifold theorem \cite{Jost2008}, Corollary~6.3.1)]\label{thmGSMTF}
The stable manifolds $W^\tts(s)$ for $s\in\cS$ of the flow
associated to $F$ \eqref{edefmodflow} are immersed $C^1$-manifolds
of dimension $n-k$, where $k$ is the index of $s$.
\end{theo}
In the present case of a gradient flow, the result can be strengthened
to the following.%
\begin{theo}[(Global stable manifold theorem for gradient systems \cite{Jost2008}, Corollary~6.4.1)]\label{thmGSMTGF}
The stable manifolds $W^\tts(s)$ for $s\in\cS$ of the gradient flow
associated to $H$ are embedded $C^1$-submanifolds of dimension $n-k$,
where $k$ is the index of $s$.
\end{theo}
\begin{pf}
We have to modify the proof of \cite{Jost2008}, Corollary~6.3.1, since
$\nabla H$ can have superlinear growth. Instead, considering the
gradient flow w.r.t. $H$, we consider the equivalent flow $\tilde\phi
_t$ of Theorem~\ref{thmGFrepara}. We have to observe two additional
facts, which we postpone to the end of the proof.
\begin{longlist}[(b)]
\item[(a)] The flow has no nonconstant homoclinic orbits, that is,
nonconstant orbits with $\lim_{t\to-\infty} \tilde\phi_t(x) = \lim_{t\to\infty}\tilde\phi_t(x)$ (cp. \cite{Jost2008}, Lemma~6.4.3).
\item[(b)] For each $x$, holds $\llvert \nabla H(\tilde\phi_t(x))\rrvert  \to0$ as
$t\to\infty$ and either $\llvert \nabla H(\tilde\phi_t(x))\rrvert  \to0$ or
$H(\tilde\phi_t(x)) \to\infty$ as $t\to-\infty$ (cp. \cite{Jost2008}, Lemma~6.4.4).
\end{longlist}
This allows us to complete the proof by first applying Theorem~\ref
{thmGSMTF} to $F(x) = -\nabla H(x) /(1+\llvert \nabla H(x)\rrvert )$. Every point
$x\in\R^n$ is contained in a unique trajectory $\phi_t(x)$ by
Corollary~\ref{coruniqueGF}. However, a trajectory is typical not
compact. In (b) we show that limit points in $\R^n$ are critical points
of $H$. The local situation around critical points is given by the
local stable manifold theorem~\ref{thmLSMTf}, which provides a local
chart around the critical point. Selfintersection of trajectory is
excluded by the observation in (a). Hence, the immersion of Theorem~\ref
{thmGSMTF} is an embedding.

We still have to show (a) and (b):

Ad (a): The energy also decreases w.r.t. to the reparameterized flow
\begin{equation}
\label{eqestmodflow} \partial_t H\bigl(\tilde\phi_t(x)\bigr) =
- \nabla H \cdot\partial_t \tilde \phi_t(x) = -
\frac{\llvert \nabla H(\tilde\phi_t(x))\rrvert ^2}{1+ \llvert \nabla H(\tilde\phi_t(x))\rrvert } \leq 0 .
\end{equation}
Hence, for a trajectory either holds $\llvert \nabla H\rrvert =0$ or $\llvert \nabla H\rrvert >0$
for all $t$, which gives (a).

Ad (b): Integrating the identity \eqref{eqestmodflow}, we obtain for
$t_2 > t_1$
\[
\label{eqestenergymodflow} H\bigl(\tilde\phi_{t_1}(x)\bigr) - H\bigl(\tilde
\phi_{t_2}(x)\bigr) = \int_{t_1}^{t_2}
\frac{\llvert \nabla H(\tilde\phi_t(x))\rrvert ^2}{1+ \llvert \nabla H(\tilde \phi_t(x))\rrvert } \,\dx{t} \geq\int_{t_1}^{t_2}
\bigl\llvert \nabla H\bigl(\tilde \phi_t(x)\bigr)\bigr\rrvert \,\dx{t}
.
\]
Since $H$ is bounded from below, we get that $H(\tilde\phi_{\infty
}(x))>-\infty$. Hence,
\[
H\bigl(\tilde\phi_{t_1}(x)\bigr)-H\bigl(\tilde\phi
_{t_2}(x)\bigr)<\infty
\]
 for all $t_2>t_1$ and we immediately deduce
from \eqref{eqestmodflow} that $\tilde\phi_\infty(x)\in\cS$
showing the first part of (b). If $H(\tilde\phi_{-\infty}(x))<\infty
$, then by the same argument $\phi_{-\infty}(x)\in\cS$. Hence, we
have shown the dichotomy (b).
\end{pf}

\subsection{The boundary of the basin of attraction}

%
\begin{lem}\label{lempartLebesgueGF}
The set $ \{W^\tts(m) \}_{m\in\cM}$ is a partition of $\R
^n$ upon
Lebesgue null sets, denoted by
\begin{equation}
\label{epartLebesgueGF} \R^n = \biguplus_{m\in\cM} W^{\tts}(m).
\end{equation}
Moreover, it holds
\begin{equation}
\label{epartBoundaryGF} \bigcup_{m\in\cM} \partial
W^{\tts}(m) = \bigcupdot_{y\in\cS
\setminus\cM} W^{\tts}(y).
\end{equation}
\end{lem}
\begin{pf}
For \eqref{epartLebesgueGF}, we observe that $W^\tts(y)$ for $y\in
\cS\setminus\cM$ are Lebesgue null sets, since they are $(n-k)$-dimensional $C^1$-submanifolds for $1\leq k \leq n$ (cf. Theorem~\ref
{thmGSMTGF}).

 Theorem~\ref{thmGSMTGF} proves in
particular, that for each $m\in\cM$ the embedded submanifold $W^\tts
(m)$ is open in $\R^n$, hence $\partial W^\tts(m) \cap W^\tts(m)
=\varnothing$. Therewith, the second statement \eqref
{epartBoundaryGF} follows from Corollary~\ref{corpartStateSpace}.
\end{pf}

%
\begin{theo}[(The boundary of the basin of attraction)]\label{thmBoundStabMfd}
Let $m\in\cM$ be a local minimum of $H$. There exists a set
$S_m\subset\cS\setminus\cM$ of $k$-saddles with $k\geq1$ such that
\[
\partial W^\tts(m) = \bigcupdot_{y\in S_m} W^\tts(y).
\]
\end{theo}
\begin{pf}
We define a critical point $y\in\cS$ to be in $S_m$ if for each open
neighborhood $U(y)$ holds $U(y)\cap W^\tts(m) \ne\varnothing$.
From~\ref{epartBoundaryGF} follows that $y\in S_m$ cannot be another
local minimum and hence $S_m \subset\cS\setminus\cM$. Now, we take
$x_n \to x \in\partial W^\tts(m)$. From \eqref{epartBoundaryGF}
follows that $x\in W^\tts(y)$ for some $y\in\cS\setminus\cM$. We
have to prove that $y\in S_m$. There exists an open neighborhood $U(x)$
such that $x_n \in U(x)$ for $n>N$. Then for any open neighborhood
$U(y)$ of $y$ exists $T>0$ such that $\phi_T(x) \in U(y)$. By
existence of the flow $\phi_t$ for positive time, it follows hat $\phi
_T(U(x)) \cap U(y)=:U(\phi_T(x))$ is an open neighborhood of $\phi
_T(x)$. By stability of the flow on finite time intervals [cf.
Lemma~\ref{lempropGF}(iii)], it follows $ \phi_T(x_n)\to\phi
_T(x)$, hence $\phi_T(x_n)\in U(\phi_T(x))$ for $n$ large enough,
which shows that $W^\tts(m)\cap U(\phi_T(x))\ne\varnothing$ and
finally $y\in S_m$.
\end{pf}
\begin{pf*}{Proof of Theorem~\ref{thmNeumann}}
Let $m$ be a local minimum of $H$. By Theorem~\ref{thmBoundStabMfd}
the boundary of $W^\tts(m)$ is the union of $C^1$-submanifolds. The
relevant submanifolds for integration, are the $(n-1)$-dimensional ones.
By Theorem~\ref{thmGSMTGF}, these $(n-1)$-dimensional submanifolds
correspond to stable manifolds of saddle points of index $1$, denoted
by $\cS_1$. Hence, for $\cH^{n-1}$-a.e. $x\in\partial W^\tts(m)$
exists a $1$-saddle $y\in S_m\cap\cS_1$ such that $x\in W^\tts(y)$.
Therefore, the normal on $W^\tts(m)$ exists $\cH^{n-1}$-a.e., which
gives enough regularity to integrate for $f,g\in H^1(\mu|_\Omega)$ by parts
\[
\int_\Omega f (-L g) \,\dx{\mu} = \eps\int
_{W^\tts(m)} \langle\nabla f , \nabla g \rangle \,\dx{\mu} - \eps
\sum_{y\in S_m\cap\cS
_1}\int_{W^\tts(y)} f \nabla g
\cdot n \cH^{n-1}(\dx{\mu}) .
\]
%
By the assumption $\nabla g \parallel \nabla H$, it is enough to show that
$\nabla H(x) \cdot n = 0$ for $\cH^{n-1}$-a.e. $x\in\partial W^\tts
(m)$. This is proven by contradiction for $x\in\partial W^\tts(m)$.
Assume that $x\notin\cS$, that is, $\nabla H(x) \ne0$ and $\nabla
H(x) \cdot n \ne0$. Then for some $\eps>0$ there exists $t^\ast\in
(-\eps,\eps)$ such that $\phi_{t^\ast}(x) \in W^\tts(m)$. By
definition of $W^\tts(m)$ and global existence of the trajectory
$ \{\phi_t(x) \}_{t\geq t^\ast}$ from Lemma~\ref
{lempropGF}(ii) follows
$x\in W^\tts(m)$, which contradicts \eqref{epartBoundaryGF} and
Corollary~\ref{corpartStateSpace}.
\end{pf*}

\section{Auxiliary results from Section~\texorpdfstring{\protect\ref{CMEANDIFF}}{4}}

\subsection{Partial Gaussian integrals}\label{appendsubecpartGauss}
This section is devoted to proof the representation for partial or
incomplete Gaussian integrals. Lemma \eqref{appendlempartGauss} is
an ingredient to evaluate the weighted transport cost in Section~\ref{safftransp}.

%
\begin{lem}[(Partial Gaussian integral)]\label{appendlempartGauss}
Let $\Sigma^{-1}\in\R^{n\times n}_{\sym,+}$ be a symmetric positive
definite matrix and let $\eta\in S^{n-1}$ be a unit vector. Therewith,
$ \{r\eta+z^* \}_{r\in\R}$ is for $z^* \in\R^n$ with
$ \langle\eta , z^*  \rangle=0$ an affine subspace of $\R
^n$. The integral of a centered
Gaussian w.r.t. to this subspace evaluates to
\begin{eqnarray}
\label{appendequpartGauss} \int_{\R} \exp \biggl(-\frac{1}{2}
\Sigma^{-1}\bigl[r\eta+ z^*\bigr] \biggr)\,\dx {r} = \frac{\sqrt{2\pi}}{\sqrt{\Sigma^{-1}[\eta]}}
\exp { \bigl(-{\tilde\Sigma^{-1}\bigl[z^*\bigr]} \bigr)},
\nonumber
\\
\eqntext{\displaystyle \mbox{with } \tilde\Sigma^{-1} = \Sigma^{-1} -
\frac
{\Sigma^{-1}\eta\otimes\Sigma^{-1}\eta}{\Sigma^{-1}[\eta]} .} %
\end{eqnarray}
\end{lem}
\begin{pf}
To evaluate this integral on an one-dimensional subspace of $\R^n$, we
have to expand the quadratic form $\Sigma^{-1}[r\eta+ z^*]$ and
arrive at the relation
\begin{eqnarray*}
&&\int_{\R} \exp \biggl(-\frac{1}{2}\Sigma
^{-1}\bigl[r\eta+ z^*\bigr] \biggr)\,\dx{r}
\\
&&\qquad= \exp{ \biggl(-\frac{1}{2}\Sigma^{-1}\bigl[z^*\bigr]
\biggr)} \int_{\R
}\exp \biggl(-\frac{r^2}{2}
\Sigma^{-1}[\eta] + r \bigl\langle\eta , \Sigma^{-1}z^* \bigr
\rangle \biggr) \,\dx{r}
\\
&&\qquad= \exp{ \biggl(-\frac{1}{2}\Sigma^{-1}\bigl[z^*\bigr]
\biggr)} \frac
{\sqrt{2\pi}}{\sqrt{\Sigma^{-1}[\eta]}} \exp \biggl(\frac {
\langle\eta , \Sigma^{-1}z^*  \rangle^2}{2\Sigma^{-1}[\eta
]} \biggr)
\\
&&\qquad= \frac{\sqrt{2\pi}}{\sqrt{\Sigma^{-1}[\eta]}} \exp { \biggl(-\frac{1}{2}{ \biggl(
\Sigma^{-1}-\frac{\Sigma ^{-1}\eta
\otimes\Sigma^{-1}\eta}{\Sigma^{-1}[\eta]} \biggr)\bigl[z^*\bigr]} \biggr)} ,
\end{eqnarray*}
which concludes the hypothesis.
\end{pf}

\subsection{Subdeterminants, adjugates and inverses}\label
{appendsubecsubdet}

Let $A\in\R^{n\times n}_{\sym,+}$, then define for $\eta\in
S^{n-1}$ the matrix
\begin{equation}
\label{appendtildeAdef} \tilde A := A - \frac{A\eta\otimes A\eta}{A[\eta]} .
\end{equation}
The matrix $\tilde A$ has at least rank $n-1$, since we subtracted from
the positive definite matrix $A$ a rank-$1$ matrix. Further, from the
representation, it is immediate that $\tilde A$ has rank $n-1$ if and
only if $\eta$ is an eigenvector of $A$. In this case, $\ker A=\vspan
{\eta}$. It immediately follows $\tilde A > 0$ on $V:=\vspan \{
\eta \}^\perp$, which is the $(n-1)$-dimensional subspace
perpendicular to $\eta$. Then for a matrix $A\in\R^{n\times n}_{\sym
,+}$ we want to calculate the determinant of $A$ restricted to this
subspace $V$. This determinant is obtained by first choosing $Q\in
SO^n$ such that $Q( \{0 \}\times\R^{n-1}) = V$ and then evaluating
the determinant of the minor consisting of the $(n-1)\times(n-1)$
lower right submatrix of $Q^\top A Q$ denoted by $\det_{1,1}(Q^\top
AQ)$. Hence, we have
\[
\label{appenddefnrestrictdet} \det_{1,1}\bigl(Q^\top A Q\bigr)  \qquad
\mbox{with } Q\in SO(n)\dvtx  Q^\top\eta= e^1 = (1,0,
\dots,0)^\top.
\]
Since $V=\vspan \{\eta \}^\perp$, it follows that the
first column
of $Q$ is given by $\eta$ and we can decompose $Q^\top A Q$ into
\[
Q^\top A Q = %
\pmatrix{ A[\eta] & \widehat{Q^\top
A\eta}
\cr
\widehat{Q^\top A\eta}{}^\top& \widehat{Q^\top
AQ} } %
,
\]
where for a matrix $M$, $\widehat{M}$ is the lower right $(n-1)\times
(n-1)$ submatrix of $M$ and for a vector $v$, $\widehat{v}$ the
$(n-1)$ lower subvector of $v$. Therewith, we find a similarity
transformation which applied to $Q^\top AQ$ results in
\begin{eqnarray*}
\det A &=& \det Q^\top AQ = \det  \left(\pmatrix{ A[\eta] &
\widehat{Q^\top A\eta}
\cr
\widehat{Q^\top A
\eta}{}^\top & \widehat{Q^\top AQ} } \pmatrix{   1 & -
\dfrac{\widehat{Q^\top A\eta
}}{A[\eta]}
\cr
0 & \Id_{n-1} }  \right)
\\
&=& \det %
\pmatrix{ A[\eta] & 0
\cr
\widehat{Q^\top A
\eta}{}^\top& \widehat{Q^\top AQ} - \dfrac{\widehat
{A\eta}\otimes\widehat{A\eta}}{A[\eta]} }
\\
&=& A[\eta] \det_{1,1} \biggl({Q^\top AQ} -
\frac{{Q^\top A\eta
}\otimes{Q^\top A\eta}}{A[\eta]} \biggr) .
\end{eqnarray*}
The determinant of the minor is given by
\[
\det_{1,1} \biggl({Q^\top AQ} - \frac{{Q^\top A\eta}\otimes {Q^\top
A\eta}}{A[\eta]} \biggr)
= \det_{1,1} \biggl(Q^\top \biggl(A - \frac
{A\eta\otimes A \eta}{A[\eta]}
\biggr)Q \biggr) .
\]
Hence, by the definition \eqref{appendtildeAdef} of $\tilde A$ and
the subdeterminant, we found the identity
\begin{equation}
\label{appendtildeAdetident} \det A = A[\eta] \det_{1,1}\bigl(Q^\top\tilde
A Q\bigr) .
\end{equation}

\subsection{A matrix optimization}

%
\begin{lem}\label{coarselemmatrixopt}
Let $B\in\R^{n \times n}_{\sym,+}$, then it holds
\[
\inf_{A\in\R^{n \times n}_{\sym,+}} \biggl\{\frac{\det A}{\sqrt {\det (2 A-B )}} \dvtx  2A>B \biggr\} =
\sqrt{\det B}
\]
and for the optimal $A$ holds $A=B$.
\end{lem}
\begin{pf}
We note that
\[
\frac{\det A}{\sqrt{\det (2 A-B )}}= \frac{1}{\sqrt {\det(
A^{-1})\det(2\Id- A^{-\fraca{1}{2}} B A^{-\fraca{1}{2}})}} .
\]
Therewith, it is enough to maximize the radical of the root. Therefore,
we substitute $ A^{-\fraca{1}{2}}=C B^{-\fraca{1}{2}}$ with $C>0$ not
necessarily symmetric and observe that $ A^{-\fraca{1}{2}}= B^{-\fraca
{1}{2}}C^\top$. We obtain
\[
\det\bigl( A^{-1}\bigr)\det\bigl(2\Id- A^{-\fraca{1}{2}} B
A^{-\fraca{1}{2}}\bigr) = \det\bigl( B^{-1}\bigr)\det
\bigl(CC^\top\bigr)\det\bigl(2\Id- CC^\top\bigr) .
\]
Note that $CC^\top\in\R^{n\times n}_{\sym,+}$ and it is enough to calculate
\[
\sup_{\tilde C\in\R^{n\times n}_{\sym,+}} \bigl\{\det(\tilde C)\det (2\Id- \tilde C) \dvtx {
\tilde C<2\Id} \bigr\} .
\]
From the constraint $0<\tilde C<2 \Id$, we can write $\tilde C=\Id
+D$, where $D$ is symmetric and satisfies $-\Id< D<\Id$ in the sense
of quadratic forms. From here, we finally observe
\[
\det(\tilde C)\det(2\Id- \tilde C) = \det(\Id+D)\det(\Id-D) = \det\bigl(\Id-
D^2\bigr) .
\]
Since $D^2\geq0$, we find the optimal $\tilde C$ given by $\Id$,
which yields that $ A= B$.
\end{pf}

%

\subsection{Jacobi matrices}\label{secappJacMat}

For a smooth function $f \dvtx  \R^n \to\R^n$ denotes $Df(x)$ the \emph
{Jacobi matrix} of the partial derivatives of $f$ in $x\in\R^n$ given by
\[
Df(x) := \biggl({\frac{\textup{d} f_i}{\textup{d} x_j}}(x) \biggr)_{i,j=1}^n .
\]
%
%
\begin{lem}
Let $A,B\in\R^{n\times n}$ and $f\dvtx \R^n\to\R^n$ smooth, then it holds
\begin{eqnarray}
\nabla\bigl\llvert Ax + f(B x)\bigr\rrvert &=& \bigl(A + Df(x) B
\bigr)^\top \frac{Ax +
f(Bx)}{\llvert  Ax+f(Bx)\rrvert }, \label{matrixjacf1}
\\
D \frac{f(x)}{\llvert  f(x)\rrvert } &=& \frac{1}{\llvert
f(x)\rrvert } \biggl(\Id- \frac{f(x)}{\llvert  f(x)\rrvert }
\otimes\frac{f(x)}{\llvert  f(x)\rrvert } \biggr) Df(x) \label{matrixjacf2}.
\end{eqnarray}
\end{lem}
\begin{pf}
Let us first check the relation \eqref{matrixjacf1} and calculate
the partial derivative
\begin{equation}
\label{matrixjacf1p1} %
\quad {\frac{\textup{d} \llvert  Ax + f(B x)\rrvert }{\textup{d}
x_i}} = \frac{1}{2\llvert  Ax + f(B x)\rrvert } \sum
_{j} {\frac{\textup{d} }{\textup{d} x_i}} \biggl(\sum
_{k} A_{jk}x_k +
f_j(Bx) \biggr)^2. %
\end{equation}
The inner derivative of \eqref{matrixjacf1p1} evaluates to
\begin{eqnarray}
\label{matrixjacf1p2} &&{\frac{\textup{d} }{\textup{d} x_i}} \biggl(\sum_{k}
A_{jk}x_k + f_j(Bx) \biggr)^2\nonumber
\\[-8pt]
\\[-8pt] &&\qquad =
2 \biggl(\sum_{k} A_{jk}x_k
+ f_j(Bx) \biggr) \biggl(A_{ji} + {\frac
{\textup{d} f_j(Bx)}{\textup{d} x_i}}
\biggr) .
\nonumber
\end{eqnarray}
The derivative of $f_j(Bx)$ becomes
\begin{eqnarray}
\label{matrixjacf1p3} {\frac{\textup{d} f_j(Bx)}{\textup{d} x_i}} &=& {\frac{\textup{d}
f_j (\sum_{k} B_{1k} x_k, \dots,\sum_{k} B_{nk} x_k
)}{\textup{d} x_i}}
\nonumber
\\[-8pt]
\\[-8pt]
&=& \sum_{k=1}^n \partial_k
f_j(Bx) B_{ki} = \bigl(Df(Bx) B\bigr)_{ji}.
\nonumber
\end{eqnarray}
Hence, a combination of \eqref{matrixjacf1p1}, \eqref
{matrixjacf1p2} and \eqref{matrixjacf1p3} leads to
\begin{eqnarray}
\label{matrixjacf1p4} %
{\frac{\textup{d} \llvert  Ax + f(B x)\rrvert }{\textup{d}
x_i}} &=& \frac{1}{\llvert  Ax + f(B x)\rrvert } \sum
_{j} \bigl((Ax)_j +
f_j(Bx) \bigr) \bigl(A_{ji} \bigl(Df(Bx)B
\bigr)_{ji} \bigr)
\nonumber
\\
&=& \sum_{j} \bigl(A + Df(Bx)B
\bigr)^\top_{ij} \frac{(Ax + f(Bx))_j}{\llvert  Ax + f(B x)\rrvert } ,\nonumber %
\end{eqnarray}
which shows \eqref{matrixjacf1}. For the equation \eqref
{matrixjacf2}, let us first consider the Jacobian of the function
$F(x)= \frac{x}{\llvert  x\rrvert }$, which is given by
\[
D F(x) = \frac{1}{\llvert  x\rrvert } \biggl(\Id- \frac
{x}{\llvert  x\rrvert }\otimes
\frac{x}{\llvert  x\rrvert } \biggr) .
\]
Then, by the chain rule, we observe that
\[
D \frac{f(x)}{\llvert  f(x)\rrvert } = D (F\circ f) (x) = DF\bigl(f(x)\bigr) Df(x),
\]
which is just \eqref{matrixjacf2}.
\end{pf}
\end{appendix}

\section*{Acknowledgments}
This work is part of the Ph.D. thesis of the second author, written
under the supervision of Stephan Luckhaus at the University of Leipzig.

The authors are greatly indebted to Felix Otto for drawing their
attention to the low temperature regime and for several helpful
discussions.

Moreover, they want to thank the anonymous referees for pointing out
the role of the Neumann boundary conditions in the Lyapunov argument.

They wish to thank the Max-Planck-Institute for Mathematics in the
Sciences in Leipzig where the paper was written under great working
conditions. 


%

\printaddresses

\end{document}